%% file: statSolsNSE.tex
\documentclass[11pt,a4paper]{elsarticle}

\usepackage[a4paper, margin=1.1in]{geometry}
\usepackage{amsmath,amsfonts}
\usepackage{hyperref}
\usepackage{amssymb}
\usepackage{amsthm}
\usepackage{mathrsfs}
\usepackage{mathtools}
\usepackage{esint}
\usepackage{graphicx, subcaption}
\usepackage{dsfont}
\usepackage{url}
\usepackage{textcomp}
\usepackage{braket}
\usepackage{eufrak}
\usepackage[shortlabels]{enumitem}
\usepackage[ruled]{algorithm2e}

\input{macros.tex}

\newcommand{\figsDir}{}

\makeatletter
\def\ps@pprintTitle{%
 \let\@oddhead\@empty
 \let\@evenhead\@empty
 \def\@oddfoot{}%
 \let\@evenfoot\@oddfoot}
\makeatother

\begin{document}

\begin{abstract}
Statistical solutions, which are time-parameterized probability measures on spaces of square-integrable functions, have been established as a suitable framework for global solutions of incompressible Navier-Stokes equations (NSE).
We compute numerical approximations of statistical solutions of NSE on two-dimensional domains with \emph{non-periodic} boundary conditions and empirically investigate the convergence of these approximations and their observables. For the numerical solver, we use Monte Carlo sampling with an H(div)-FEM based deterministic solver.
Our numerical experiments for high Reynolds number turbulent flows demonstrate that the statistics and observables of the approximations converge. We also develop a \emph{novel} algorithm to compute \emph{structure functions} on unstructured meshes.
\end{abstract}

\begin{frontmatter}
\title{Numerical approximation of statistical solutions of the incompressible Navier-Stokes Equations}

\author[eth]{Pratyuksh Bansal\corref{cor1}}
\cortext[cor1]{Corresponding author.}
\ead{pratyuksh.bansal@sam.math.ethz.ch}
\address[eth]{Seminar for Applied Mathematics, ETH Z\"urich, R\"amistrasse 101, CH-8092 Z\"urich}


\end{frontmatter}
%
%
\section{Introduction} 
%
The flow of a viscous, incompressible Newtonian fluid is described by the (incompressible) Navier-Stokes equations:
\begin{align*}
 \der{\bu}{t} + \div_{\bx}(\bu \otimes \bu) - \nu \Delta_{\bx} \bu + \nabla_{\bx}{p} = \bb{f},\quad
  \div_{\bx}{\bu} = 0.
\end{align*}
Here $\bu$ is the velocity, $p$ is the pressure, which plays the role of a Lagrange multiplier to enforce the divergence-free constraint on the velocity,
$\nu$ is the kinematic viscosity and $\bb{f}$ represents the effects of external forces, for example, gravity and buoyancy. This system also needs to be supplied with initial and boundary conditions. 

Navier-Stokes equations play a fundamental role in many applications and have been studied extensively for more than a century.
The existence of global-in-time weak solutions, both in two and three spatial dimensions, can be traced back to the pioneering works of Leray~\cite{Le1934} and Hopf~\cite{Ho1951}. In two spatial dimensions, the uniqueness of weak solutions of NSE has been established. However, in three spatial dimensions, the question of uniqueness remains unanswered and has been designated as a Millennium Prize Problem by the Clay Mathematics Institute.

Keeping the uniqueness question aside, it is known that, for small viscosity $\nu$ (or, equivalently, large Reynolds number $Re \sim \nu^{-1}$), i.e. when inertial forces are much stronger than viscous diffusion, the fluid flow is very sensitive to initial conditions and characterized by chaotic motions.
Therefore, in a deterministic framework, measurement errors in the problem data could have a drastic affect on the solutions of NSE, and so, it is not physically meaningful to describe these turbulent flows as individual solutions.
On the other hand, there is ample experimental evidence that the statistical observables, for example, mean and variance, can be inferred reliably for turbulent flows.
Hence, for practical applications, it is more useful to consider the evolution of NSE under uncertainties on the problem data in a suitable probabilistic framework. 
In their seminal work \cite{Fo1972, FPr1976}, Foia\c{s} and Prodi proposed the framework of statistical solutions of incompressible Navier-Stokes equations, in which, given a probability distribution on the initial velocity, velocity ensembles are evolved according to the NSE and described by a time-parameterized family of probability measures on the function space of initial velocity. 
The existence of such statistical solutions, both in two and three spatial dimensions, has been well established \cite{Fo1972, FMRTe2001}. The global uniqueness in three spatial dimensions is an open problem. In two dimensions, statistical solutions are unique and they are defined as the push-forward of the probability measure on the initial velocity data.
The authors in \cite{FMRTe2001} show that several results of the conventional theory of turbulence for NSE attributed to the ground-breaking work of Kolmogorov \cite{Ko1941a, Ko1941b} can be recovered with statistical solutions, thus, providing evidence of the importance of this solution framework.

The computation of statistical solutions of NSE is a special case of Uncertainty Quantification (UQ) in Computational Fluid Dynamics (CFD), for which different methods are available in the literature (see \cite{BLMSc2013} and the references therein). Loosely speaking, the different methods used to solve UQ problems can be categorized into two classes: stochastic Galerkin methods and stochastic collocation methods. Stochastic Galerkin methods consider the Ritz-Galerkin formulation of the underlying PDEs in the stochastic parameter space, thus, these methods are highly intrusive and not amenable for implementation from a practical viewpoint \cite[Chap.~12]{Su2015}. On the other hand, stochastic collocation methods are well-suited for large scale applications as they are non-intrusive \cite[Chap.~13]{Su2015}. 
For applications, Monte Carlo (MC) and Quasi-Monte Carlo sampling can be seen as a subclass of stochastic collocation methods. Multi-level Monte Carlo methods have also received a lot of attention in research in recent times (see \cite{LMSc2016} and the references therein).

In the discussion above, the non-intrusive methods require a deterministic solver that numerically approximates the solution of the underlying system of PDEs (in our case the NSE).
There are many numerical methods available in the literature that can efficiently approximate the NSE. For problems with periodic boundary, spectral methods are appealing because of their efficiency and high resolution \cite{Le2018}.
Finite difference methods with Leray projection were first introduced in \cite{KLa1966, Ch1967, Ch1969}. Later, the authors in \cite{BCGl1989} proposed a finite-volume scheme that was an efficient variant of \cite{Ch1969}. 
In \cite{LMSc2016}, the authors develop a finite-difference scheme for the vorticity formulation of the NSE with uniformly stable bounds with respect to viscosity; this scheme is a variant of the one proposed in \cite{LTa1997} for incompressible Euler equations. 
Finite element methods are another class of methods to solve NSE \cite{GRa1986}. 
While standard mixed FEM are not pressure-robust, i.e. pressure approximation influences the velocity approximation, (weakly) divergence-free mixed methods are (see \cite{JLMNRe2017} and the references therein). 
In these weakly divergence-free mixed methods, if the velocity field belongs to the Hilbert space $H^{1}$, then the a priori error bounds are non-uniform with respect to viscosity. However, if H(div)-conforming spaces are used instead, then the a priori velocity error estimates are robust with respect to viscosity \cite{SLe2018}. 
Moreover, the divergence-free constraint on the velocity field is satisfied point-wise in the spatial domain.
The authors in \cite{LSc2016} proposed a computationally efficient hybridized variant of H(div) FEM, called the H(div)-HDG method, where they introduce hybridization for the tangential velocity components. 
This method was further improved in \cite{LLSc2018} and used to study turbulent flows in \cite{FKLLSc2019, SJLLLSc2019}.
\subsection{Contributions}
We restrict ourselves to two-dimensional open, bounded spatial domains and choose the framework of statistical solutions for NSE.
Given the discussion above, to numerically approximate these solutions, we use Monte Carlo sampling for the stochastic space and H(div)-FEM as the deterministic solver.
We can easily use QMC instead of MC, but it is not necessary as our main focus is on the presence or absence of convergence of the approximate solutions and their statistics, not the rates. We choose the H(div) scheme and not the more efficient H(div)-HDG scheme because of the ease of implementation. In the deterministic solver, we use the implicit Euler method for time integration as there is no CFL restriction on the time-step size (see \cite{MPRo1990, LSc2016} for higher-order time integration methods). 
To the best of my knowledge, this is the first computational effort in this direction.

Structure functions are an important observable in turbulence literature \cite{Fr1995}. In the context of incompressible Euler equations, a uniform decay of the structure functions ensures the convergence of approximate statistical solutions in the very recent work \cite{LMPa2021a}. The behaviour of structure functions is also crucial in the study of the vanishing viscosity limit of incompressible Navier-Stokes in \cite{LMPa2021b}.
In this light, a novel contribution of this work is the development of an algorithm for approximating structure functions on unstructured meshes. We compute the structure functions only at the final time, which is a good indicator of their behaviour over the whole time interval as seen in \cite{LMPa2021a}.
\subsection{Outline}
After covering some preliminaries, we state the initial boundary value problem (IBVP) and introduce the concept of weak solutions for NSE.
In {\S}\ref{s:incompNSStatSols}, we describe statistical solutions of NSE and define structure functions and Wasserstein distances. 
We present the deterministic H(div)-conforming scheme in {\S}\ref{s:uqIncompNSNumScheme}.
In {\S}\ref{s:uqIncompNSImplement}, we describe the algorithm for approximating structure functions on unstructured meshes.
We present the results of our numerical experiments in {\S}\ref{s:uqIncompNSNumExp} and conclusions in {\S}\ref{s:uqIncompNSConclusions}.

%
\section{Preliminaries}\label{s:uqIncompNSPrelims}
%
Consider an open time interval $\Dt := (0, T)$, where $T > 0$ denotes a finite time horizon, and let $\Dx \subset \R^2$ be an open and bounded polygon with Lipschitz boundary $\Dxb$ and with outward unit normal $\nDx$. Define the space-time cylinder $Q$ as the Cartesian product of the spatial domain $\Dx$ and the time domain $\Dt$, i.e. $Q := \Dx \times \Dt$. Let the scalar $t \in \Dt$ and the vector $\bx = (x_{1}, x_{2})^{\top} \in \Dx$ denote the time coordinate and the spatial Cartesian coordinates, respectively, then the space-time Cartesian coordinates are $(\bx,t)\in Q$.

We highlight vector fields by bold face fonts and tensor fields by bold face fonts with an underline.

Given $m,p \in \N$ and a vector $\bw = (w_{1},w_{2},\ldots, w_{m})^{\top} \in \R^{m}$,
the norm $\norm{\bw}_{p} = \left(\sum_{i=1}^{m}w_{i}^{p}\right)^{\frac{1}{p}}$, 
we also denote the Euclidean norm by $\abs{\bw} = \norm{\bw}_{2}$.

Given $m,n\in\N$ and two matrices $A, B \in \R^{m\times n}$ with entries $A_{i,j}, B_{i,j}$, $i\in\{1,\ldots,m\}$ and $j\in\{1,\ldots,n\}$,
we denote by $A:B$ the Frobenius inner product of the two matrices.
%
%
\subsection{Function spaces}\label{ss:uqIncompNSFes}
For any open, bounded domain $\frakD \subseteq \Dx$, we use the Lebesgue spaces $L^{p}(\frakD)$ for scalar-valued functions with the associated norm $\norm{\cdot}_{L^{p}(\frakD)}$, for $1 \leq p \leq \infty,\ p\in \N$. 
We also use the standard Hilbert spaces $H^{k}(\frakD)$ with associated norms $\norm{\cdot}_{H^{k}}$ and seminorms $\abs{\cdot}_{H^{k}}$, for $k\in\N$.
For vector-valued functions of size $m$, we indicate these spaces as $L^{p}(\frakD)^{m}$ and $H^{k}(\frakD)^{m}$. Spaces and norms for tensor-valued functions are indicated with bold font.
We use the standard notation $C(\frakD)$ for the space of continuous functions on $\frakD$ and $C^{k}(\frakD)$, $k\in\N$, for the space of functions on $\frakD$ that are $k$-times differentiable with continuous $k$-th derivative.
The space
\begin{align}\label{eq:HdivSpace}
\hdiv{\frakD} &:= \{\btau \in \ltwo{\frakD}^{2}: \nabla_{\bx}\cdot\btau\in\ltwo{\frakD}\}.
\end{align}
For a Hilbert space $X$ and $s\in\Nz$, we use the standard notation $H^{s}(\Dt;X)$ for Bochner spaces and define $\ltwo{\Dt;X} := H^{0}(\Dt;X)$.
\subsection{Probability}\label{ss:probability}
Let $(\Omega, \boldsymbol{\Omega}, P)$ be a probability space \cite{Sc2005}. 
For any random variable $X \in \ltwo{\Omega}$, its mean $\E(X)$ and variance $\mrm{Var}(X)$ are given by
\begin{subequations}
 \begin{align}
 \E[X] &= \int_{\Omega} X(\omega) dP(\omega),\\
 \mrm{Var}[X] &= \E[(X - \E[X])^{2}] = \E[X^{2}] - \E[X]^{2}.
 \end{align} 
\end{subequations}
We assume that $(\Omega, \boldsymbol{\Omega}, P)$ is complete.

Given a topological space $X$, we denote by $\calB(X)$ the Borel $\sigma$-algebra of $X$ and by $\calP(X)$ the space of all probability measures on $\calB(X)$.

\subsection{Meshes, mesh faces, averages and jumps}\label{ss:uqIncompNSPrelimsMeshes}
We choose a mesh $\Dxmh$ of domain $\Dx$ such that it is shape-regular and does not allow hanging nodes. We allow both simplicial meshes and meshes with quadrilateral elements,
which we also refer to as triangular and quadrilateral meshes, respectively.
The mesh-size of $\Dxmh$ is denoted by $\hx := \max_{\Kx \in \Dxmh} \hKx$.

We define the following sets and unions for mesh faces:
\begin{itemize}[nolistsep]
 \item $\F$ denotes the set of all faces of $\Dxmh$,
 \item $\F^{\mrm{int}}$ denotes the set of all interior faces of $\Dxmh$,
 \item $\F^{\mrm{bdr}}$ denotes the set of all boundary faces of $\Dxmh$,
 \item the mesh skeleton $\Dxmfh := \bigcup_{F \in \F} F$,
 \item the interior mesh skeleton $\Dxmfhi := \bigcup_{F \in \F^{\mrm{int}}} F$,
 \item the boundary mesh skeleton $\Dxmfhb := \bigcup_{F \in \F^{\mrm{bdr}}} F$.
\end{itemize}
Further, we divide $\F^{\mrm{bdr}}$ into two disjoint sets $\F^{\mrm{bdr},\mrm{D}}$ and $\F^{\mrm{bdr},\mrm{out}}$ that correspond to the boundaries with Dirichlet and outflow boundary conditions, respectively, in NSE \eqref{eq:ibvpIncompNS}.
Using these sets, we also define the unions $\DxmfhbD := \bigcup_{F \in \F^{\mrm{bdr},\mrm{D}}} F$ and $\Dxmfhbo := \bigcup_{F \in \F^{\mrm{bdr},\mrm{out}}} F$.

For $\Kx \in \Dxmh$, we denote the outward-pointing unit normal vector on its boundary $\Kxb$ by $\nxKx \in \R^2$.
Let $K_{1}, K_{2} \in \Dxmh$ be distinct.
Then, for a vector field $\bw$ and a tensor field $\underline{\btau}$, which are element-wise continuous on $\Dxmh$, we define the averages $\avg{\cdot}$ and jumps $\jump{\cdot}$ on the interior mesh face $F = \Kb_{1} \cap \Kb_{2} \neq \emptyset$ as follows:
\begin{equation}
 \begin{aligned}
  \avg{\bw} &:= \frac{\bw|_{K_{1}} + \bw|_{K_{2}}}{2},\\
  \avg{\underline{\btau}} &:= \frac{\underline{\btau}|_{K_{1}} + \underline{\btau}|_{K_{2}}}{2},\\
  \jump{\bw\otimes\bn_{F}} &:= \bw|_{K_{1}} \otimes \bn_{K_{1}} + \bw|_{K_{2}}\otimes \bn_{K_{2}},\\
  \jump{\bw\cdot\bn_{F}} &:= \bw|_{K_{1}} \cdot \bn_{K_{1}} + \bw|_{K_{2}}\cdot \bn_{K_{2}},\\
  \jjump{\bw} &:= \bw|_{K_{1}} - \bw|_{K_{2}}.
 \end{aligned}
\end{equation}
Here $\bn_{F}$ is the unit normal at face $F$ and we set $\bn_{F} = \bn_{K_{1}}$ by convention.
We extend the definitions of the jumps and averages for boundary faces.
Let $\Kx\in \Dxmh$ such that $F = \Kxb \cap \Dxmfhb \neq \emptyset$, then we define
\begin{equation}
 \begin{aligned}
  \avg{\bw}     &:= \bw|_{\Kx},&
  \avg{\underline{\btau}} &:= \underline{\btau}|_{\Kx},\\
  \jump{\bw\otimes\bn_{F}} &:= \bw|_{\Kx} \otimes \bn_{\Kx},&
  \jump{\bw\cdot\bn_{F}} &:= \bw|_{\Kx} \cdot \bn_{\Kx}.
 \end{aligned}
\end{equation}
We frequently drop the subscript in $\bn_{F}$ for convenience.
\subsection{Raviart-Thomas spaces}\label{ss:RTSpaces}
We revisit the definition of Raviart-Thomas spaces and some important results associated with them, which are presented in \cite[Chapter~3]{Ga2014} and \cite[Chapter~2]{BBFo2013}. We use these spaces to construct element-wise polynomial subspaces of $\hdivDx$ in our numerical scheme.

Given $k \in \Nz$, we denote by $\bbP_{k}$ the space of polynomials of total degree at most $k$,
and we define the discrete polynomial space
\begin{equation}\label{eq:WkSpace}
 \Wk := \Wk(\Dxmh) := \left\{ q \in \ltwo{\Dx} : q|_{\Kx} \in \bbP_{k}(\Kx)\ \forall \Kx \in \Dxmh \right\}.
\end{equation}
We denote Raviart-Thomas spaces of degree $k$ by $RT_{k}$, they are defined as follows:
\begin{definition}[Raviart-Thomas spaces]\label{def:RTSpace}
 For any element $\Kx\in\Dxmh$, the local Raviart-Thomas space of degree $k \in \Nz$ is defined as
 \begin{equation}\label{eq:localRTspace}
  \rt{\Kx}{k} := \bbP_{k}(\Kx)^{2} + \bbP_{k}(\Kx)\bx.
 \end{equation}
 Let $\globRt{k}$ denote the global Raviart-Thomas space, it is defined as
 \begin{equation}\label{eq:globalRTspace}
  \globRt{k} := \left\{\bv \in \hdivDx : \bv|_{\Kx} \in \rt{\Kx}{k}\ \forall\Kx \in \Dxmh\right\}.
 \end{equation}
\end{definition}
\begin{lemma}[{\cite[{\S}2.5.2]{BBFo2013}}]\label{lem:hdivProps}
 Any function $\bv\in\globRt{k}$ is continuous in the normal direction at all interior faces of the mesh, i.e.
 \begin{equation}\label{eq:vnContinuous}
  \jump{\bv\cdot\bn} = 0,\ \forall F \in \F^{\mrm{int}}.
 \end{equation}
 %
 Moreover,
 the following holds true:
 \begin{equation}\label{eq:divRT}
  \div_{\bx}\left(\globRt{k}\right) = \Wk(\Dxmh).
 \end{equation}
\end{lemma}
Consider the equation $\int_{\Dx} \left(\div_{\bx}(\bu) q\right)d\bx = 0$, for $\bu \in \globRt{k}$ and for all $q \in \Wk(\Dxmh)$, which appears in the numerical scheme we use in {\S}\ref{ss:HdivScheme}.
From Lemma~\ref{lem:hdivProps}, we deduce that $\div_{\bx}(\bu) = 0$ in $\Dx$. Thus, the divergence-free constraint in NSE can be imposed easily with the pair of FE spaces $\globRt{k}$ and $\Wk(\Dxmh)$.
\section{Navier-Stokes equations}\label{s:incompNS}
%
We assume that the boundary $\Dxb$ can be divided into Dirichlet boundary $\DxbD$ and outflow boundary $\Dxbo$, where either $\DxbD$ or $\Dxbo$ may be empty. 
We also assume that $\DxbD$ and $\Dxbo$ consist of entire segments of $\Dx$.

The IBVP for the time-dependent incompressible Navier-Stokes equations in the divergence form reads as follows:
\begin{subequations}\label{eq:ibvpIncompNS}
 \begin{align}
  \bu_t + \div_{\bx}(\bu \otimes \bu) - \nu \Delta_{\bx} \bu + \nabla_{\bx}{p} &= \bb{f} &\mathrm{in} &\quad Q,\label{eq:ibvpIncompNSEq1}\\
  \div_{\bx}{\bu} &= 0 &\mathrm{in} &\quad Q,\label{eq:ibvpIncompNSEq2}\\
  \bu(\cdot,0) &= \bu_0 &\mathrm{in} &\quad \Dx,\\
  \bu &= \bb{g} & \mathrm{on} &\quad \DxbD\times\Dt,\\
  (\nu\nabla_{\bx}\bu -p \underline{\bb{I}})\cdot\nDx &= 0 & \mathrm{on} &\quad \Dxbo\times\Dt \label{eq:ibvpIncompNSOutflowBC},
 \end{align}
\end{subequations}
where $\bu$ is the velocity, $p$ is the pressure,
$\nu$ is the kinematic viscosity, $\bu_{0}$ is the initial velocity, $\bb{f} \in \R^{2}$ is the forcing, $\bb{g} \in \R^{2}$ is the boundary data and $\underline{\bb{I}}$ is the identity tensor.
\begin{remark}\label{rem:stableOutflow}
For very large Reynolds numbers, a stabilized outflow boundary condition is needed to ensure numerical stability in the case of back-flow at the outlet, we refer to \cite{BCBBr2018} for a review of the different stabilization methods available in the literature and in particular to \cite{DKCh2014, Do2015} for some popular choices. 
In the present work, see {\S}\ref{s:uqIncompNSNumExp}, our numerical experiments for problems with outflow preclude any backward-flow at the outlet. 
Therefore, the outflow condition \eqref{eq:ibvpIncompNSOutflowBC} is sufficient for our purpose.
\end{remark}

%
\subsection{Weak solutions}\label{ss:incompNSWeak}
We assume homogeneous Dirichlet BCs, i.e. $\Dxbo = \emptyset$ and $\bb{g} = 0$.
We choose the functional setting used in \cite[{\S}II.5]{FMRTe2001} and define the following Hilbert spaces:
\begin{align}
 H := H_{nsp} &:= \{\bv \in \ltwo{\Dx}^{2} : \div_{\bx}\bv = 0,\ \bv\cdot\nDx|_{\Dxb} = 0\},\label{eq:spaceHnsp}\\
 V := V_{nsp} &:= \{\bv \in \hone{\Dx}^{2} : \div_{\bx}\bv = 0,\ \bv|_{\Dxb} = \bb{0}\}.\label{eq:spaceVnsp}
\end{align}
On these spaces, we define the inner products
\begin{align}
 (\bv, \bw)_{H} &:= \int_{\Dx} \bv\cdot\bw~d\bx,\ \mrm{for}\ \bv,\bw \in H,\label{eq:innerProdHnsp}\\
 (\bv, \bw)_{V} &:= \int_{\Dx} \nabla_{\bx}\bv : \nabla_{\bx}\bw~d\bx,\ \mrm{for}\ \bv,\bw \in V,\label{eq:innerProdVnsp}
\end{align}
and the associated norms
\begin{align}\label{eq:normHnspVnsp}
 \norm{\bv}_{H} := \left((\bv, \bv)_{H}\right)^{\half},\ \mrm{for}\ \bv\in H,\quad
 \norm{\bv}_{V} := \left((\bv, \bv)_{V}\right)^{\half},\ \mrm{for}\ \bv\in V.
\end{align}

The weak formulation reads as:
find $\bu \in \linf{\Dt; H}\cap \ltwo{\Dt; V}$ such that, $\forall\bv \in V$,
\begin{align}\label{eq:incompNSWeak}
 \frac{d}{dt}(\bu,\bv)_{H} + \nu (\bu, \bv)_{V} + (\bu\cdot\nabla_{\bx}\bu, \bv)_{\ltwo{\Dx}} = (\bb{f}, \bv)_{H},
\end{align}
with the viscosity $\nu > 0$, the initial velocity $\bu_{0} \in H$ and the forcing $\bb{f} \in \ltwo{\Dt; H}$.
This weak formulation can be traced back to the seminal work of Leray \cite{Le1933, Le1934}. 

It is well-known that in two-dimensional spatial domains with sufficiently smooth boundaries, for any $\nu > 0$, any initial velocity $\bu_{0} \in H$ and any forcing $\bb{f}\in\ltwo{\Dt; H}$, weak solutions of NSE exist and they are unique. We state this result in the following (cf. \cite[Chapter II, Theorem~7.3 and Remark~7.2]{FMRTe2001}):
\begin{theorem}[Existence and uniqueness of weak solutions of NSE]\label{thm:incompNSExistUnq}
Assume that the spatial domain $\Dx$ has $C^{2}$ boundary $\Dxb$, and that
\begin{align}
 \bu_{0} \in H,\quad \bb{f} \in \ltwo{\Dt; H}.
\end{align}
Then, for every $\nu > 0$, there exists a unique solution $\bu$ of \eqref{eq:incompNSWeak} such that
\begin{align}\label{eq:incompNSReg}
 \bu &\in \ltwo{\Dt; V} \cap C(\overline{\Dt}; H),
\end{align}
and, for all $t \in \Dt$, it satisfies the energy equation
\begin{align*}
 \half\norm{\bu(t)}_{H}^{2} + \nu \int_{0}^{t} \norm{\bu(s)}_{V}^{2} ds = \half \norm{\bu_{0}}_{H}^{2} + \int_{0}^{t} (\bb{f}(s), \bu(s))_{\Dx}.
\end{align*}
Moreover, there exists a continuous solution operator $\calS^{\nu}: H \rightarrow H$ such that $\bu(t) = \calS^{\nu}(t)\bu_{0}$.
\end{theorem}

%
\section{Statistical solutions} \label{s:incompNSStatSols}
%
The discussion in this section is based on the monograph \cite[Chapter~V]{FMRTe2001} and it partially uses the notation prescribed in \cite[{\S}2.2]{LMSc2016}.
We consider statistical solutions of NSE \eqref{eq:ibvpIncompNS} with fixed, no-slip boundary, i.e. $\Dxbo = \emptyset$ and $\bb{g} = 0$, 
where the boundary $\Dxb$ is of class $C^{2}$.

Given a probability distribution $\mu_{0} \in \calP(H)$ on the initial velocity data, i.e.
\begin{align*}
 P(\{\bu_{0} \in E\}) = \mu_{0}(E),\ \mrm{for}\ E \in \calB(H),
\end{align*}
the main idea of statistical solutions, as introduced by Foia\c{s}-Prodi \cite{Fo1972, FPr1976}, is to describe the evolution of velocity ensembles by a time-parameterized family of probability measures $\mu^{\nu} = (\mu^{\nu}_{t})_{t\in\overline{\Dt}}$, such that $\mu^{\nu}_{0} = \mu_{0}$ (up to null sets), and for every $t > 0$, $\mu^{\nu}_{t} \in \calP(H)$.
This solution framework includes the individual weak solutions of NSE as a special case with the probability measure $\mu_{0} = \delta_{\bu_{0}}$, where $\delta_{\bu_{0}}$ is the Dirac measure (also known as the unit mass) at $\bu_{0}$, cf. \cite[{\S}4.7]{Sc2005}.

From Theorem~\ref{thm:incompNSExistUnq}, we know that the individual weak solutions of NSE are unique and the solution operator $\calS^{\nu}$ exists. 
As a result, for any time-independent forcing $\bb{f}\in H$, unique solutions $\mu^{\nu}_{t}$ are given by transporting the initial probability measure $\mu_{0}$ under the solution operator $\calS^{\nu}$, i.e.
\begin{align}
 \mu^{\nu}_{t}(E) = \mu_{0}(\left(\calS^{\nu}(t)\right)^{-1} E),\ \forall E \in \calB(H).
\end{align}
Concretely, we state this existence and uniqueness result in the following (cf. \cite[Chapter~V, Theorems~1.1 and 1.2]{FMRTe2001}):
\begin{theorem}\label{thm:statsIncompNSExistUnq}
 Let $\mu_{0}$ be a probability measure on $H$ with finite kinetic energy, i.e.
 \begin{align*}
  \int_{H} \norm{\bu}_{H}^{2} d\mu_{0}(\bu) < \infty.
 \end{align*}
 Assume that the forcing term $\bb{f} \in \ltwo{\Dt; H}$. Then, for every $\nu > 0$ and homogeneous Dirichlet boundary conditions, there exists a statistical solution $\mu^{\nu} = (\mu_{t}^{\nu})_{t\in\overline{\Dt}}$ of the incompressible Navier-Stokes equations \eqref{eq:ibvpIncompNS} on $H$.
 
 Moreover, if $\mu_{0}$ has bounded support in $H$ and $\bb{f} \in H$ is independent of time, then the statistical solution $\mu^{\nu}$ is unique and is given by $\mu^{\nu}_{t} = \calS^{\nu}(t) \mu_{0}$, where $\calS^{\nu}(t)$ is the solution operator of the incompressible Navier-Stokes equations, cf. Theorem~\ref{thm:incompNSExistUnq}.
\end{theorem}
\subsection{Structure functions} \label{s:incompNSStatSolsStructure}
For a solution $\mu^{\nu}_{t}$, other than its common statistics like mean and variance, its structure functions are of great interest to us. 
The structure functions are defined as follows (we refer to \cite{Ly2020}):
\begin{definition}[Structure functions]\label{def:structureFn}
 Let $r \in \R^{+}$, $p \in \N$ and $t \in \overline{\Dt}$. 
 Let $\bu: \Omega\times Q \rightarrow \R^{2}$ be a random field such that $\bu(\omega; \cdot,t) \in L^{p}(\Dx)^{2}$, for any $\omega \in \Omega$, and $\mu_{t} \in \calP(L^{p}(\Dx)^{2})$ be the probability distribution of $\bu$ at time $t$.
 Then, the structure function $S^{p}_{r,t}$ associated with $\mu_{t}$ is given by:
\begin{align}\label{eq:structureFn}
S^{p}_{r,t}(\mu_{t}) = \left( \int_{L^{p}(\Dx)^{2}} \int_{\Dx} \fint_{B_{r}(\bx)} \norm{\bu(\cdot;\bx, t) - \bu(\cdot;\by, t)}_{p}^{p} d\by\ d\bx\ d\mu_{t}(\bu) \right)^{\frac{1}{p}},
\end{align}
where $B_{r}(\bx) \subset \R^{2}$ denotes the ball of radius $r$ centered at a point $\bx \in \Dx$.
\end{definition}
We use the structure functions $S^{p}_{r,t}(\mu^{\nu}_{t})$ as a measure of the regularity of the velocity field $\bu$, and we investigate the scaling behaviour of $S^{p}_{r,t}$ with respect to the offset $r$ in our numerical experiments.
\begin{remark}\label{rem:structureFnScaling}
 If $\bu$ is Lipschitz-continuous in $\Dx$, then straightforward calculations show that $S^{p}_{r,t}(\mu^{\nu}_{t}) = O(r)$. We expect the same scaling behaviour for the velocity field in NSE due to regularity \eqref{eq:incompNSReg} under the assumptions of Theorem~\ref{thm:incompNSExistUnq}.
\end{remark}

\subsection{Wasserstein distances} \label{s:incompNSStatSolsWasserstein}
In the next section, we describe a numerical method to approximate statistical solutions with ensembles of approximate individual velocity solutions. To study the convergence of these velocity ensembles, it is important to define a notion of distance between two ensembles.
For this purpose, we use Wasserstein distance (cf. \cite[Definition~2.2]{FLMi2017}):
\begin{definition}
 Let $X$ be a separable Banach space, and let $\mu, \rho$ be probability measures on $X$ with finite $p$-th moments, i.e. $\int_{X} \abs{x}^{p}d\mu(x) < \infty$ and $\int_{X} \abs{x}^{p}d\rho(x) < \infty$. Then, the $p$-Wasserstein distance between $\mu$ and $\rho$ is defined as
 \begin{align*}\label{eq:WassersteinDist}
  W_{p}(\mu, \rho) := \left(\inf_{\gamma \in \Gamma(\mu,\rho)} \int_{X^{2}} \abs{x - y}^{p}~d\gamma(x, y) \right)^{\frac{1}{p}},
 \end{align*}
 where the infimum is taken over the set $\Gamma(\mu, \rho) \subset \calP(X^{2})$ of all transport plans from $\mu$ to $\rho$, i.e. those $\gamma \in \calP(X^{2})$ that satisfy
 \begin{align*}
 \int_{X^{2}} \left(F(x) + G(y)\right) d\gamma(x,y) = \int_{X} F(x) d\mu(x) + \int_{X} G(y) d\rho(y),\quad \forall F,G \in C_{b}(X).
 \end{align*}
 Here $C_{b}(X)$ denotes the space of bounded, continuous, real-valued functionals on $X$.
\end{definition}

Given $\mu, \rho \in \calP(H)$, we can define $p$-Wasserstein distance for the $1$-point distribution at $\bx$ and the $2$-point correlation at $\bx,\by$, respectively, as:
\begin{align}
 W_{p}(\mu(\bx), \rho(\bx)) 
 &:= \left(\inf_{\gamma \in \Gamma(\mu,\rho)} \int_{H^{2}} \norm{\bu(\bx) - \bv(\bx)}^{p}_{2} d\gamma(\bu, \bv) \right)^{\frac{1}{p}},\\
 W_{p}(\mu(\bx,\by), \rho(\bx,\by)) 
 &:= \left(\inf_{\gamma \in \Gamma(\mu,\rho)} \int_{H^{2}} 
 \norm{ \begin{bmatrix}
         \bu(\bx)\\ \bu(\by)
        \end{bmatrix}
        - 
        \begin{bmatrix}
         \bv(\bx)\\ \bv(\by)
        \end{bmatrix}
      }^{p}_{2} d\gamma(\bu, \bv) \right)^{\frac{1}{p}}.
\end{align}
For the whole domain $\Dx$, we can define
\begin{align}\label{eq:WassersteinDistW12}
 W_{p}^{1}(\mu,\rho) &:= \int_{\bx\in \Dx} W_{p}(\mu(\bx), \rho(\bx))~d\bx,\\
 W_{p}^{2}(\mu,\rho) &:= \int_{(\bx,\by)\in \Dx^{2}} W_{p}(\mu(\bx,\by), \rho(\bx,\by))~d\by d\bx.
\end{align}
In this work, we always use $1$-Wasserstein distances $W_{1}^{1}$ and $W_{1}^{2}$, so we simply refer to them as Wasserstein distances and drop the subscripts.

%
\section{Numerical approximation of statistical solutions}\label{s:uqIncompNSNumScheme}
%
In this section, our goal is to approximate the statistical solutions $\mu^{\nu}_{t}$ described in the previous section. To this end, we assume that initial probability measure $\mu_{0}$ is given as the law of a random field $\bu_{0} \in \ltwo{\Omega; H}$ defined on the underlying probability space $(\Omega, \boldsymbol{\Omega}, P)$.
This allows us to draw random samples from the initial distribution, which can then be evolved using a numerical solution operator. We estimate the statistical solutions using \eqref{eq:approxStatSol}. Concretely, this Monte Carlo type method is described in Algorithm~\ref{alg:monteCarlo}.
\begin{algorithm}[!htb]
\SetAlgoLined
\KwData{Initial distribution $\mu_{0}$, number of samples $M \in \N$, numerical solution operator $\calS^{\nu}_{h}$, mesh-size $h > 0$, time $t$}
\KwResult{Approximate statistical solution $\mu^{\nu}_{t;h,M}$}
 Sample independent random variables $\{\bu_{0;m}\}_{m=1,\ldots,M}$ from the distribution $\mu_{0}$.\\
 \For {$m=1,2,\ldots,M$} {
    Evolve $m$-th sample to the time $t$, $\bu_{h;m} = \calS^{\nu}_{h}(t)\bu_{0;m}$.
  }
  Compute the empirical measure
  \begin{align}\label{eq:approxStatSol}
   \mu^{\nu}_{t;h,M} = \frac{1}{M}\sum_{m=1}^{M} \delta_{\bu_{h; m}(\omega_{m};\cdot,t)}.
  \end{align}
\caption{Monte Carlo algorithm}
\label{alg:monteCarlo}
\end{algorithm}
Here $\calS^{\nu}_{h}$ is a numerical approximation of the operator $\calS^{\nu}$ such that $\displaystyle\lim_{h\rightarrow 0} \calS^{\nu}_{h} = \calS^{\nu}$.

To approximate an individual solution, we choose an H(div)-conforming method based on the work of Sequeira et al. \cite{GSSe2017} and Cockburn et al. \cite{CKSc2007}; this scheme has also been considered in \cite{SLe2018, SLu2018}.

For spatial discretisation \S\ref{ss:HdivScheme}, we use the H(div)-conforming scheme designed in \cite{GSSe2017} for the incompressible Euler equations along with the symmetric interior penalty discretisation for the viscous diffusion given in \cite{CKSc2007} and \cite[\S4.2.2]{DEr2012}.
This semi-discrete scheme is pressure-robust and $Re$-semi-robust, we refer to \cite{SLe2018} for details, which means that the velocity error bounds have no explicit dependence on the pressure approximation and the Reynolds number. 
To obtain a fully discrete scheme {\S}\ref{ss:HdivSchemeImpEuler}, we combine the H(div) spatial discretisation with implicit Euler time-stepping as formulated in \cite{GSSe2017}.

\subsection{H(div)-conforming spatial discretisation}\label{ss:HdivScheme}
We use the function space $\Vk$, which is defined as 
\begin{equation}\label{eq:fesRT}
  \Vk := \Vk(\Dxmh) := \left\{\bv \in \globRt{k} : \bv\cdot\nDx = 0\ \mathrm{on}\ \DxbD\right\},
 \end{equation}
for the velocity field.
Here the space $\globRt{k}$ is given by \eqref{eq:globalRTspace}.
In the case of purely Dirichlet boundary, we denote $\Vk$ by $\Vzk$.
We use the function space $\Wk$ \eqref{eq:WkSpace} for the pressure field.
\begin{remark}\label{rem:zeroMeanPressure}
 For purely Dirichlet boundary, i.e. $\Dxbo = \emptyset$, if $(\bu, p)$ is a solution of \eqref{eq:ibvpIncompNS} such that $p \in \ltwo{\Dx}$, then $(\bu, p + p_{0})$, for any constant $p_{0} \in \R$, is also a solution. Therefore, an additional constraint is needed to ensure a unique pressure solution. A common choice is to impose vanishing mean for the pressure (cf., for example, \cite{GSSe2017, JLMNRe2017}) using the function space
\begin{subequations}\label{eq:WkzSpace}
 \begin{align}
 \Wzk := \Wzk(\Dxmh) &:= \left\{ q \in \ltwoz{\Dx} : q|_{\Kx} \in \bbP_{k}(\Kx)\ \forall \Kx \in \Dxmh \right\},\\
 \mrm{where}\ \ltwoz{\Dx} &:= \left\{q \in \ltwo{\Dx} : \int_{\Dx} q~d\bx = 0\right\}.
 \end{align}
\end{subequations}
\end{remark}
Given initial velocity data $\bu_{0}\in H$, let $\bu_{0,h}$ denote the $L^{2}$-orthogonal projection of $\bu_{0}$ into $\Vk$. Then, the H(div)-conforming spatial discretisation of the IBVP \eqref{eq:ibvpIncompNS} reads as follows:

Given the initial data $\bu_{0,h}$, the forcing $\bb{f}$ and the boundary data $\bb{g}$, 
find $\buh\in\hone{\Dt;\Vk}$ and $p_{h}\in\ltwo{\Dt;\Wk}$ such that, for all $\bvh\in\Vk$ and for all $q_h\in\Wk$,
\begin{subequations}\label{eq:hdivFem}
\begin{align}
 \bfMass{\der{\buh}{t}}{\bvh} + \bfConvU{\buh}{\buh}{\bvh} + \nu\bfDiffSip{\buh}{\bvh} + \bfDiv{\bvh}{p_{h}} &= \lfSource{\bvh},\label{eq:hdivFem1}\\
 \bfDiv{\buh}{q_{h}} &= 0\label{eq:hdivFem2},
\end{align}
\end{subequations}
where the mass bilinear form
\begin{equation}
\label{eq:bfMass}
 \bfMass{\bu}{\bv} := \int_{\Dxmh} \bu \cdot \bv~ d\bx,\quad \mrm{for}\ \bu,\bv \in \Vk,
\end{equation}
the convection trilinear form with upwind flux
\begin{equation}\label{eq:bfConvU}
\begin{aligned}
 \bfConvU{\bw}{\bu}{\bv} 
 &:=\ \sum_{\Kx \in \Dxmh} \int_{\Kx} (\bw \cdot \nabla_{\bx}{\bu})\cdot \bv~d\bx 
 - \int_{\Dxmfhi} (\bw \cdot \bn) \jjump{\bu}\cdot \avg{\bv}~d\bx & \\
 &\ + \int_{\Dxmfhi} \abs{\bw \cdot \bn}\jjump{\bu}\cdot \jjump{\bv}~d\bx, \quad \mathrm{for}\ \bw, \bu, \bv \in \Vk,
\end{aligned}
\end{equation}
where $\bw\cdot\bn|_{F} = \bw_{K_{1}}\!\cdot\bn_{K_{1}} = -\bw_{K_{2}}\!\cdot\bn_{K_{2}}$ (cf. \eqref{eq:vnContinuous}),
the symmetric interior penalty diffusion bilinear form
\begin{equation}\label{eq:bfDiff}
\begin{aligned}
\bfDiffSip{\bu}{\bv} &:= \sum_{\Kx \in \Dxmh}\int_{\Kx} \nabla_{\bx}{\bu} : \nabla_{\bx}{\bv}~d\bx\\
 &\ - \sum_{F \in \F^{\mrm{int}} \cup \F^{\mrm{bdr},\mrm{D}}} 
 \int_{F} \jump{\bv \otimes \bn} : \avg{\nabla_{\bx} \bu}~d\bx \\
 &\ - \sum_{F \in \F^{\mrm{int}} \cup \F^{\mrm{bdr},\mrm{D}}} 
 \int_{F} \jump{\bu \otimes \bn} : \avg{\nabla_{\bx} \bv}~d\bx \\
 &\ + \sum_{F \in \F^{\mrm{int}} \cup \F^{\mrm{bdr},\mrm{D}}}\hspace{-5mm}
 \penParam\hFx^{-1} \int_{F} \jump{\bu \otimes \bn} : \jump{\bv \otimes \bn}~d\bx,\quad \mathrm{for}\ \bu,\bv \in \Vk,
\end{aligned}
\end{equation}
where $\penParam \in \R^{+}$ is the jump penalization parameter and $\hFx$ is the size of face $F$,
the divergence bilinear form
\begin{equation}\label{eq:bfDiv}
 \bfDiv{\bu}{q} := \sum_{\Kx \in \Dxmh} \int_{\Kx} q \div_{\bx}(\bu)~d\bx,\quad \mrm{for}\ \bu \in \Vk,\ q \in \Wk,
\end{equation}
and the linear form
\begin{equation}\label{eq:lfSource}
\begin{aligned}
 \lfSource{\bv} &:= \int_{\Dxmh} \bb{f} \cdot \bv\ d\bx\\
                        &\ - \nu \sum_{F \in \F^{\mrm{bdr},\mrm{D}}} \int_{F} (\bb{g} \otimes \bn) : \nabla_{\bx}{\bv}\ ds\\
                        &\ + \nu~\penParam \sum_{F \in \F^{\mrm{bdr},\mrm{D}}} \hFx^{-1}\int_{F} (\bb{g} \otimes \bn) : (\bv \otimes \bn)~d\bx,\quad \mathrm{for}\ \bv \in \Vk.
\end{aligned}
\end{equation}
\begin{remark}
 In the formulation \eqref{eq:hdivFem}, the outflow boundary conditions are imposed weakly through the boundary integrals on $\Dxbo$ that appear in the bilinear forms $\bfDiffSip{\buh}{\bvh}$ and $\bfDiv{\bvh}{p_{h}}$.
\end{remark}
\begin{remark}\label{rem:hdivFemErrorViscous}
 For homogeneous Dirichlet boundary conditions, the scheme \eqref{eq:hdivFem} is the same as the one analyzed in \cite{SLe2018}, except for the slight difference in the convection trilinear term that does not affect the error estimates.
 Therefore, for $\nu > 0$, from Theorem~5.3, Theorem~5.6 and Corollary~5.9 in \cite{SLe2018},
 we deduce that the error $\norm{(\bu - \buh)(\cdot,T)}_{\ltwo{\Dx}^{2}}$ for the discrete solution $\buh$ of \eqref{eq:hdivFem} is bounded by $O(h^{k})$, $k\geq 1$, for $\bu$ with sufficient regularity.
\end{remark}
%

%
\subsection{Fully-discrete scheme}\label{ss:HdivSchemeImpEuler}
As mentioned, we combine the semi-discrete formulation \eqref{eq:hdivFem} with implicit Euler time-stepping to obtain a fully discrete scheme {\S}\ref{ss:HdivSchemeImpEuler}, by following the procedure described in \cite{GSSe2017}.

For $N \in \N$, let $\{t_{n}\}_{n=0,1,\ldots,N}$ be a set of $N+1$ points in the time interval $\Dt$ such that $t_{n} = n \Delta t$ with $\Delta t = T / N$.
Given the solution $(\bu, p)$ of the IBVP \eqref{eq:ibvpIncompNS}, 
we define the notation $\bu^{n} := \bu(\cdot, \tn)$ and $p^{n} := p(\cdot, \tn)$, and we denote by $\buh^{n}$ and $p_{h}^{n}$ the numerical approximations of $\bu^{n}$ and $p^{n}$, respectively.

The Taylor series expansion of $\bu(\cdot, t_{n-1})$ at $\tn$ is given by
\begin{align}\label{eq:unm1TaylorSeries}
 \bu(\cdot, t_{n-1}) &= \bu(\cdot, \tn) - \Delta t \der{\bu}{t}(\cdot,\tn) + \frac{(\Delta t)^{2}}{2}\dder{\bu}{t}(\cdot,\tn) + \ldots\nonumber\\
 &= \bu(\cdot, \tn) + O(\Delta t),
\end{align}
which yields the time derivative
\begin{align}\label{eq:impDudt}
 \der{\bu}{t}(\cdot, \tn) &= \frac{\bu(\cdot, t_{n}) - \bu(\cdot, t_{n-1})}{\Delta t} + \bb{E}_{0}(\tn),
\end{align}
with the truncation error
\begin{align}\label{eq:impDudtTruncation}
 \bb{E}_{0}(\tn) &:= \frac{\Delta t}{2}\dder{\bu}{t}(\cdot,\tn) + \ldots = O(\Delta t).
\end{align}
Substituting \eqref{eq:impDudt} in \eqref{eq:ibvpIncompNS}, 
at any discrete time $\tn$, we obtain
\begin{align*}
 \bu^{n} + \Delta t \left(\bu^{n}\cdot\nabla_{\bx}\bu^{n} - \nu \Delta_{\bx} \bu^{n} + \nabla_{\bx}{p^{n}}\right) &= \bu^{n-1} + \Delta t\ \bb{f}(\cdot,\tn) + \Delta t\ \bb{E}_{0}(\tn),\\
  \div_{\bx}{\bu^{n}} &= 0,
\end{align*}
which can be written equivalently as
\begin{equation}\label{eq:incompNSTimeDisc}
\begin{aligned}
 \bu^{n} + \Delta t \left(\bu^{n-1}\cdot\nabla_{\bx}\bu^{n} - \nu \Delta_{\bx} \bu^{n} + \nabla_{\bx}{p^{n}}\right) &= \bu^{n-1} + \Delta t\ \bb{f}(\cdot,\tn)\\ 
  &\quad + \underbrace{\Delta t((\bu^{n-1}-\bu^{n})\cdot\nabla_{\bx}\bu^{n} + \bb{E}_{0}(\tn))}_{=:\ \tilde{\bb{E}}_{0}(\tn)},\\
  \div_{\bx}{\bu^{n}} &= 0.
\end{aligned}
\end{equation}
Using \eqref{eq:unm1TaylorSeries} and \eqref{eq:impDudtTruncation}, we deduce that the term $\tilde{\bb{E}}_{0}(\tn) = O([\Delta t]^{2})$. Dropping the truncation term $\tilde{\bb{E}}_{0}(\tn)$ in \eqref{eq:incompNSTimeDisc}, the remaining equation is linear with respect to the variable $\bu^{n}$. We introduce H(div)-conforming spatial discretisation in the truncated system \eqref{eq:incompNSTimeDisc} to arrive at a \emph{fully discrete method}, which reads as follows:

Given the initial data $\bu_{0,h}$, the forcing $\bb{f}$ and the boundary data $\bb{g}$, 
find $\buh^{n}\in\Vk$ and $p_{h}^{n}\in\Wk$ such that, for all $\bvh\in\Vk$ and for all $q_h\in\Wk$,
\begin{subequations}\label{eq:hdivFemImpEuler}
\begin{align}
 \bfMass{\buh^{n}}{\bvh} &+ \Delta t \left(\bfConvU{\buh^{n-1}}{\buh^{n}}{\bvh} + \nu\bfDiffSip{\buh^{n}}{\bvh} + \bfDiv{\bvh}{p_{h}^{n}}\right)\nonumber\\ 
 &= \bfMass{\buh^{n-1}}{\bvh} + \Delta t~ \lfSource{\bvh},\\
 \bfDiv{\buh^{n}}{q_{h}} &= 0,
\end{align}
\end{subequations}
for $n=1,\ldots, N$, with $\buh^{0} = \bu_{0,h}$.
\begin{remark}
 The discrete formulation \eqref{eq:hdivFemImpEuler} satisfies the following:
 \begin{align*}
  \norm{\buh^{n}}_{\ltwo{\Dx}^{2}}^{2} + &\Delta t \left(\abs{\buh^{n}}_{\buh^{n-1},up}^{2} + \nu C_{\penParam} \norm{\buh^{n}}_{e}^{2}\right) \nonumber\\
  &\leq (\buh^{n-1}, \buh^{n})_{\ltwo{\Dx}^{2}} + \Delta t(\bb{f}, \buh^{n})_{\ltwo{\Dx}^{2}},
 \end{align*}
 for $n=1,\ldots, N$. Then, for the forcing $\bb{f} = \bb{0}$, we have 
 \begin{align*}
  \norm{\buh^{n}}_{\ltwo{\Dx}^{2}} \leq \norm{\buh^{n-1}}_{\ltwo{\Dx}^{2}}.
 \end{align*}
 Therefore, the scheme is $\ltwo{\Dx}$-stable.
\end{remark}
%
%
\begin{remark}\label{rem:hdivFemImpEulerErrorViscous}
 For viscosity $\nu > 0$, to the best of my knowledge, error estimates of our numerical scheme are not available, and it is beyond the scope of this work to derive them.
 However, based on the reasoning in the proofs of \cite[Lemma~5.5, Theorem~5.6]{SLu2018} and \cite[Lemma~3.4, Theorem~3.6]{SLe2018}, 
 we expect the velocity error to be bounded by $O(h^{k} + \Delta t)$.
\end{remark}

%
\section{Implementation details}\label{s:uqIncompNSImplement}
%
In this section, we discuss important details about our implementation of the Monte Carlo algorithm \ref{alg:monteCarlo}, which we refer to as MC-FEM solver; 
here all the samples are evolved in parallel using the self-scheduling algorithm \cite[{\S}3.6]{GLSk2014}.
In particular, we describe the deterministic solver that is used for a single Monte Carlo sample and present the algorithms used to approximate structure functions and Wasserstein distances.

We denote vectors by a lower-case font with an underline and matrices by upper-case font or bold font with an underline, only for the discussion on the deterministic solver.
Given bases of $\Vk(\Dxmh)$ and $\Wk(\Dxmh)$, we denote the corresponding vectors of basis coefficients by $\vec{u} \in \R^{\dim(\Vk)}$ and $\vec{p} \in \R^{\dim(\Wk)}$.
Then, the finite element assembly for the discrete formulation \eqref{eq:hdivFemImpEuler} yields the following block-structured linear system of equations:
\begin{align}\label{eq:assemblyIncompNSHdiv}
 \begin{bmatrix}
    {M}_{h} + \Delta t ({C}_{h} + \nu {A}_{h})
    & \Delta t {B}_{h}^{\top}
    \\
    {B}_{h}
    & \mat{0}
 \end{bmatrix}
 \begin{bmatrix}
    \vec{u}
    \\
    \vec{p}
 \end{bmatrix}
  = \begin{bmatrix}
    \vec{f}_{h}
    \\
    \vec{0}
 \end{bmatrix},
\end{align}
where the matrices ${M}_{h}$, ${C}_{h}$, ${A}_{h}$ and ${B}_{h}$ are assembled from the mass form $\bfMasse$, the convection form $\bfConvUe$, the diffusion form $\bfDiffSipe$ and the divergence form $\bfDive$, respectively. Here the right-hand side vector $\vec{f}_{h}$ is assembled from the linear form $\ell_{h}$. Note that we must use the Piola transformation in our finite element assembly.

We need to assemble the matrices ${M}_{h}$, ${A}_{h}$ and ${B}_{h}$ only once, but the matrix ${C}_{h}$ has to be re-assembled for each discrete time $\tn$ because of its dependence on the velocity solution of the previous time step $t_{n-1}$. Using direct solvers for the linear system \eqref{eq:assemblyIncompNSHdiv} would involve factorizing the re-assembled system matrix at each time step, which is computationally expensive. Therefore, we use suitable iterative solvers with preconditioning. Re-scaling the pressure solution as $\widetilde{\vec{p}} := (\Delta t)\vec{p}$, the linear system \eqref{eq:assemblyIncompNSHdiv} can be re-written in the form
\begin{align}\label{eq:assemblyIncompNSHdivRescaled}
 \begin{bmatrix}
    {M}_{h} + \Delta t ({C}_{h} + \nu {A}_{h})
    & {B}_{h}^{\top}
    \\
    {B}_{h}
    & \mat{0}
 \end{bmatrix}
 \begin{bmatrix}
    \vec{u}
    \\
    \widetilde{\vec{p}}
 \end{bmatrix}
  = \begin{bmatrix}
    \vec{f}_{h}
    \\
    \vec{0}
 \end{bmatrix},
\end{align}
which is more suitable for iterative solvers in practice. 

The linear system \eqref{eq:assemblyIncompNSHdivRescaled} is in saddle point form. We refer to the monograph \cite{BGLi2005} for an extensive review of the literature on the numerical solution of saddle point problems. Note that the matrices ${M}_{h}$ and ${A}_{h}$ are symmetric, but ${C}_{h}$ is non-symmetric, which implies that the linear system \eqref{eq:assemblyIncompNSHdivRescaled} is also non-symmetric. So, we choose the GMRES method \cite{SSc1986} to solve our non-symmetric linear system, with the state-of-the-art block triangular preconditioner described in \cite[{\S}{10.1.2}]{BGLi2005}.

When needed, see Remark~\ref{rem:zeroMeanPressure}, we impose the vanishing mean constraint for pressure through the iterative solver via an intermediate preconditioning matrix
\begin{align}\label{eq:zeroMeanPressurePr}
\begin{bmatrix}
 \mat{I} & \mat{0}\\
 \mat{0} & \mat{I} - \frac{1}{\abs{\Dx}}({M}_{h}^{p})^{-1} \vec{b}_{h}\vec{b}_{h}^{\top}
\end{bmatrix},
\end{align}
where $\abs{\Dx}$ is the volume of $\Dx$, $M_{h}^{p}$ denotes the matrix that corresponds to the bilinear form $\int_{\Dxmh}p_{h} q_{h} d\bx$, for $p_{h}, q_{h} \in \Wk(\Dxmh)$, $\vec{b}_{h}$ denotes the vector that corresponds to the linear form $\int_{\Dxmh} q_{h} d\bx$, for $q_{h} \in \Wk(\Dxmh)$ and $\mat{I}$ is the identity matrix.

With the components presented above, we develop our MC-FEM solver in the programming language C++, using the finite element library MFEM \cite{mfem}, version 4.1, with MPI-parallelized domain decomposition.
We use the mesh generation routines provided by MFEM and the open-source simplicial mesh generator Gmsh \cite{Gmsh}, version 2.10.0. 
The codes developed can be found at \url{https://github.com/pratyuksh/NumHypSys}.

\subsection{Computing the structure functions}\label{ss:computeStructureFns}
We want to compute the structure functions \eqref{eq:structureFn} for the approximate statistical velocity solutions $\mu^{\nu}_{t; h,M}$. 
An efficient algorithm to compute $S^{p}_{r,t}(\mu^{\nu}_{t;h,M})$ on uniform grids was presented by Lye \cite{Ly2020}. 
In the present work, our objective is to compute structure functions on unstructured meshes for rectangular domains $\Dx = (x_{1}^{\mrm{left}}, x_{1}^{\mrm{right}})\times(x_{2}^{\mrm{left}}, x_{2}^{\mrm{right}}) \subset \R^{2}$, 
with $x_{1}^{\mrm{left}} < x_{1}^{\mrm{right}}$ and $x_{2}^{\mrm{left}} < x_{2}^{\mrm{right}}$, which we commonly use in our numerical experiments in {\S}\ref{s:uqIncompNSNumExp}.

By definition, for any $t \in \overline{\Dt}$, we have that
\begin{align*}
(S^{p}_{r,t}(\mu^{\nu}_{t; h,M}))^{p} = \int_{L^{p}(\Dx)^{2}} \int_{\Dx} \fint_{B_{r}(\bx)} \norm{\bu(\cdot;\bx, t) - \bu(\cdot;\by, t)}_{p}^{p} d\by\ d\bx\ d\mu^{\nu}_{t; h,M}(\bu),
\end{align*}
and using the fact that $\mu^{\nu}_{t; h,M}$ is a sum of Dirac masses (see \eqref{eq:approxStatSol}),
\begin{align*}
(S^{p}_{r,t}(\mu^{\nu}_{t; h,M}))^{p} 
&= \int_{L^{p}(\Dx)^{2}} \int_{\Dx} \fint_{B_{r}(\bx)} 
\norm{\bu(\cdot;\bx, t) - \bu(\cdot;\by, t)}_{p}^{p} d\by\ d\bx\  d\left(\frac{1}{M} \sum_{m=1}^{M} \delta_{\bu_{h, m}(\omega_{m}; \cdot, t)}\right)\\
&= \frac{1}{M} \sum_{m=1}^{M} \int_{\Dx} \fint_{B_{r}(\bx)} 
\norm{\bu_{h,m}(\omega_{m}; \bx, t) - \bu_{h,m}(\omega_{m}; \by, t)}_{p}^{p} d\by\ d\bx\\
&= \frac{1}{M} \sum_{m=1}^{M} S^{p}_{r,t}(\bu_{h,m}),
\end{align*}
where we define
\begin{equation}\label{eq:structureFnOfSampleExact}
 S^{p}_{r,t}(\bu_{h,m}) := \int_{\Dx} \fint_{B_{r}(\bx)} 
\norm{\bu_{h,m}(\omega_{m}; \bx, t) - \bu_{h,m}(\omega_{m}; \by, t)}_{p}^{p} d\by\ d\bx.
\end{equation}
The discrete solutions $\bu_{h,m}$ are only first-order accurate because of implicit Euler time-stepping, cf. Remark~\ref{rem:hdivFemImpEulerErrorViscous}. Therefore, we can further approximate $\bu_{h,m}$ with an element-wise constant function $\overline{\bu}_{h,m}$,
defined in every element $\Kx\in\Dxmh$ as
\begin{align}\label{eq:velAvgs}
 \overline{\bu}_{h,m}(\omega_{m}; \bx, t) := \frac{1}{\abs{\Kx}}\int_{\Kx} {\bu}_{h,m} (\omega_{m}; \by, t)~d\by,\ \forall \bx \in \Kx.
\end{align}
It can be shown with straightforward calculations that this element-averaged function $\overline{\bu}_{h,m}$ is a first-order approximation of the element-wise polynomial function $\bu_{h,m}$. 
In the numerical experiments {\S}\ref{s:uqIncompNSNumExp}, we compute the structure functions at time $t=T$, so we drop the subscript corresponding to time in $S_{r,T}^{p}$ in the text that follows, and by ${\bu}_{h,m}$ and $\overline{\bu}_{h,m}$ we simply refer to the approximate solutions at the final time.
Thus, our task boils down to computing $S^{p}_{r}({\bu}_{h,m}) \approx S^{p}_{r}(\overline{\bu}_{h,m})$, for $m \in \{1,2,\ldots,M\}$.

The primary component required to compute $S^{p}_{r}(\overline{\bu}_{h,m})$ is a nearest neighbor search (NNS) algorithm that helps us efficiently search for mesh elements lying within a ball $B_{r}$.
NNS problem arises in various fields, for example, statistics, machine learning and astronomy, we refer to, for example, \cite{GriSPy} and the references therein. A popular choice to solve NNS problems is the so-called ``cell-techniques'' method with fixed-radius searches, \cite{GriSPy, Fr2006}.
This fixed-radius NNS algorithm, adapted to our setting, is as follows:
\begin{enumerate}[nolistsep]
 \item Equi-partition the domain $\Dx$ into a uniform grid of cells such that the size of each cell is greater than or equal to the offset $r$ of the ball $B_{r}$. Map each element $\Kx \in \Dxmh$ to a cell in the uniform grid according to the coordinates of its centroid.
 \item Loop over all the cells in the uniform grid, and for each cell, loop over all the mesh elements mapped to it. 
 For each element $\Kx$, we find the neighbouring mesh elements lying within the ball $B_{r}$ by searching the current cell and its immediate neighbouring cells.
\end{enumerate}
Concretely, in the above-listed, the first step is implemented in Algorithm~\ref{alg:makeHashTable} using a hash table data structure, 
where the function call $\mrm{element\_number}(\Dxmh, \Kx)$ returns the element number of $\Kx \in \Dxmh$, $\mrm{centroid}(\Kx)$ returns the centroid of the mesh element $\Kx$, $\lfloor\cdot\rfloor$ is the floor function and the last entry in the tuple $(id_{\Kx}, \abs{\Kx}, \bx_{\Kx}^{c}, \vec{0})$ corresponds to the velocity data of $\Kx$ (set to zero by default). 
To compute $S^{p}_{r}(\overline{\bu}_{h,m})$, we update the hash table with the velocity data $\overline{\bu}_{h,m}$ using Algorithm~\ref{alg:updateHashTable} and then use Algorithm~\ref{alg:computeStructureSingleSample}. 
Put together, we use Algorithm~\ref{alg:computeStructure} to compute the structure function for the velocity ensemble $\{\overline{\bu}_{h,m}\}_{m=1,\ldots,M}$.
\begin{algorithm}[!htb]
\SetAlgoLined
\KwData{Domain $\Dx = (x_{1}^{\mrm{left}}, x_{1}^{\mrm{right}})\times(x_{2}^{\mrm{left}}, x_{2}^{\mrm{right}})$, mesh $\Dxmh$, grid size $N_{x_1}\!\times\!N_{x_2}$}
\KwResult{Hash table $\calM$}
 // map unstructured mesh elements to $N_{x_1}\!\times\!N_{x_2}$ uniform grid on domain $\Dx$\\
 $\calM := \{\calM_{i,j}\}$, with $\calM_{i,j}= \emptyset$, for $i=0,1,\ldots,N_{x_1}-1$ and $j=0,1,\ldots,N_{x_2}-1$.\\
 \For {$\Kx \in \Dxmh$} {
    $id_{\Kx} = \mrm{element\_number}(\Dxmh, \Kx)$\\
    $\bx_{\Kx}^{c} = (x_{\Kx;1}^{c}, x_{\Kx;2}^{c})^{\top} = \mrm{centroid}(\Kx)$\\
    Compute $\displaystyle i = \left\lfloor N_{x_1}\frac{x_{\Kx;1}^{c} - x_{1}^{\mrm{left}}}{x_{1}^{\mrm{right}} - x_{1}^{\mrm{left}}}\right\rfloor$, 
    $\displaystyle j = \left \lfloor N_{x_2}\frac{x_{\Kx;2}^{c} - x_{2}^{\mrm{left}}}{x_{2}^{\mrm{right}} - x_{2}^{\mrm{left}}}\right\rfloor$.\\
    Add the tuple $(id_{\Kx}, \abs{\Kx}, \bx_{\Kx}^{c}, \vec{0})$ to the set $\calM_{i,j}$.
  } 
\caption{make\_hash\_table$(\Dx, \Dxmh, N_{x_1}, N_{x_2})$}
\label{alg:makeHashTable}
\end{algorithm}
\begin{algorithm}[!htb]
\SetAlgoLined
\KwData{Domain $\Dx = (x_{1}^{\mrm{left}}, x_{1}^{\mrm{right}})\times(x_{2}^{\mrm{left}}, x_{2}^{\mrm{right}})$, mesh $\Dxmh$, grid size $N_{x_1}\!\times\!N_{x_2}$, element-averaged velocity sample $\overline{\bu}_{h}$, hash table $\calM$}
\KwResult{Hash table $\calM$}
 // update velocity data\\
 \For {$\Kx \in \Dxmh$} {
    $id_{\Kx} = \mrm{element\_number}(\Dxmh, \Kx)$\\
    $\bx_{\Kx}^{c} = (x_{\Kx;1}^{c}, x_{\Kx;2}^{c})^{\top} := \mrm{centroid}(\Kx)$\\
    Compute $\displaystyle i = \left\lfloor N_{x_1}\frac{x_{\Kx;1}^{c} - x_{1}^{\mrm{left}}}{x_{1}^{\mrm{right}} - x_{1}^{\mrm{left}}}\right\rfloor$, 
    $\displaystyle j = \left \lfloor N_{x_2}\frac{x_{\Kx;2}^{c} - x_{2}^{\mrm{left}}}{x_{2}^{\mrm{right}} - x_{2}^{\mrm{left}}}\right\rfloor$.\\
    Search $id_{\Kx}$ in $\calM_{i,j}$ and update the corresponding tuple to $(id_{\Kx},\cdot, \cdot, \overline{\bu}_{h}(\bx_{\Kx}^{c}, T))$.
  } 
\caption{update\_hash\_table$(\Dx, \Dxmh, N_{x_1}, N_{x_2}, \overline{\bu}_{h}, \calM)$}
\label{alg:updateHashTable}
\end{algorithm}
\begin{algorithm}[!htb]
\SetAlgoLined
\KwData{Hash table $\calM$, grid size $N_{x_1}\!\times\!N_{x_2}$, offset $r$, degree $p$}
\KwResult{Structure function $S^{p}_{r}$}
 Initialize $S^{p}_{r} = 0$.\\
 \For {$j=1,\ldots,N_{x_2}-2$} {
  \For {$i=1,\ldots,N_{x_1}-2$} {
    \For {$(\cdot, \abs{\Kx}, \bx_{\Kx}^{c}, \overline{\bv}_{\Kx}) \in \calM_{i,j}$} {
     Initialize $sum = 0$, $w = 0$.\\
     \For {$jj = j-1,j,j+1$} {
      \For {$ii = i-1,i,i+1$} {
        \For {$(\cdot, \abs{\Kx^{\prime}}, \bx_{\Kx^{\prime}}^{c}, \overline{\bv}_{\Kx^{\prime}}) \in \calM_{ii,jj}$} {
          \If {$\abs{\bx_{\Kx;1}^{c} - \bx_{\Kx^{\prime};1}^{c}} \leq r$ and $\abs{\bx_{\Kx;2}^{c} - \bx_{\Kx^{\prime};2}^{c}} \leq r$} {
            ${sum} \mathrel{+}= \abs{\Kx^{\prime}}
            \left(\abs{\overline{v}_{\Kx;1} - \overline{v}_{\Kx^{\prime};1}}^{p} + \abs{\overline{v}_{\Kx;2} - \overline{v}_{\Kx^{\prime};2}}^{p}\right)$\\
            ${w} \mathrel{+}= \abs{\Kx^{\prime}}$
          }
        }
      }
     }
     $S^{p}_{r} \mathrel{+}= \abs{\Kx} ({sum} / w)$
    }
  }
 }
\caption{compute\_structure\_function\_of\_a\_sample$(\calM, N_{x_1}, N_{x_2}, r, p)$}
\label{alg:computeStructureSingleSample}
\end{algorithm}
\begin{algorithm}[!htb]
\SetAlgoLined
\KwData{Domain $\Dx = (x_{1}^{\mrm{left}}, x_{1}^{\mrm{right}})\times(x_{2}^{\mrm{left}}, x_{2}^{\mrm{right}})$, mesh $\Dxmh$, grid size $N_{x_1}\!\times\!N_{x_2}$, 
ensemble of element-averaged velocities $\{\overline{\bu}_{m}\}_{m=1,\ldots,M}$ on $\Dxmh$, offset $r$, degree $p$}
\KwResult{Structure function $S^{p}_{r}(\overline{\mu}_{M})$, where $\overline{\mu}_{M} := \frac{1}{M}\sum_{m=1}^{M}\delta_{\overline{\bu}_{m}}$.}
 Initialize $S^{p}_{r}(\overline{\mu}_{M}) = 0$.\\
 $\calM = \mrm{make\_hash\_table}(\Dx, \Dxmh, N_{x_1}, N_{x_2})$\\
 \For {$m=1,2,\ldots,M$} {
    update\_hash\_table$(\Dx, \Dxmh, N_{x_1}, N_{x_2}, \overline{\bu}_{m}, \calM)$\\
    $S^{p}_{r}(\overline{\mu}_{M}) \mathrel{+}= \mrm{compute\_structure\_function\_of\_single\_sample}(\calM, N_{x_1}, N_{x_2}, r, p)$
  }
 Update $S^{p}_{r}(\overline{\mu}_{M}) \leftarrow \left(\frac{1}{M} S^{p}_{r}(\overline{\mu}_{M})\right)^{1/p}$.
\caption{compute\_structure\_function$(\Dx, \Dxmh, N_{x_1}, N_{x_2}, \{\overline{\bu}_{m}\}_{m=1,\ldots,M}, r, p)$}
\label{alg:computeStructure}
\end{algorithm}

We implement both serial and parallel versions of Algorithm~\ref{alg:computeStructure}.
In the parallel version, we uniformly distribute the ensemble between processors, compute the structure function for the local ensembles and then sum the values computed by all the processors. Note that the update $S^{p}_{r}(\overline{\mu}_{M}) \leftarrow \left(\frac{1}{M} S^{p}_{r}(\overline{\mu}_{M})\right)^{1/p}$ in Algorithm~\ref{alg:computeStructure} has to be done after the reduction operation in the parallel version.
\begin{remark}\label{rem:computeStructureFnsGeneralDomain}
 The algorithms presented above for computing the structure functions can also be applied, with some modifications, to more general spatial domains $\Dx$, given that $\Dx$ can be partitioned into rectangular sub-domains.
\end{remark}
\subsection{Computing the Wasserstein distances}\label{ss:computeWasserstein}
We use Algorithm~A.2.2 in Lye \cite{Ly2020} to approximate the Wasserstein distances $W^{1}$ and $W^{2}$ between two ensembles. To compute the $p$-Wasserstein distance, we use the function \texttt{ot.emd2} from the \texttt{POT} \cite{pot} module in Python.

%
\section{Numerical experiments}\label{s:uqIncompNSNumExp}
%

In the numerical experiments in \S\ref{ss:uqIncompNSTest2} and \S\ref{ss:uqIncompNSTest3}, we compute statistical solutions for the flow in a lid-driven cavity and in a rectangular channel, respectively. 
The settings in these experiments have been adopted from \cite{Le2018} with slight modifications.
We denote by $\mathcal{U}[-1,+1]$ the uniform distribution between $-1$ to $+1$, and we consider smooth perturbations of the initial data.

In our MC-FEM solver, there are two sources of error: stochastic error due to MC sampling and FE discretisation error. It is well-known that the MC stochastic error scales as $O(M^{-\half})$ with respect to the number of samples $M$ \cite{VWe1996}. The discretisation error scales as $O(\Delta t + h^{k})$ with respect the time step size $\Delta t$ and the mesh-size $h$ (cf. Remark~\ref{rem:hdivFemImpEulerErrorViscous}). 
As the main concern of the present work is the presence or absence of convergence of statistics and observables of velocity ensembles, rather than their convergence rates, we choose $\Delta t$, $h$ and $M$ such that the stochastic error dominates. To this end, we fix the polynomial degree $k=1$, which yields $O(h^{2})$ convergence for the deterministic velocity approximation of smooth solutions, in practice.

For an approximate statistical solution $\mu^{\nu}_{T;h,M}$, we compute the mean
\begin{align*}
 E_{M}[\bu_{h,T}] = \frac{1}{M}\sum_{m=1}^{M} \bu_{h;m}(\omega_{m}; \cdot, T)
\end{align*}
and the (unbiased) sample variance
\begin{align*}
 \mrm{Var}_{M}[\bu_{h,T}] = \frac{M}{M-1}\left(E_{M}[(\bu_{h,T})^{2}] - E_{M}[\bu_{h,T}]^{2} \right).
\end{align*}
We generally do not know the exact statistical solution $\mu^{\nu}_{T}$ for a given problem. Therefore, to study the convergence of the mean and variance of their approximations, we measure the Cauchy error in the $L^{2}(\Dx)$ norm, i.e., for example, 
$\norm{E_{M}[\bu_{h,T}] - E_{2M}[\bu_{h/2,T}]}_{\ltwo{\Dx}}$ in the lid-driven cavity problem.

We denote velocity ensembles by $\bb{U}^{h,M}_{T} := \{\bu_{h;m}(\omega_{m}; \cdot, T)\}_{m=1,\ldots,M}$.

%
\subsection{Lid-driven cavity}\label{ss:uqIncompNSTest2}
%
\renewcommand{\figsDir}{./}
We choose the spatial domain $\Dx := (0,1)^2$, the time horizon $T = 1$ and the Reynolds number $\nu^{-1} = Re = 3200$.
We consider the unperturbed initial velocity
\begin{equation*}
\label{eq:uqIncompNSTest2InitSol}
 \begin{aligned}
  \tilde{\bu}_{0}(\bx) = \left(x_2 - 0.5, -(x_1 - 0.5)\right)^{\top},\ \forall \bx\in\Dx.
 \end{aligned}
\end{equation*}
For $K \in \N$ such that $K$ {is an odd integer}, we consider, {for $\bx \in \Dx$}, the random perturbation function
\begin{equation*}
\begin{aligned}
f(\omega; \bx) 
&= \left(\ x_1 + \gamma_1 \sum_{k=0}^{(K-1)/2} Y_{2k}(\omega) \sin(2\pi k(x_1-0.5 + Y_{2k+1}(\omega)))\right.,\\
&\left.\qquad x_2 + \gamma_1 \sum_{k=0}^{(K-1)/2} Y_{2k+1}(\omega) \sin(2\pi k(x_2-0.5 + Y_{2k}(\omega)))\ \right),
\end{aligned}
\end{equation*}
where $Y_{j} \sim \mathcal{U}[-1,1]$ are independent, identically distributed random variables.
Then, we define the initial velocity $\bu_{0} = \tilde{\bu}_{0}\circ f$
and set $\gamma_1 = 0.025$, $\gamma_2 = 0.01$ and $K = 11$.
We fix the source term $\bb{f} = \bb{0}$.
The left, bottom and right boundaries are fixed, while for the top boundary we impose
$\bu = (1 + \gamma_2 \sin(2\pi Y_{K}),0)^{\top}$.

\begin{table}[h]
\centering
 \begin{tabular}{c | c | c}
 Resolution & Number of time steps & Number of samples\\
 \hline
 $32 \times 32$ & 100 & 32\\
 $64 \times 64$ & 200 & 64\\
 $128 \times 128$ & 400 & 128\\
 $256 \times 256$ & 800 & 256\\
 $512 \times 512$ & 1600 & 512\\
 \end{tabular}
 \caption{\small {
 The number of time steps and the number of samples used in the numerical solver for different mesh resolutions in the lid-driven cavity problem described in {\S}\ref{ss:uqIncompNSTest2}.
 }}
 \label{tab:uqIncompNSTest2Nsteps}
\end{table}
We use uniform quadrilateral meshes with different resolutions in this experiment, 
which are given in Table~\ref{tab:uqIncompNSTest2Nsteps} along with the corresponding number of time steps used in the solver. 

\begin{figure}[!htb]
\centering
  \begin{subfigure}[t]{.4\textwidth}
    \centering\includegraphics[width=\textwidth]{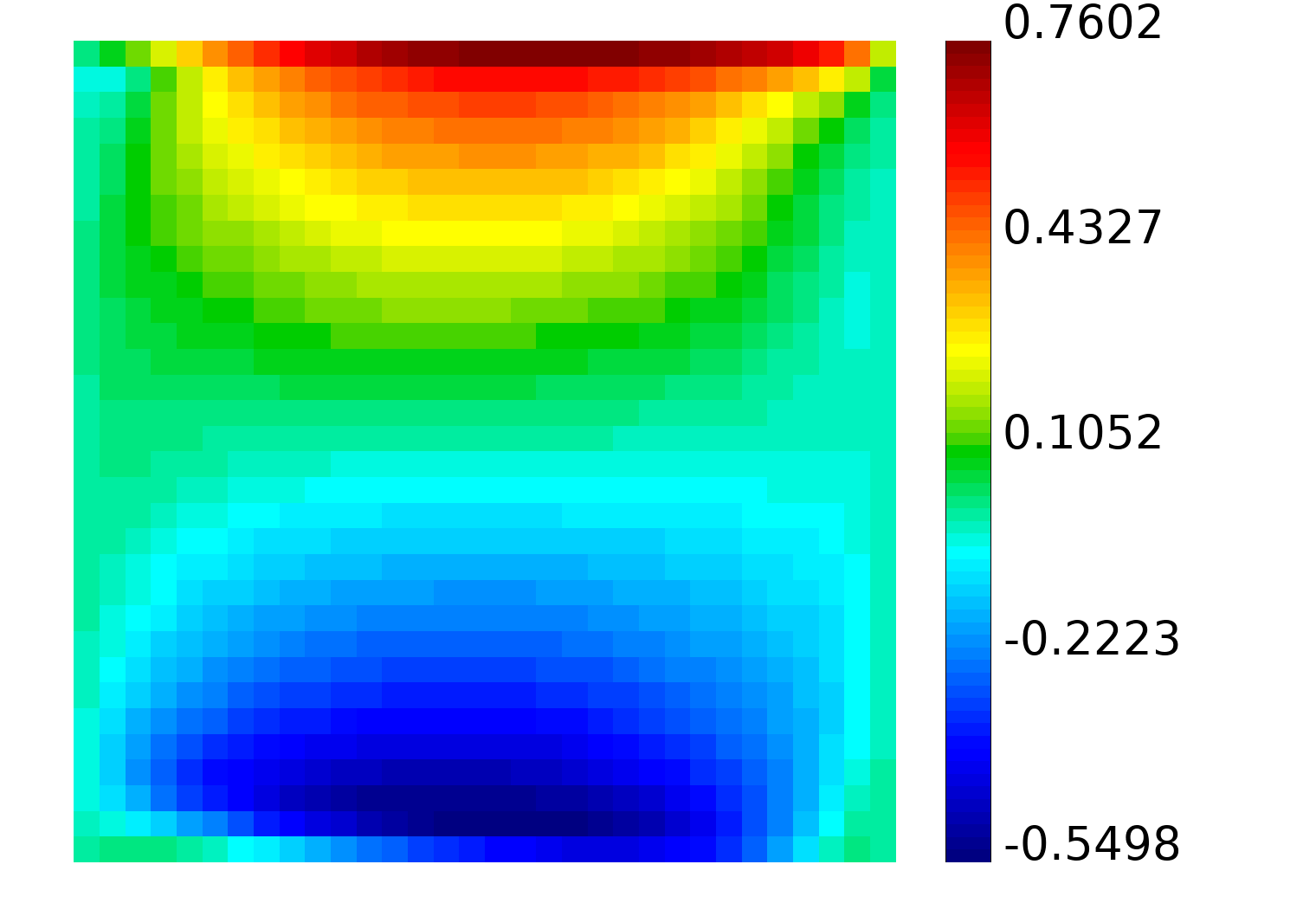}
  \end{subfigure}%
  \begin{subfigure}[t]{.4\textwidth}
    \centering\includegraphics[width=\textwidth]{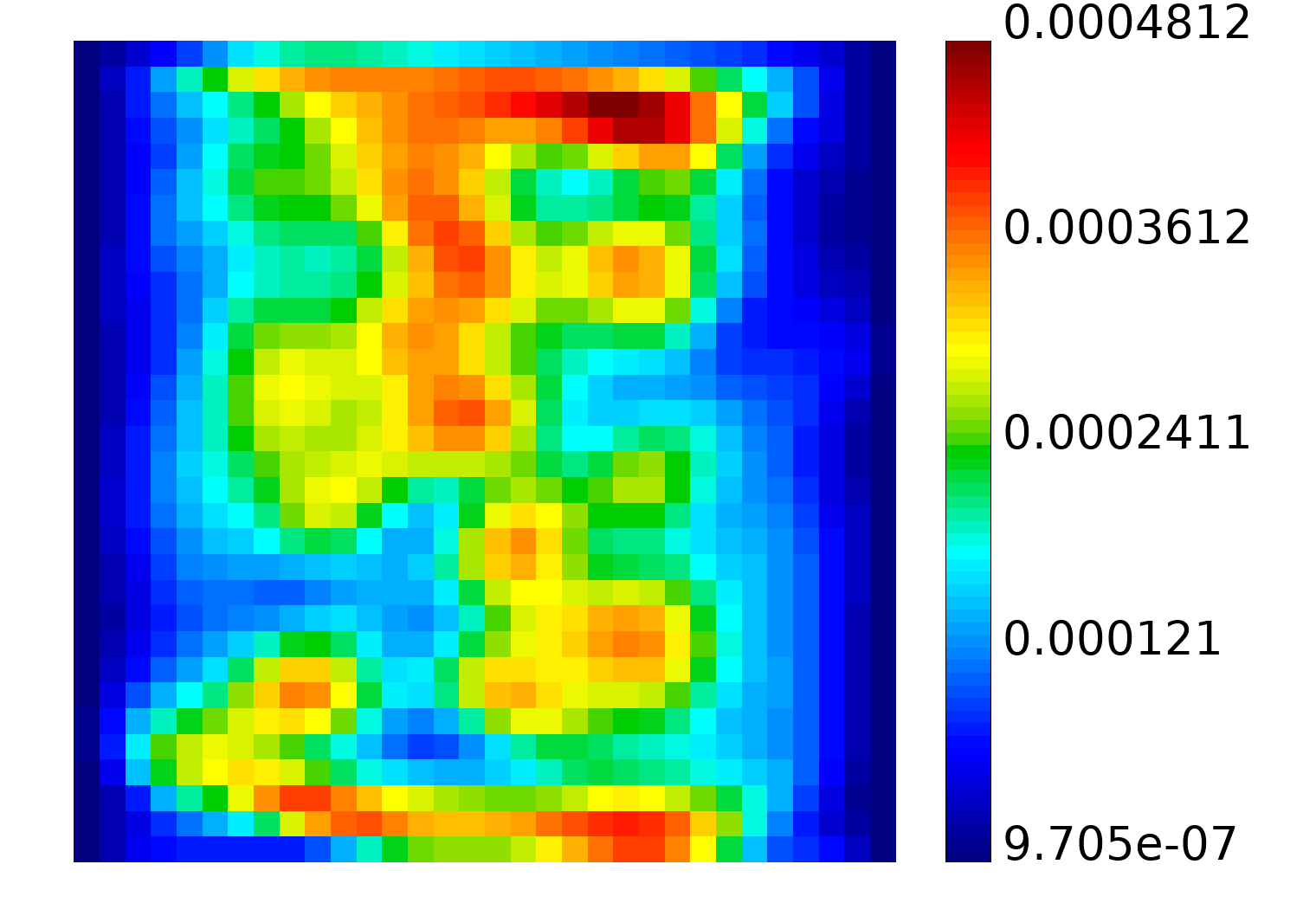}
  \end{subfigure}
  \begin{subfigure}[t]{.4\textwidth}
    \centering\includegraphics[width=\textwidth]{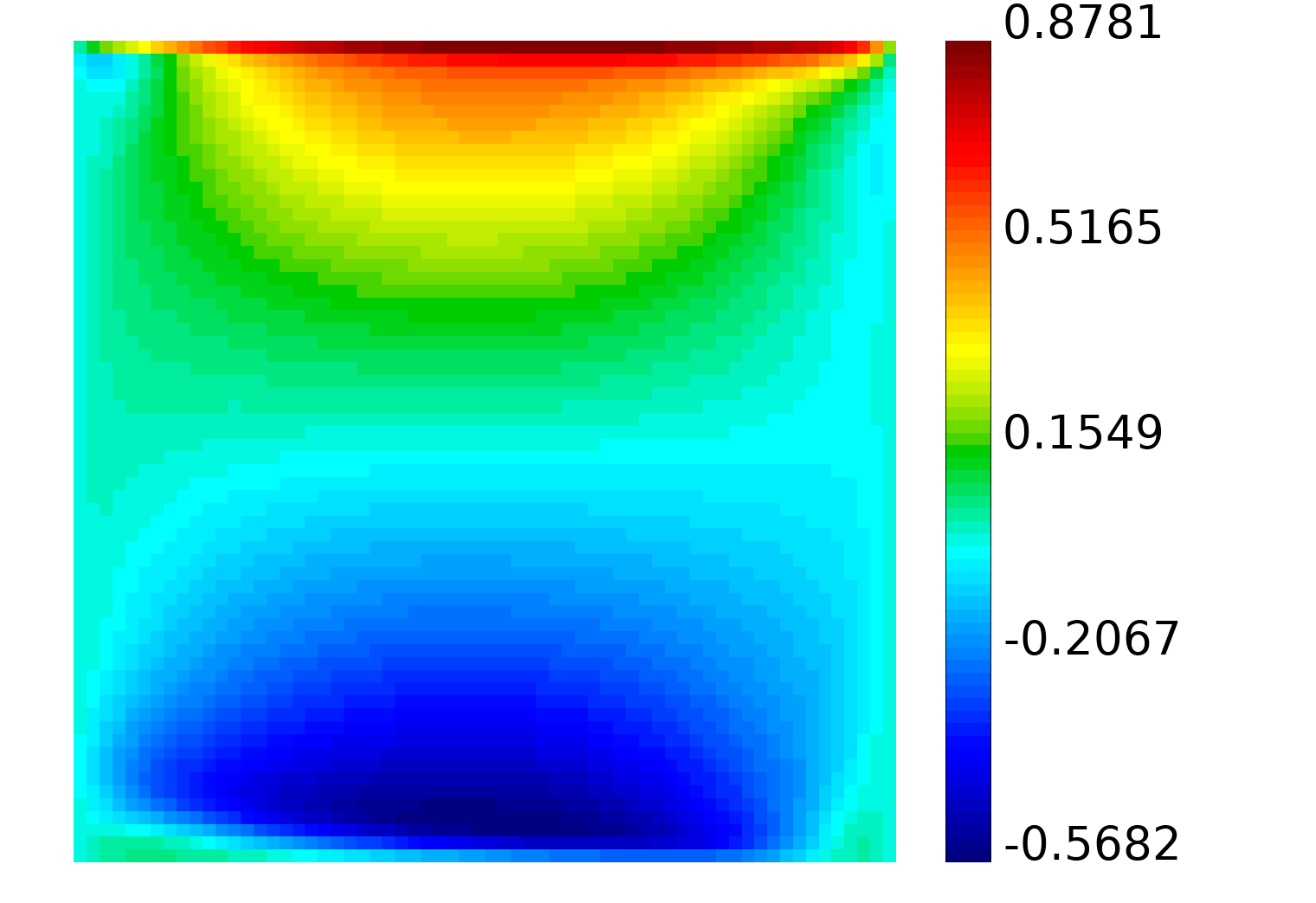}
  \end{subfigure}%
  \begin{subfigure}[t]{.4\textwidth}
    \centering\includegraphics[width=\textwidth]{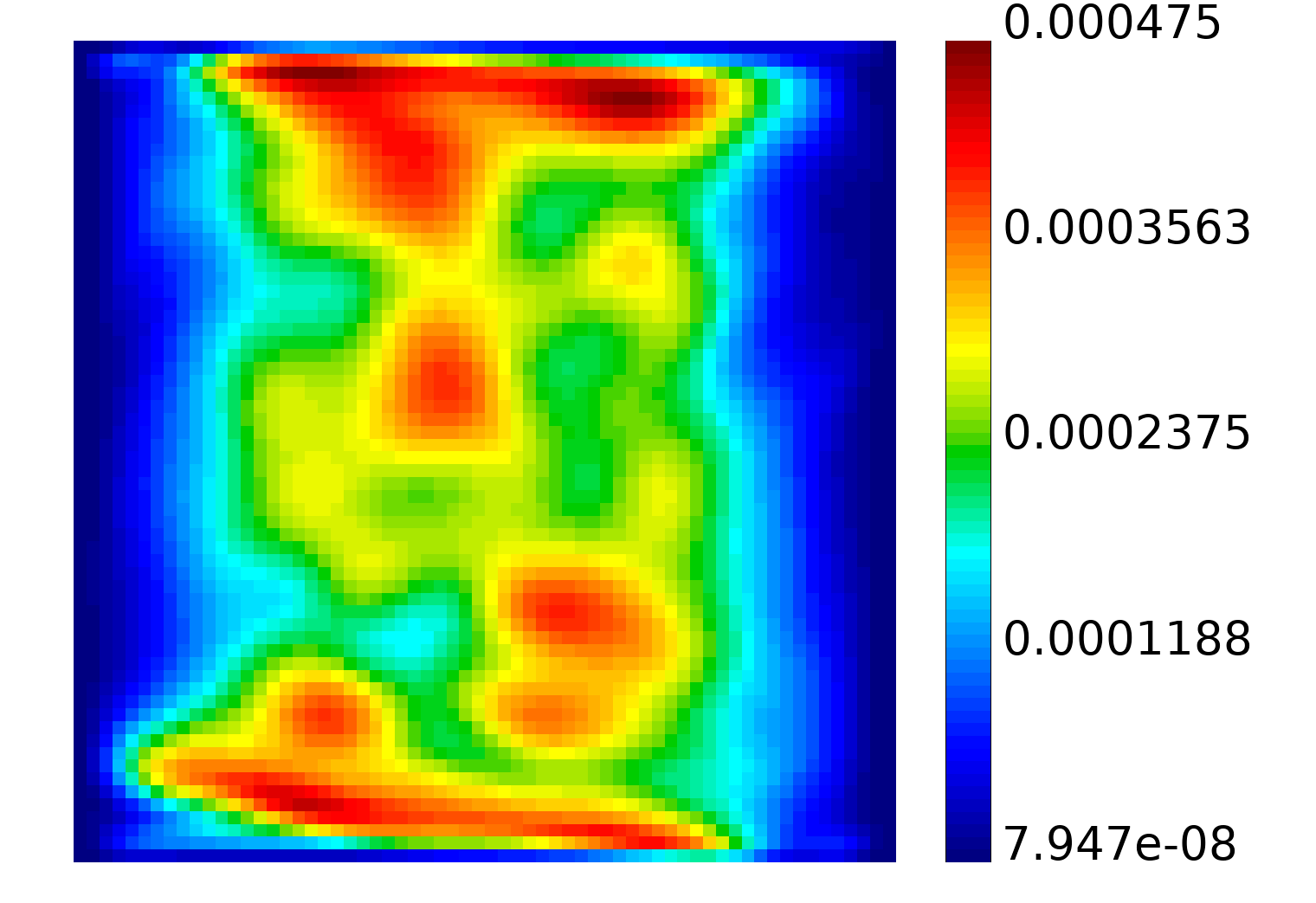}
  \end{subfigure}
  \begin{subfigure}[t]{.4\textwidth}
    \centering\includegraphics[width=\textwidth]{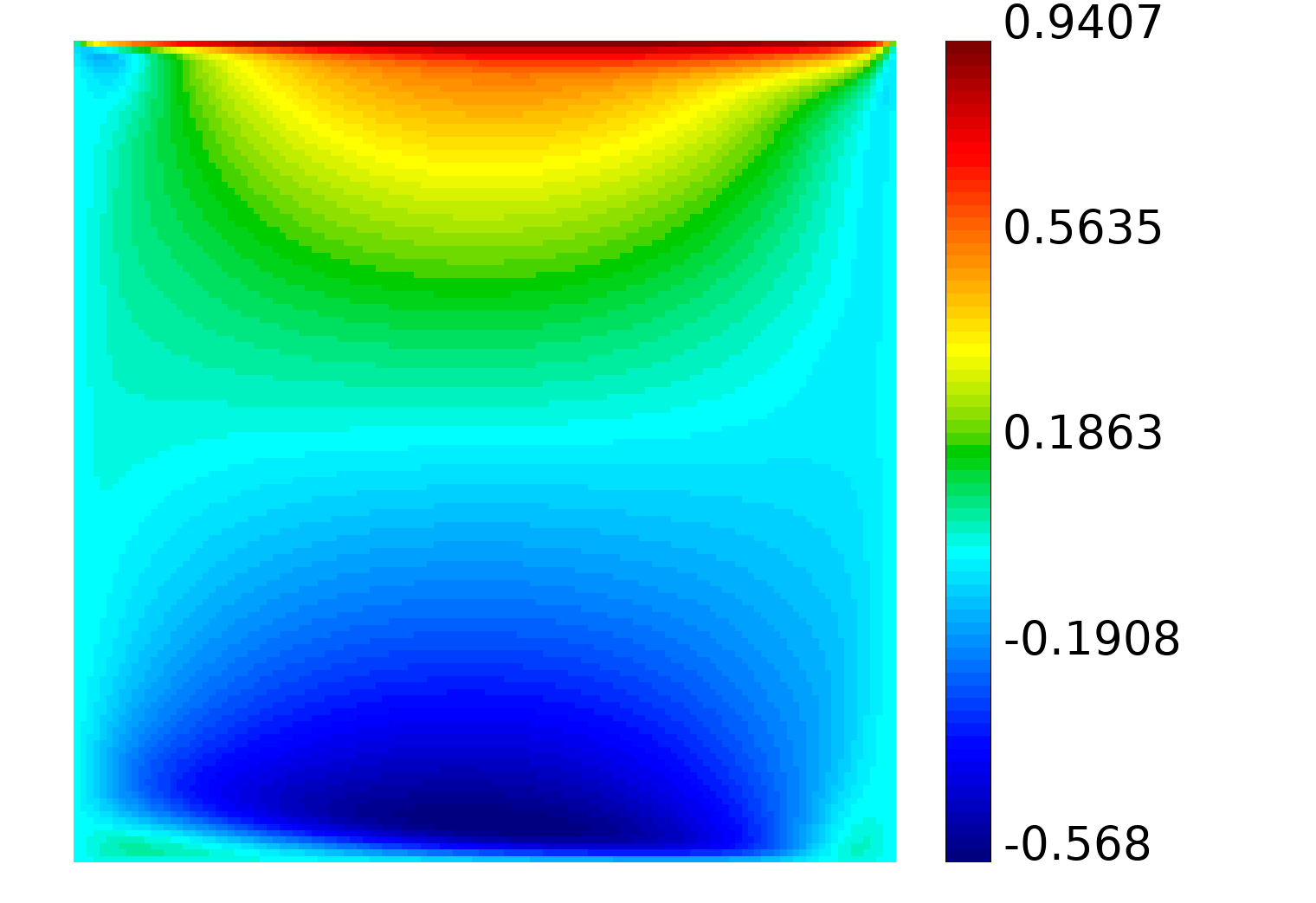}
  \end{subfigure}%
  \begin{subfigure}[t]{.4\textwidth}
    \centering\includegraphics[width=\textwidth]{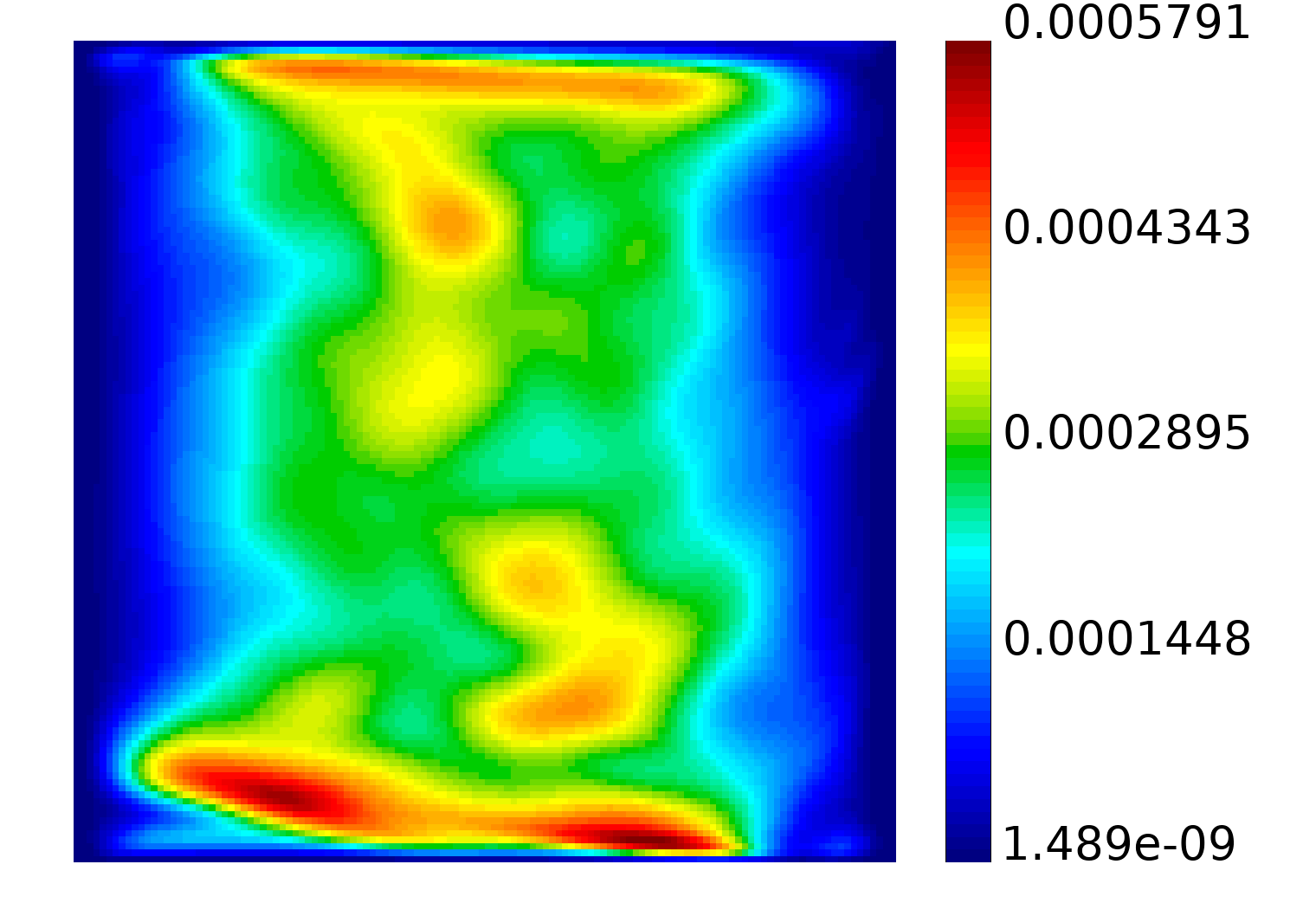}
  \end{subfigure}
  \begin{subfigure}[t]{.4\textwidth}
    \centering\includegraphics[width=\textwidth]{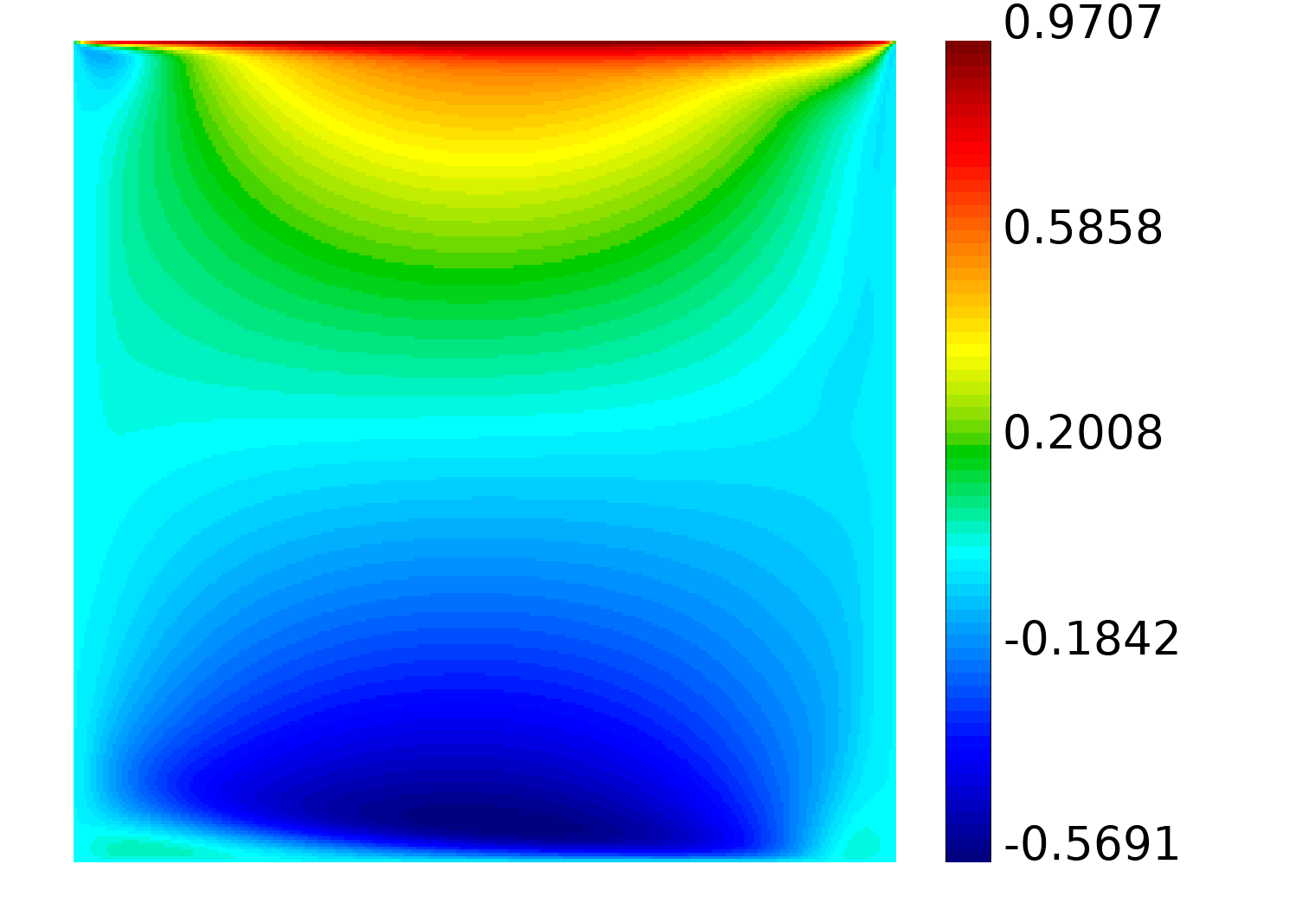}
  \end{subfigure}%
  \begin{subfigure}[t]{.4\textwidth}
    \centering\includegraphics[width=\textwidth]{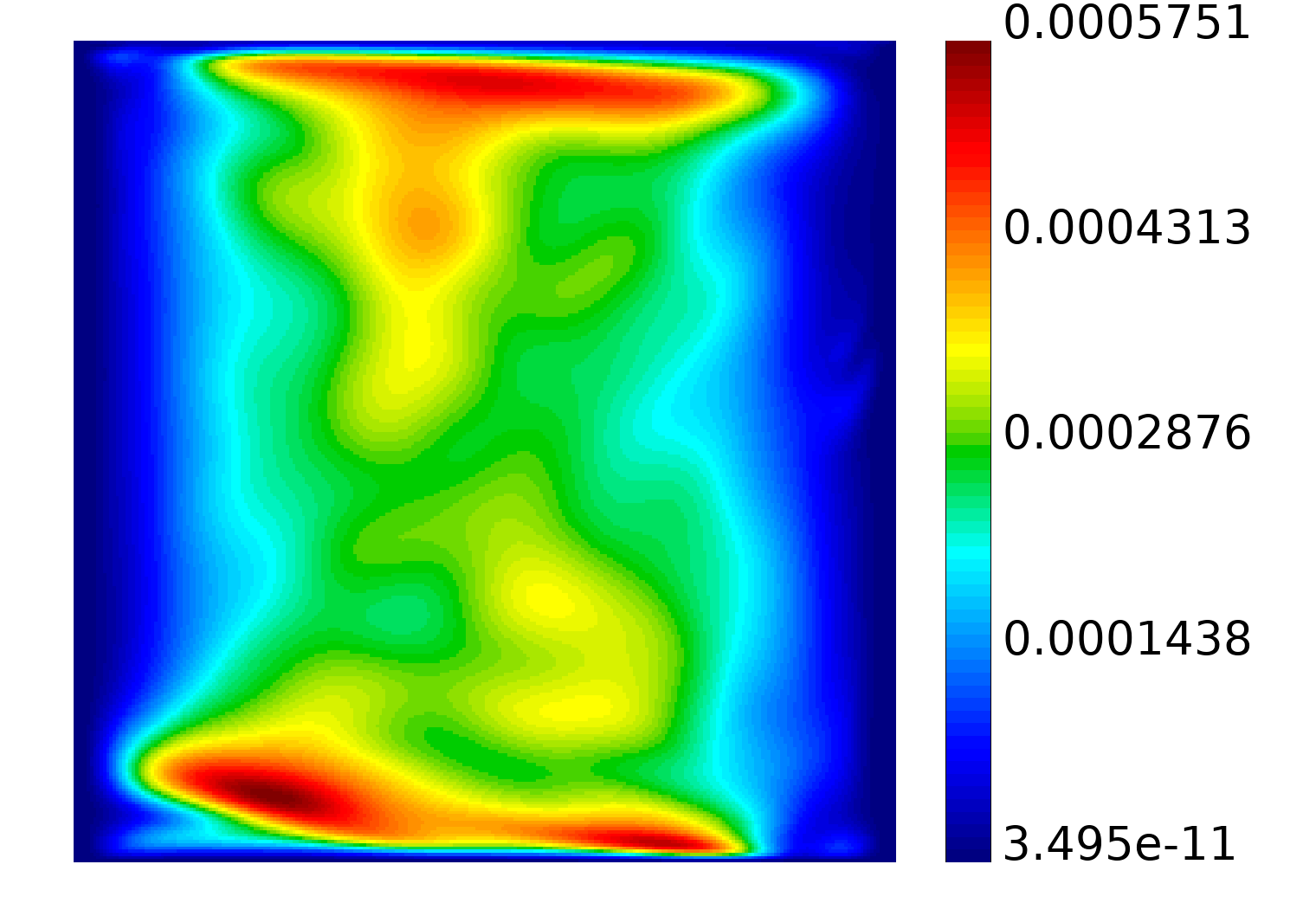}
  \end{subfigure}
  \begin{subfigure}[t]{.4\textwidth}
    \centering\includegraphics[width=\textwidth]{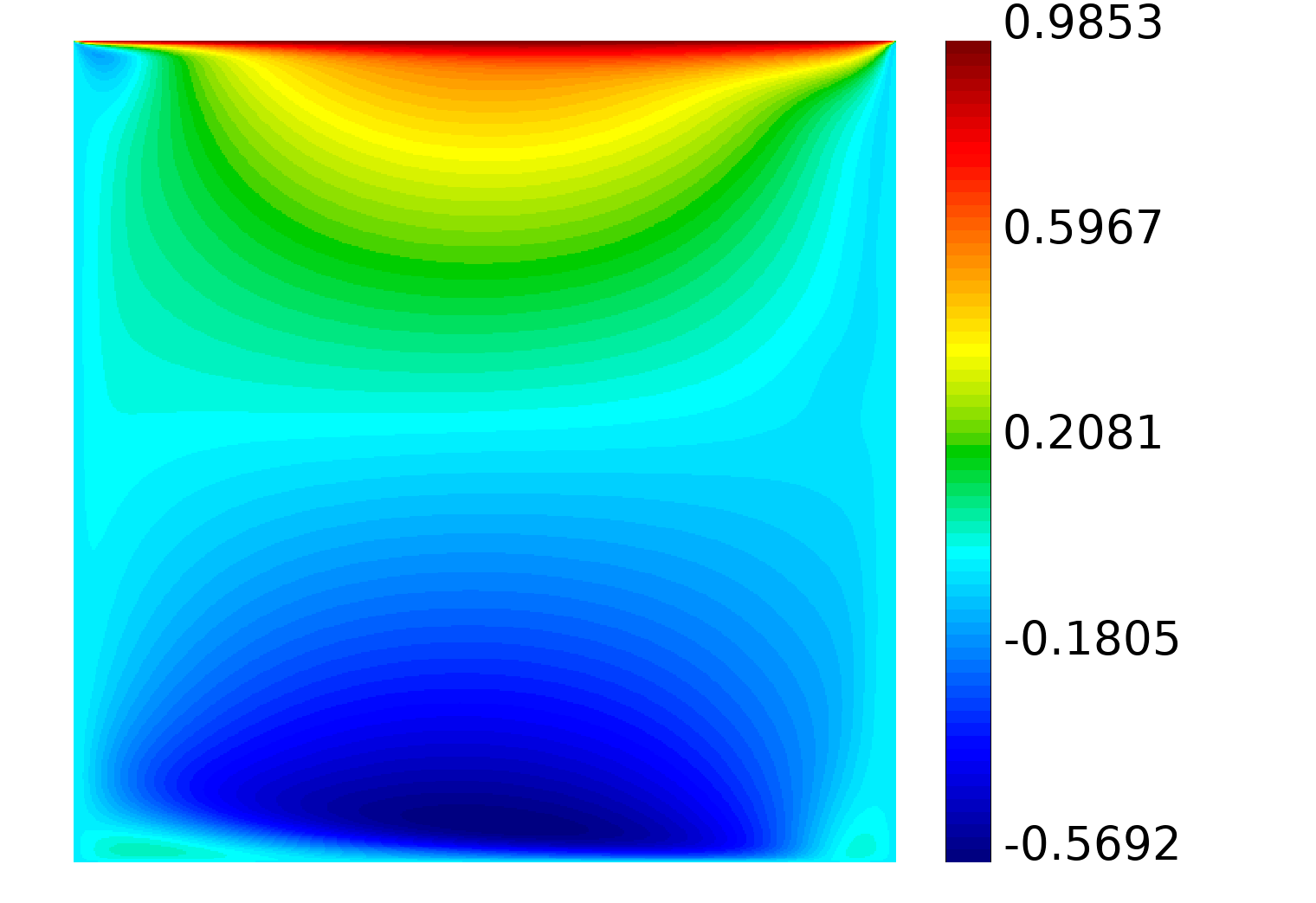}
  \end{subfigure}%
  \begin{subfigure}[t]{.4\textwidth}
    \centering\includegraphics[width=\textwidth]{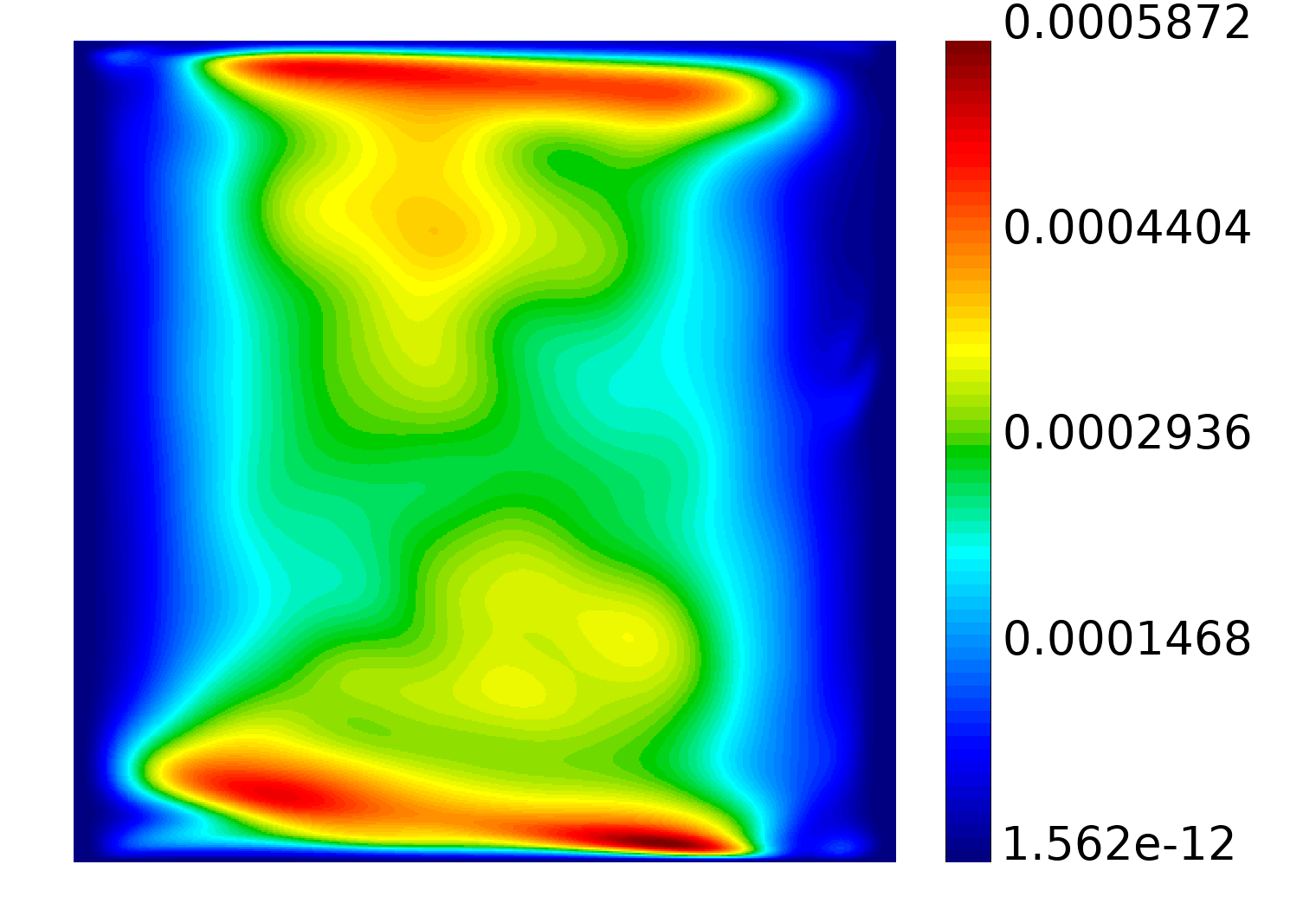}
  \end{subfigure}
  \caption[Lid-driven cavity, statistics of horizontal velocity]{\small {Lid-driven cavity, 
  as described in {\S}\ref{ss:uqIncompNSTest2};
  (left column) mean and (right column) variance of horizontal velocity at $T = 1$.
  From top to bottom, uniform quadrilateral meshes of size $32\!\times\!32$, $64\!\times\!64$, $128\!\times\!128$, $256\!\times\!256$ and $512\!\times\!512$.
  }}
 \label{fig:uqIncompNSTest2VisVx}
\end{figure}

\begin{figure}[!htb]
\centering
  \begin{subfigure}[t]{.4\textwidth}
    \centering\includegraphics[width=\textwidth]{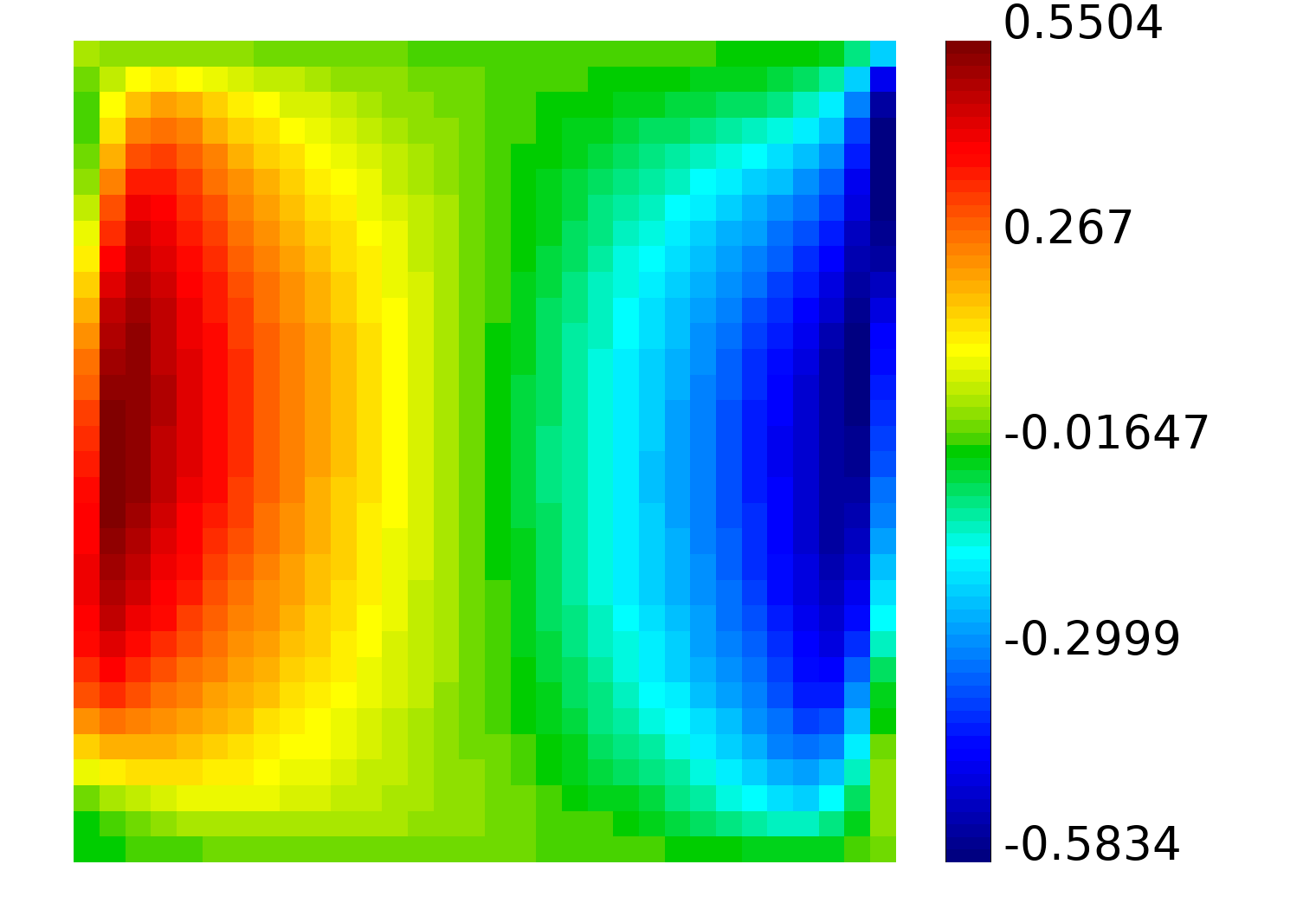}
  \end{subfigure}%
  \begin{subfigure}[t]{.4\textwidth}
    \centering\includegraphics[width=\textwidth]{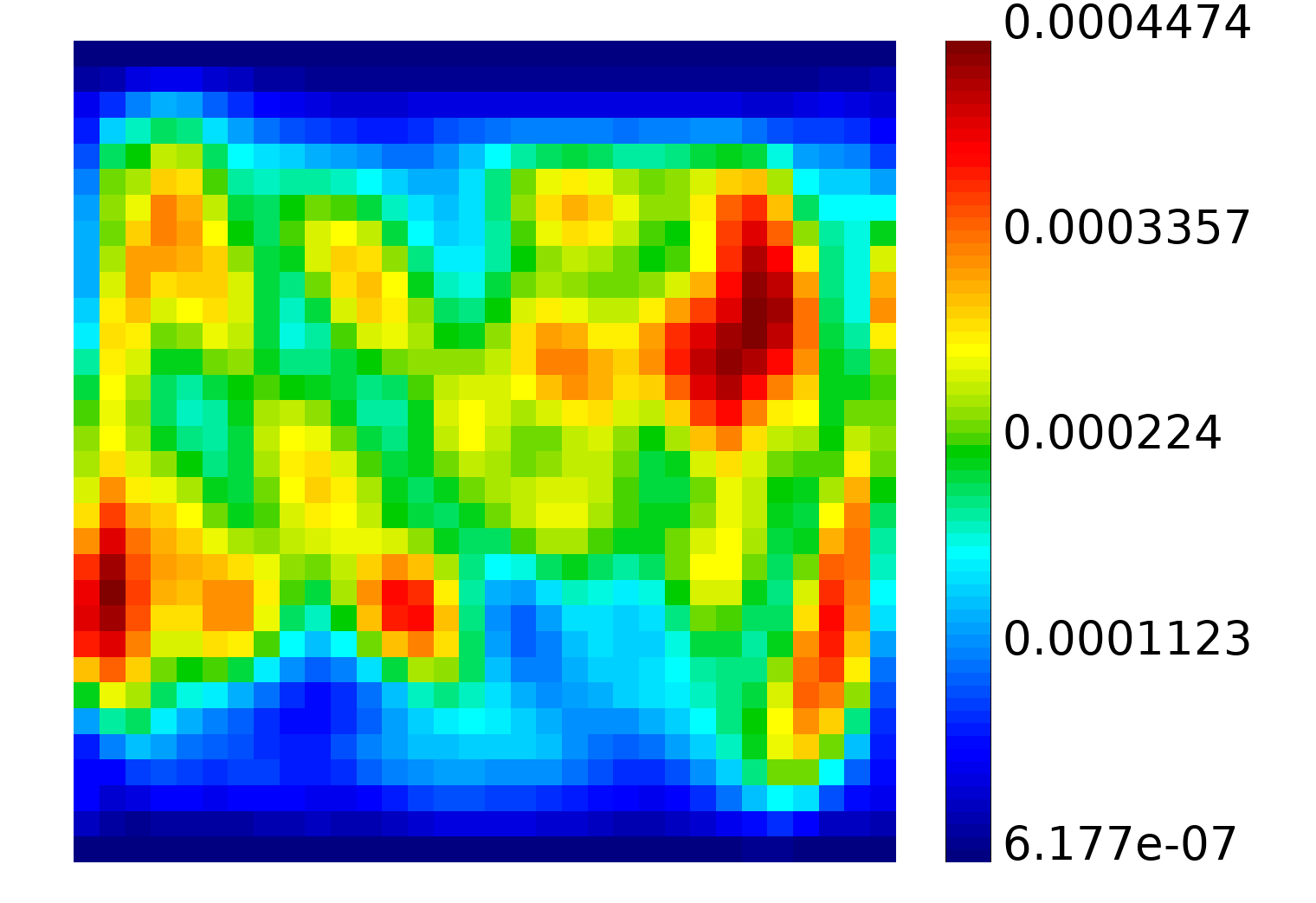}
  \end{subfigure}
  \begin{subfigure}[t]{.4\textwidth}
    \centering\includegraphics[width=\textwidth]{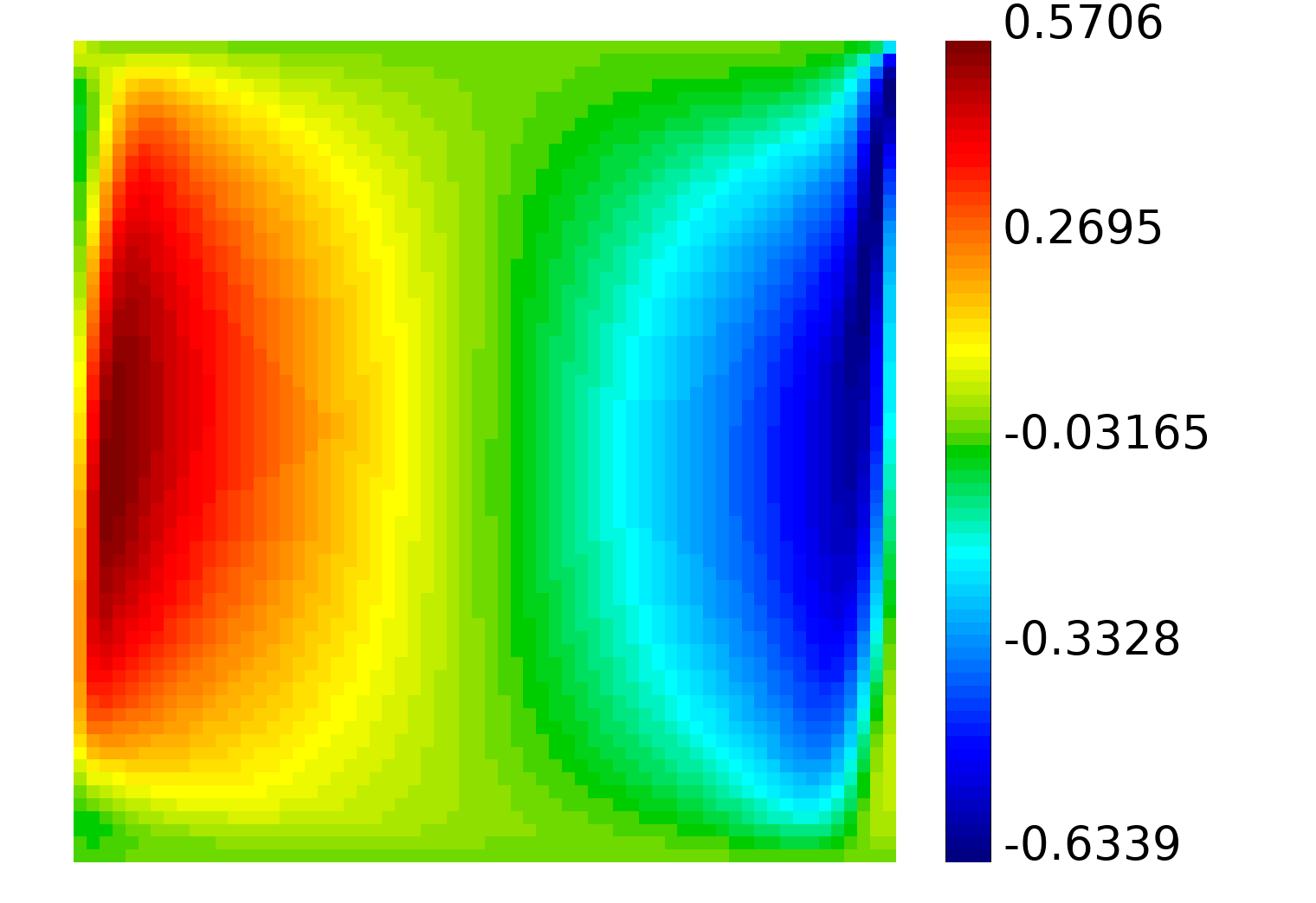}
  \end{subfigure}%
  \begin{subfigure}[t]{.4\textwidth}
    \centering\includegraphics[width=\textwidth]{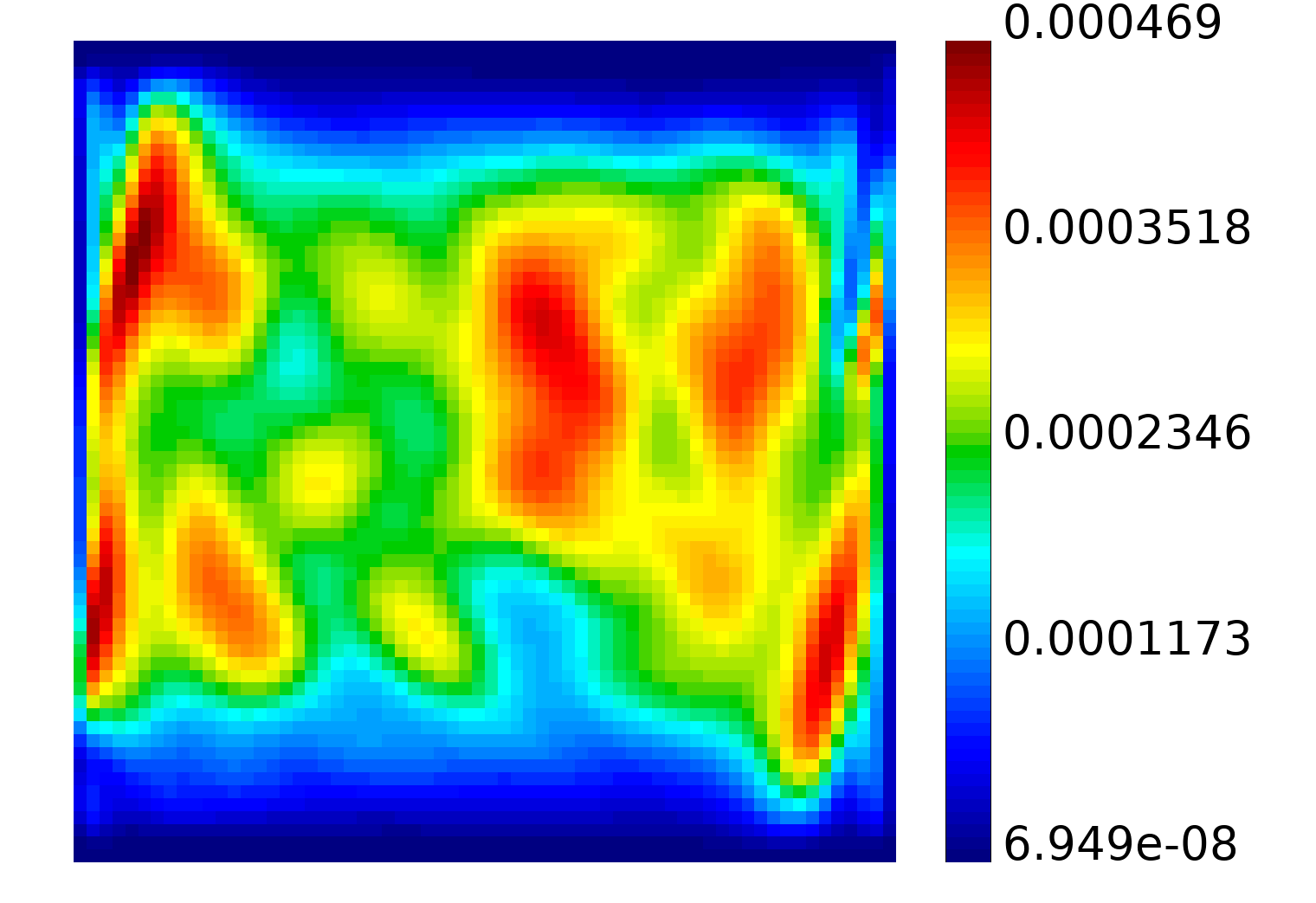}
  \end{subfigure}
  \begin{subfigure}[t]{.4\textwidth}
    \centering\includegraphics[width=\textwidth]{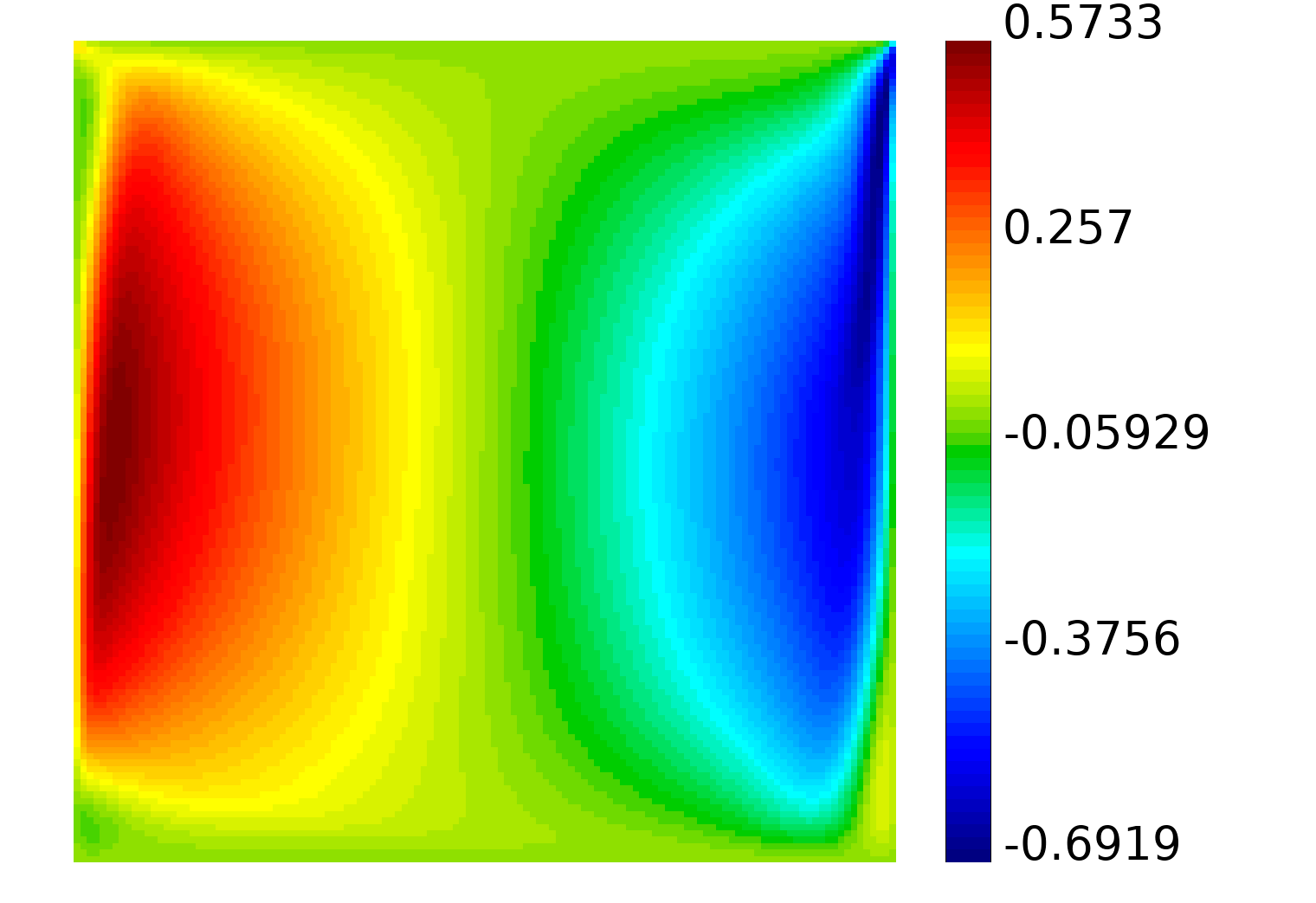}
  \end{subfigure}%
  \begin{subfigure}[t]{.4\textwidth}
    \centering\includegraphics[width=\textwidth]{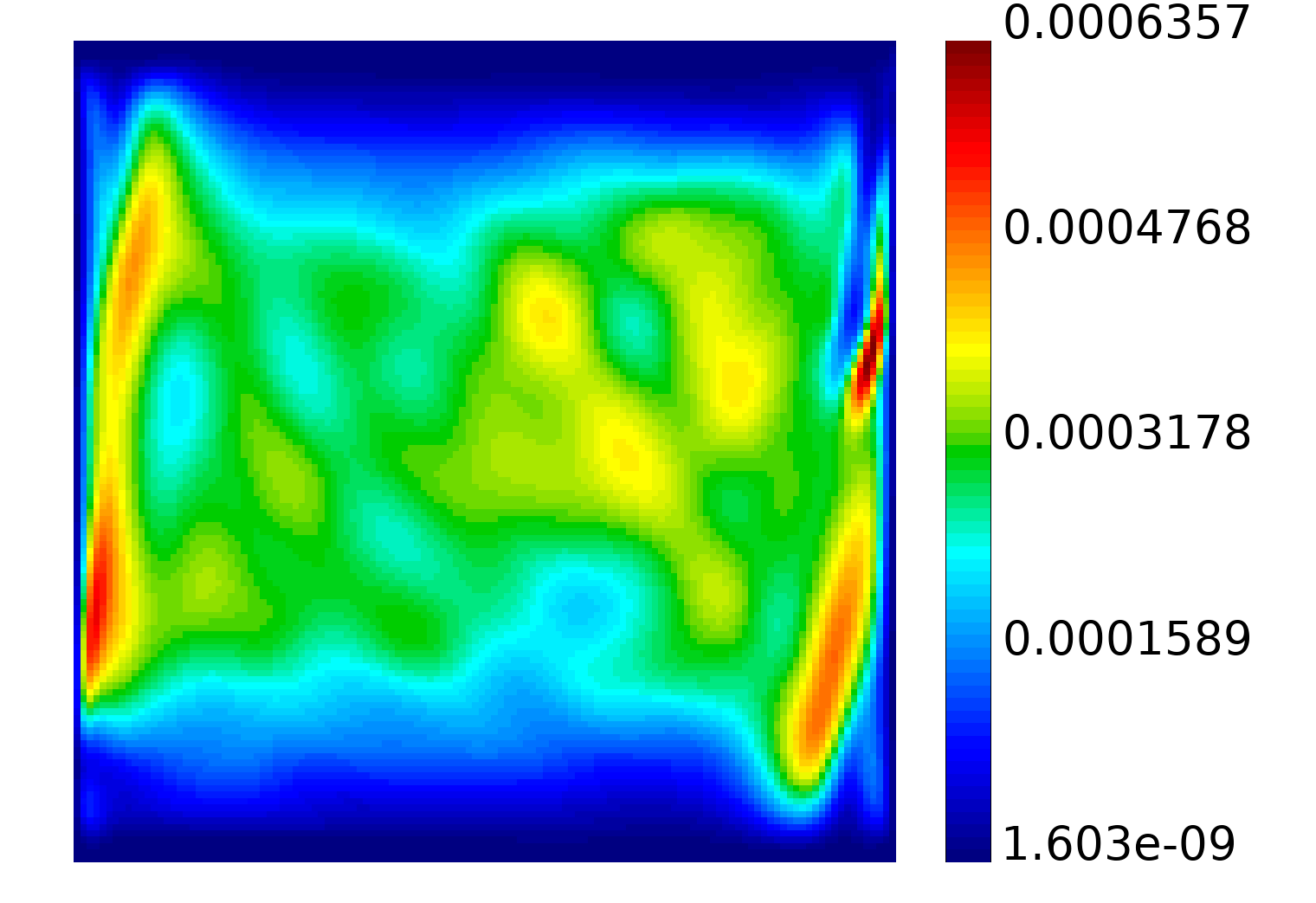}
  \end{subfigure}
  \begin{subfigure}[t]{.4\textwidth}
    \centering\includegraphics[width=\textwidth]{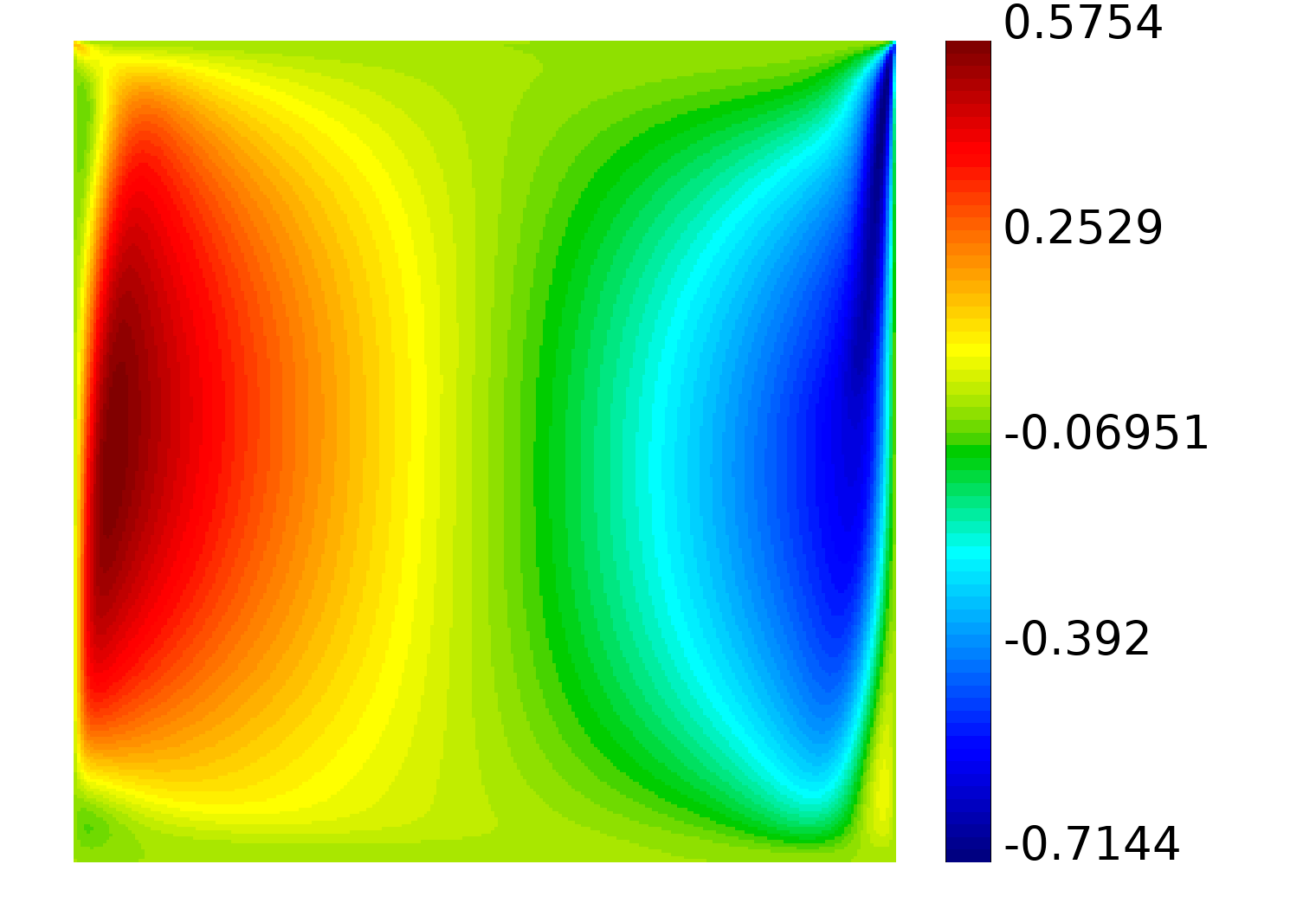}
  \end{subfigure}%
  \begin{subfigure}[t]{.4\textwidth}
    \centering\includegraphics[width=\textwidth]{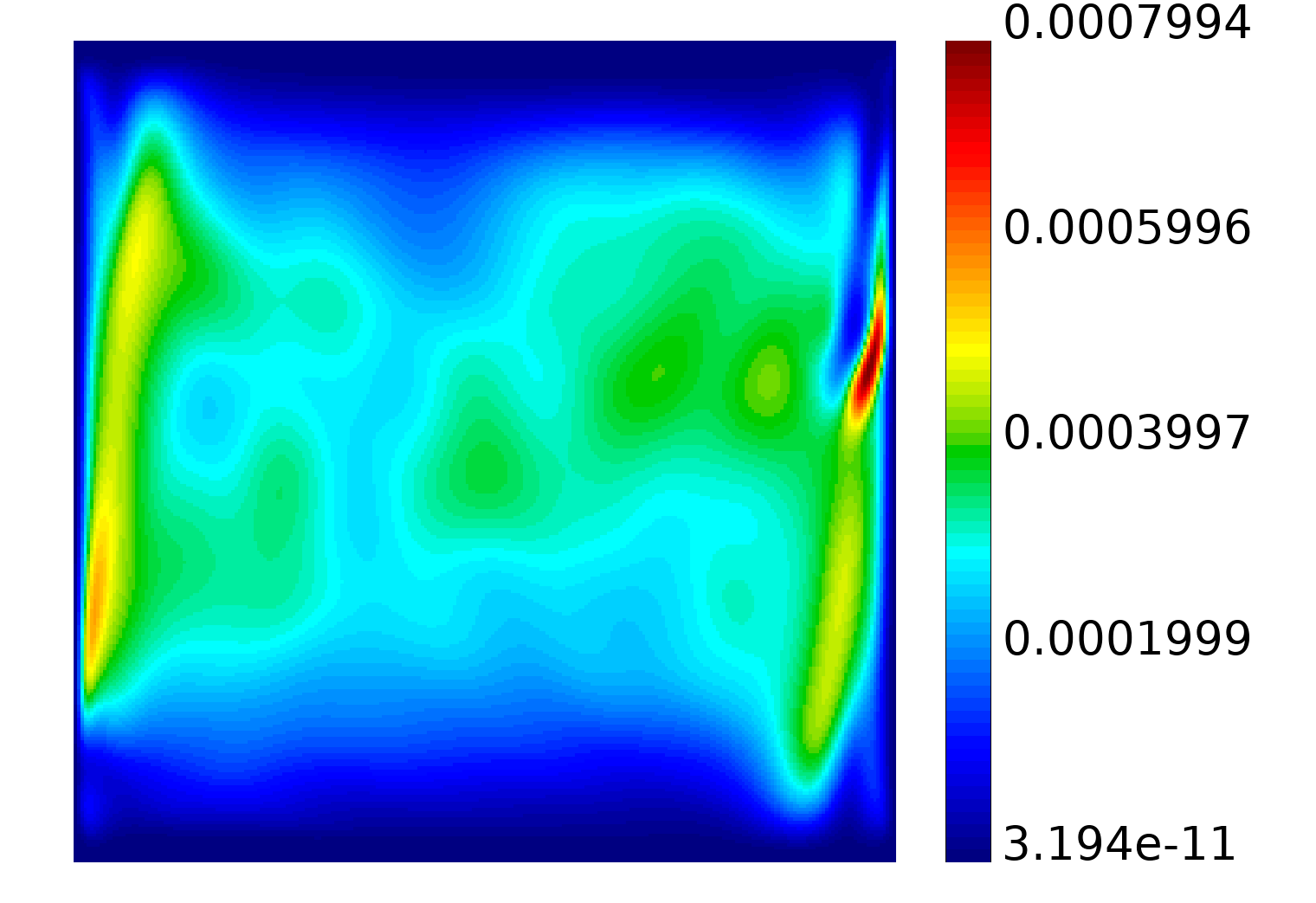}
  \end{subfigure}
  \begin{subfigure}[t]{.4\textwidth}
    \centering\includegraphics[width=\textwidth]{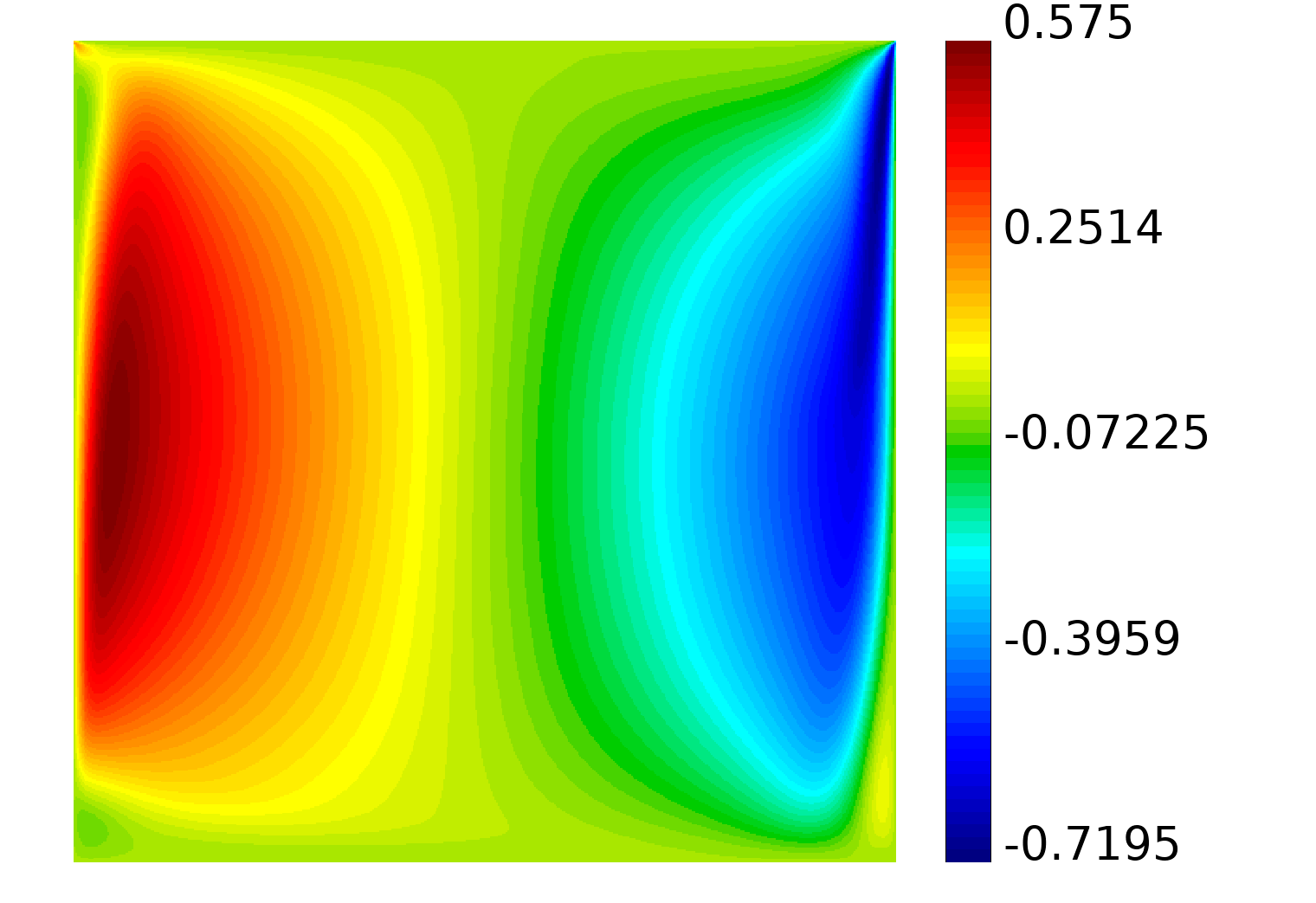}
  \end{subfigure}%
  \begin{subfigure}[t]{.4\textwidth}
    \centering\includegraphics[width=\textwidth]{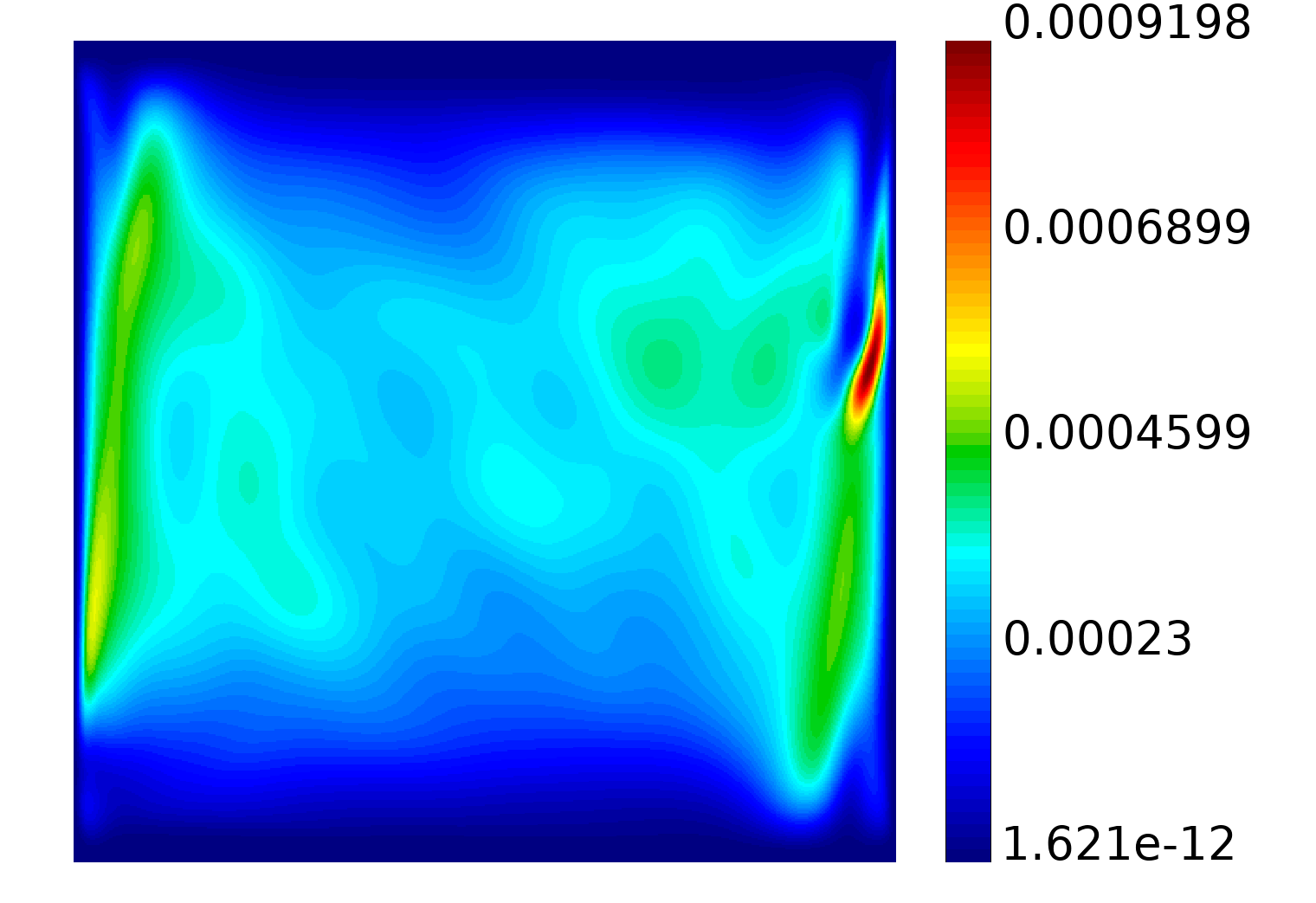}
  \end{subfigure}
  \caption[Lid-driven cavity, statistics of vertical velocity]{\small {Lid-driven cavity, 
  as described in {\S}\ref{ss:uqIncompNSTest2};
  (left column) mean and (right column) variance of vertical velocity at $T = 1$.
  From top to bottom, uniform quadrilateral meshes of size $32\!\times\!32$, $64\!\times\!64$, $128\!\times\!128$, $256\!\times\!256$ and $512\!\times\!512$.
  }}
 \label{fig:uqIncompNSTest2VisVy}
\end{figure}

\begin{figure}[!htb]
\begin{center}
 \includegraphics[width = 0.5\textwidth]{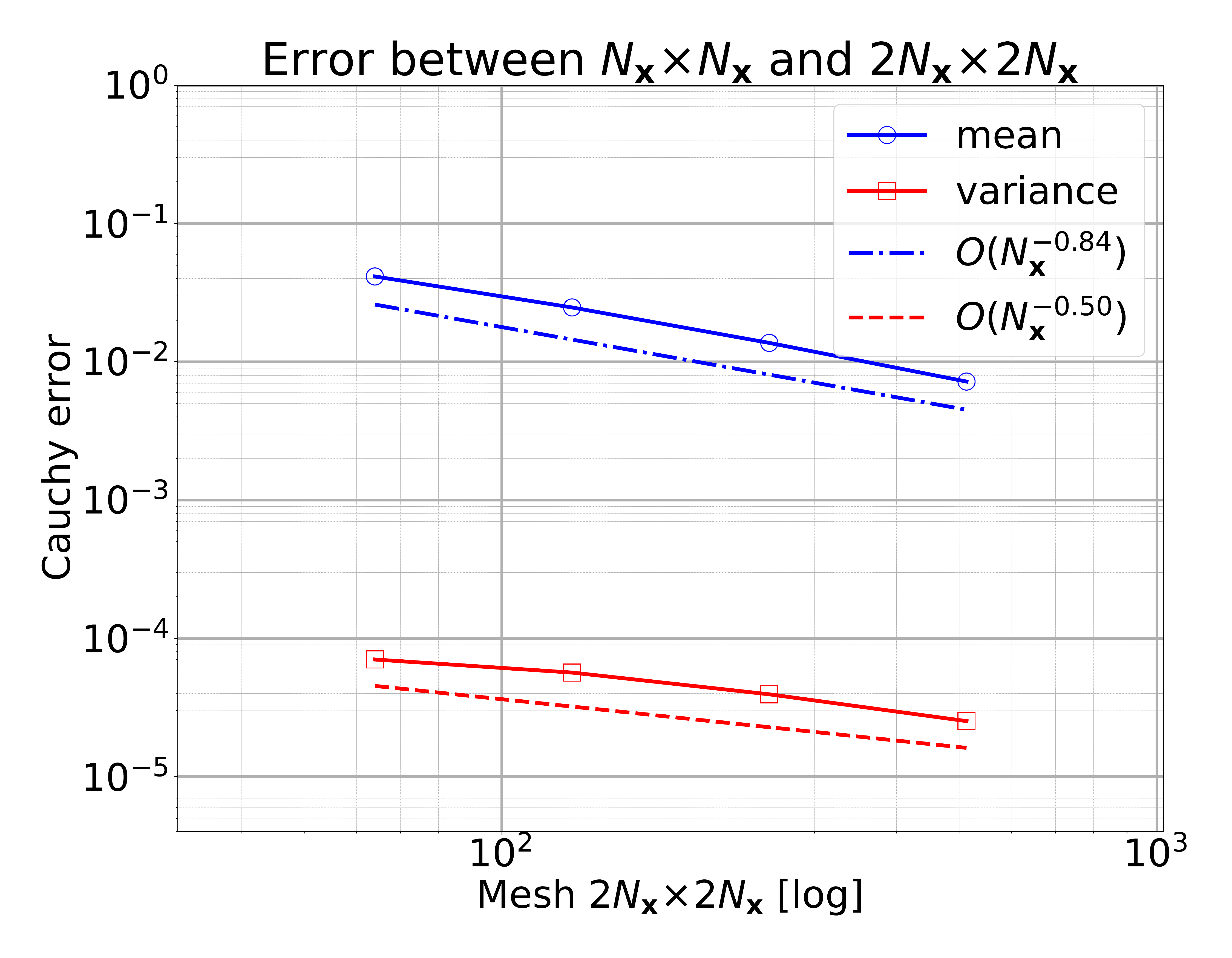}
 \caption[Lid-driven cavity, convergence of mean and variance]{\small {Lid-driven cavity, $Re=3200$, 
  as described in {\S}\ref{ss:uqIncompNSTest2};
  Cauchy error of mean and variance measured in $\ltwo{\Dx}$ norm between resolutions from $32\!\times\!32$ to $512\!\times\!512$ with uniform quadrilateral meshes, at time $T=1$.
 }}
 \label{fig:uqIncompNSTest2Convg}
\end{center}
\end{figure}

\begin{figure}[!htb]
\centering
  \begin{subfigure}[t]{.5\textwidth}
    \includegraphics[width = \textwidth]{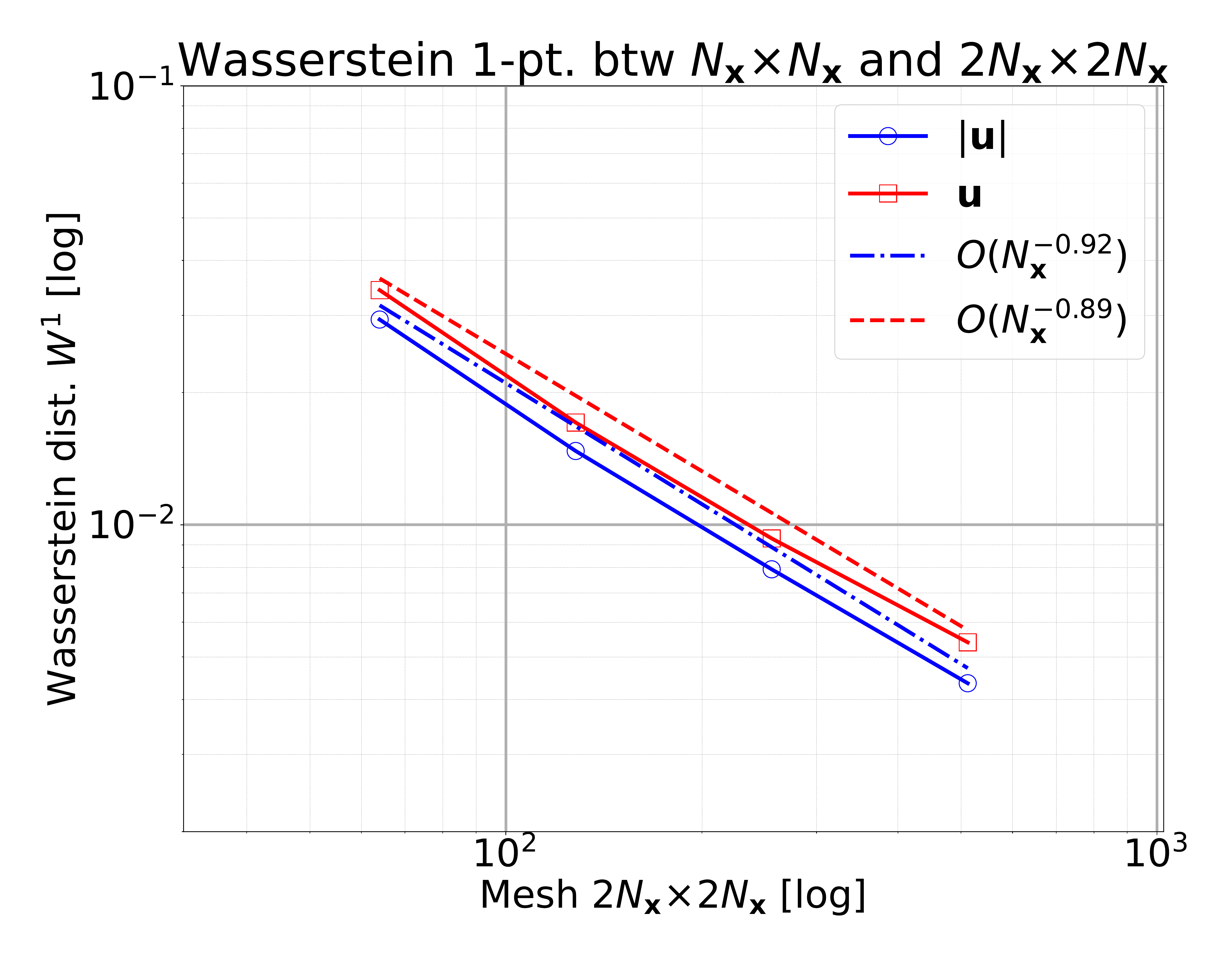}
  \end{subfigure}%
  \begin{subfigure}[t]{.5\textwidth}
    \includegraphics[width = \textwidth]{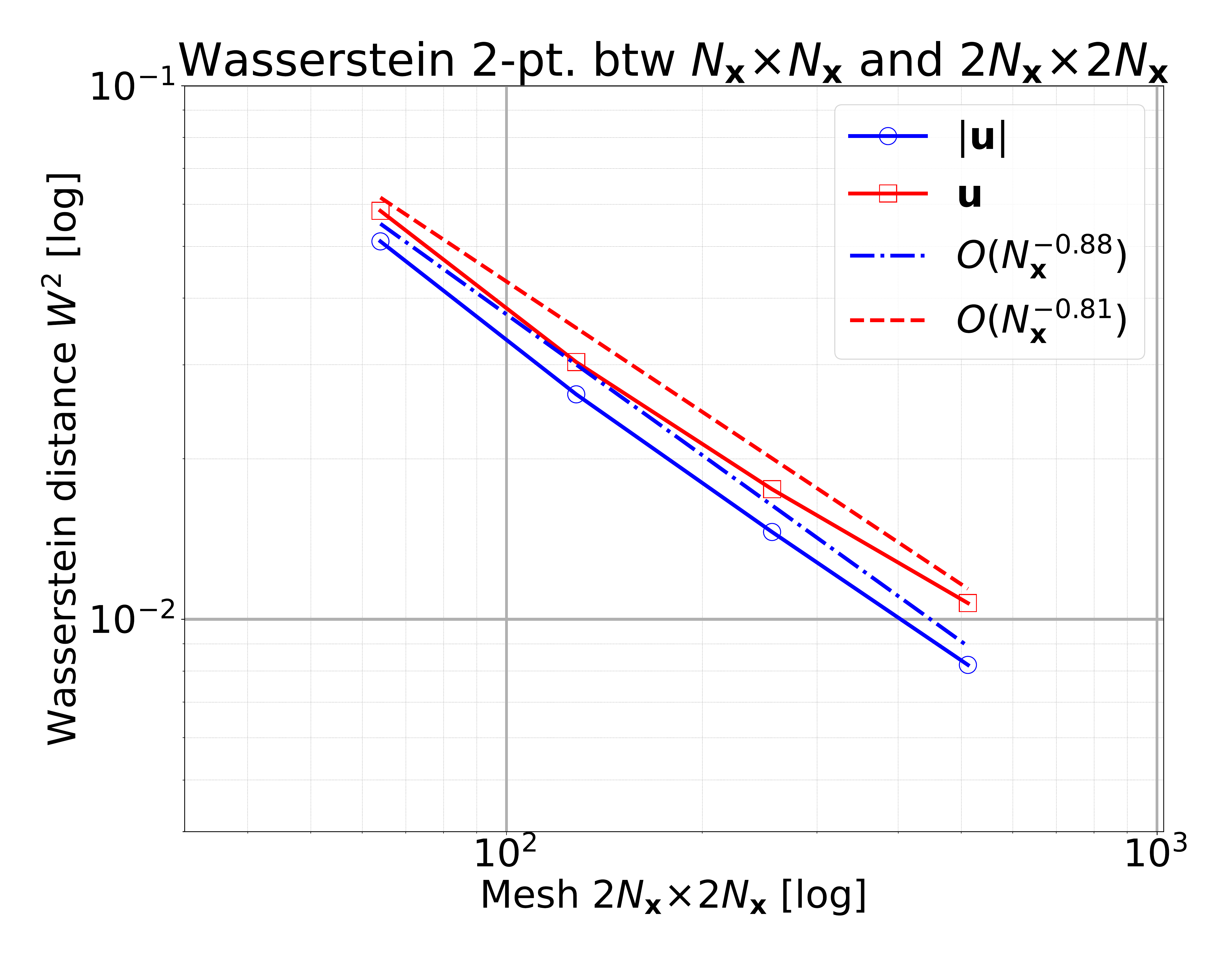}
  \end{subfigure}
 \caption[Lid-driven cavity, convergence of Wasserstein distances]{\small {Lid-driven cavity, 
  as described in {\S}\ref{ss:uqIncompNSTest2};
  Wasserstein distances of velocity field $\bu$ and its magnitude $\abs{\bu}$ measured between successive resolutions from $32\!\times\!32$ to $512\!\times\!512$ with uniform quadrilateral meshes, at time $T=1$.
  }}
 \label{fig:uqIncompNSTest2Wd}
\end{figure}

\begin{figure}[!htb]
\centering
  \begin{subfigure}[t]{.5\textwidth}
    \includegraphics[width = \textwidth]{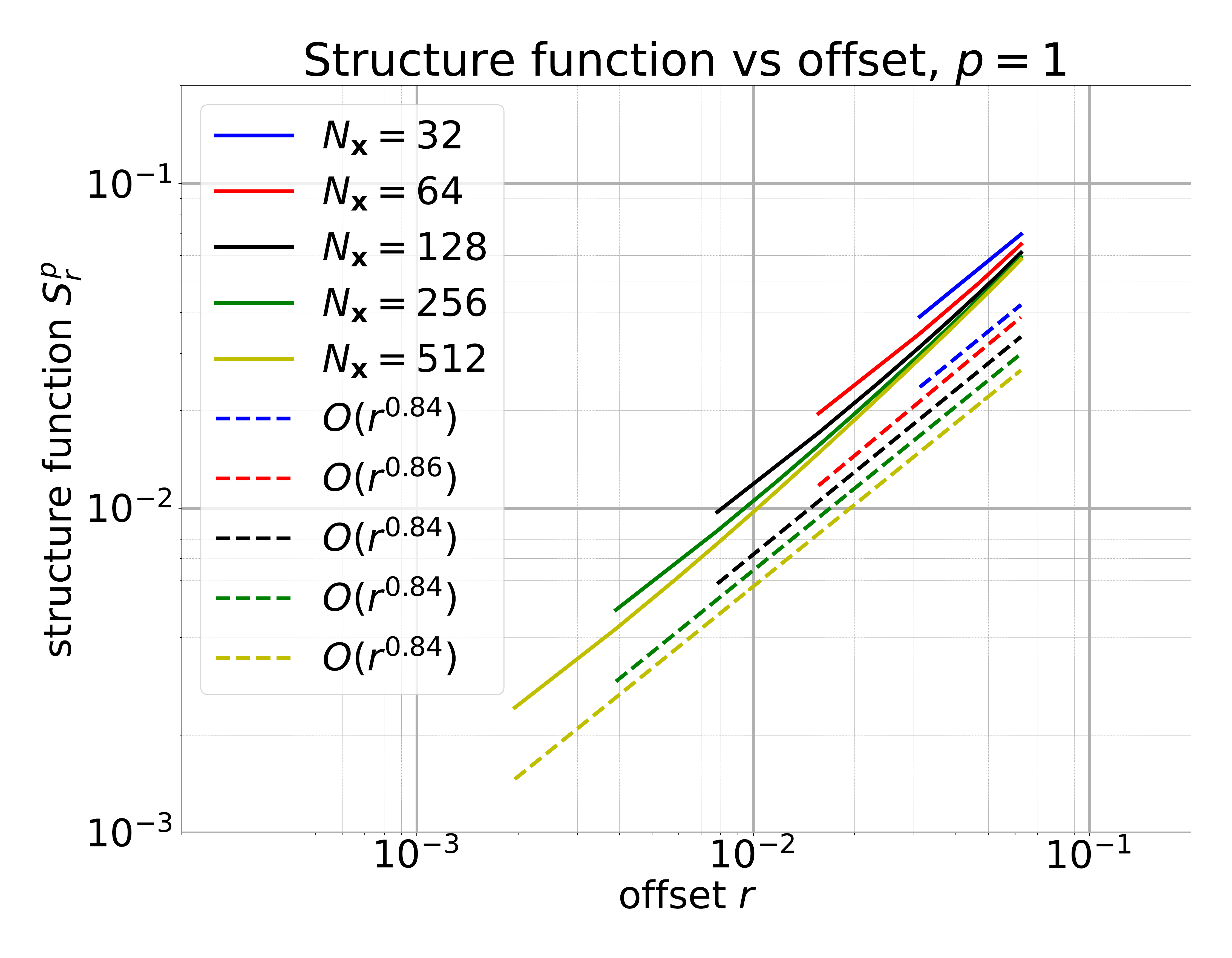}
  \end{subfigure}
  \begin{subfigure}[t]{.5\textwidth}
    \includegraphics[width = \textwidth]{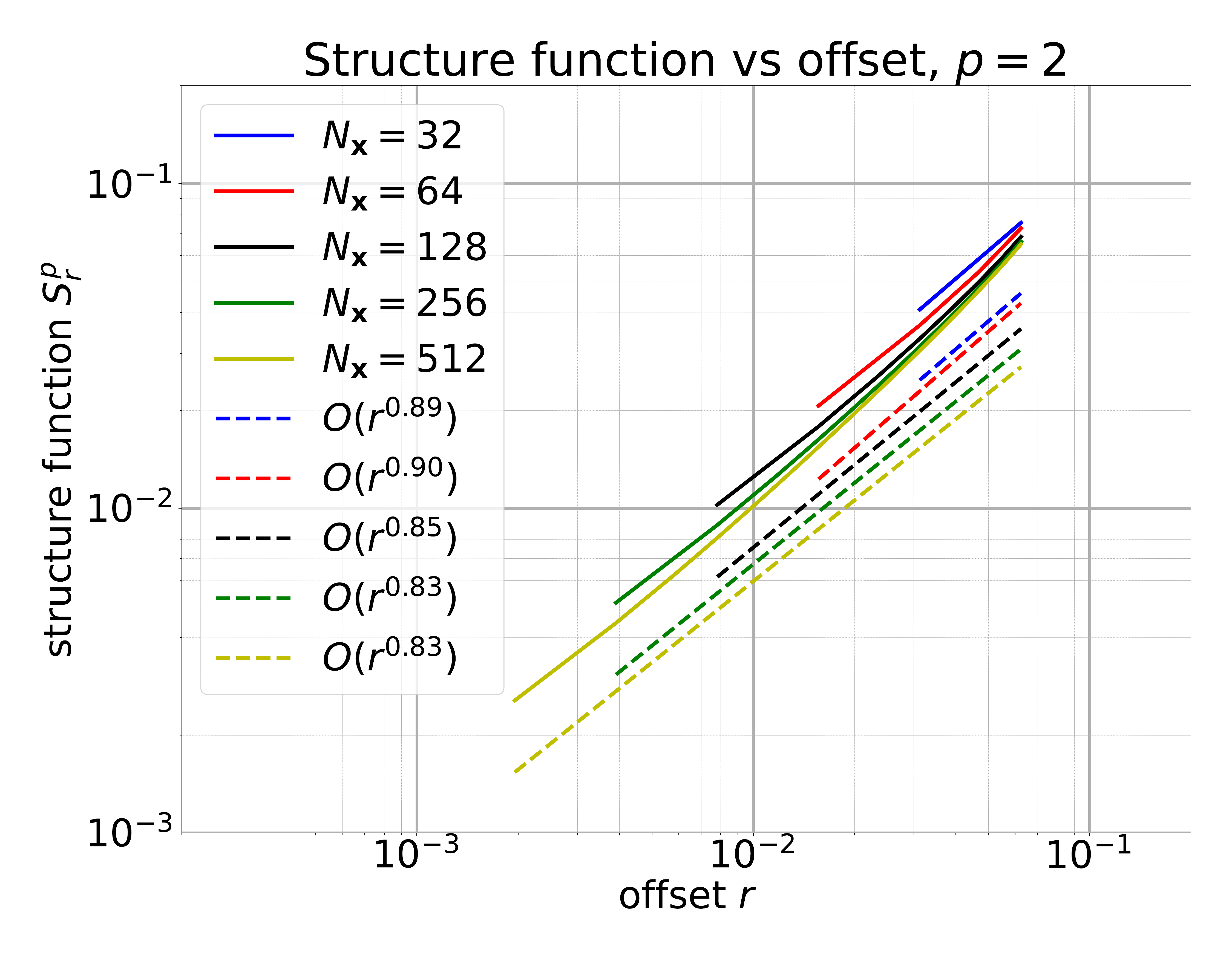}
  \end{subfigure}
  \begin{subfigure}[t]{.5\textwidth}
    \includegraphics[width = \textwidth]{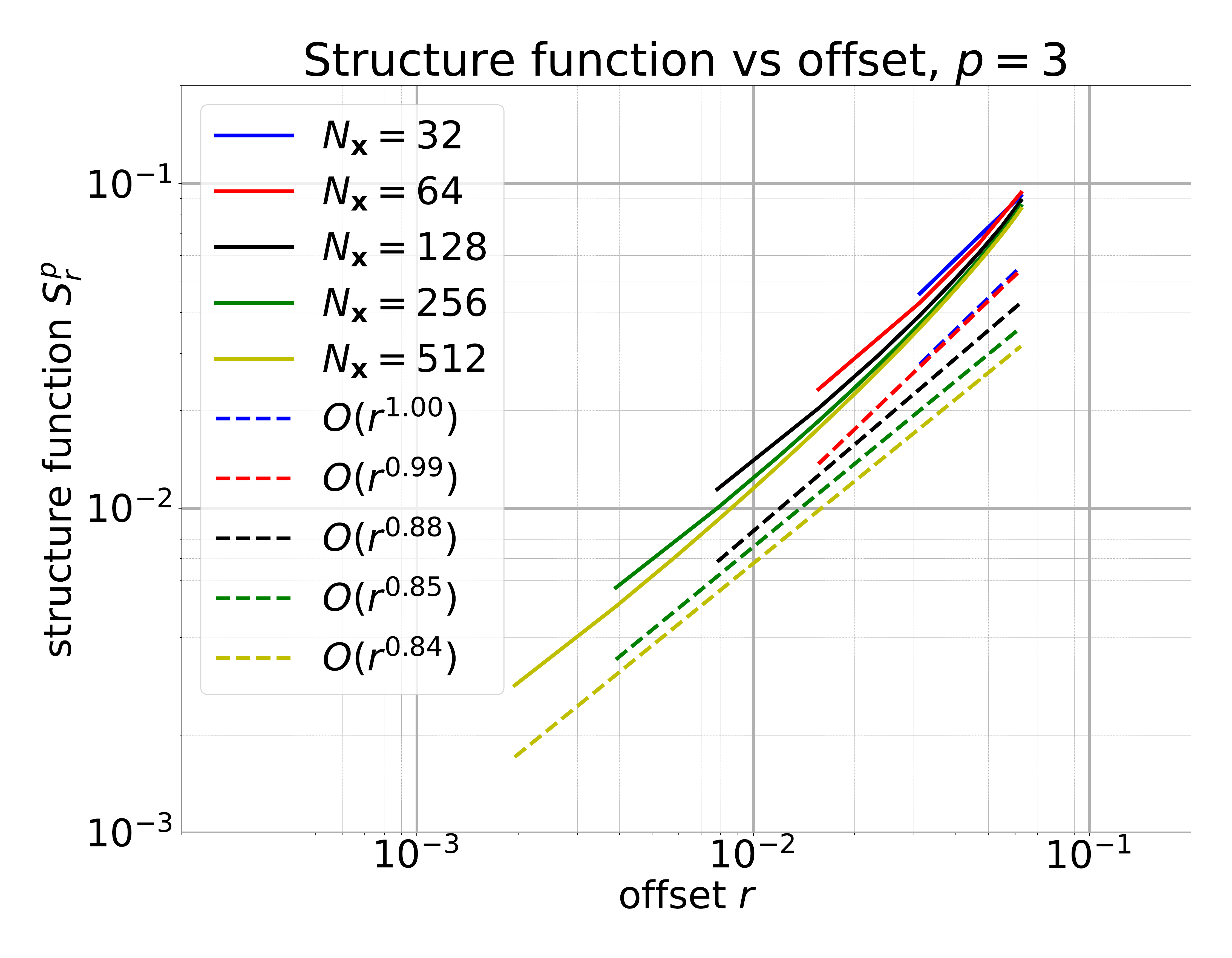}
  \end{subfigure}
 \caption[Lid-driven cavity, convergence of structure function]{\small {Lid-driven cavity, 
  as described in {\S}\ref{ss:uqIncompNSTest2};
  Approximation of the structure function of the velocity field on $N_{\bx}\!\times\!N_{\bx}$ uniform quadrilateral meshes, at time $T=1$.
 }}
 \label{fig:uqIncompNSTest2Scube}
\end{figure}

The observations are as follows:
\begin{itemize}[nolistsep]
 \item We visualize the mean and variance of the velocity solutions in figures~\ref{fig:uqIncompNSTest2VisVx}, \ref{fig:uqIncompNSTest2VisVy} and it seems that these statistics converge as the mesh resolution is increased. Furthermore, from Figure~\ref{fig:uqIncompNSTest2Convg}, it can be observed that the Cauchy error of mean and variance decays with mesh refinement.
 \item We measure Wasserstein distances $W^{1}$ and $W^{2}$ (see \eqref{eq:WassersteinDistW12}) between two ensembles $\bb{U}^{h,M}_{T}$ and $\bb{U}^{h/2,2M}_{T}$. From Figure~\ref{fig:uqIncompNSTest2Wd}, we observe that these distances decrease as we refine the mesh, which indicates convergence.
 \item In this experiment, the solution is smooth, thus, Lipschitz continuous and we expect the structure functions to converge with respect to the offset $r$ at a rate close to $1$  (cf. Remark~\ref{rem:structureFnScaling}). 
 From Figure~\ref{fig:uqIncompNSTest2Scube}, we observe that structure functions converge with a rate $0.85$ (approx.).
\end{itemize}
visualize vx
%

%
\subsection{Channel flow}\label{ss:uqIncompNSTest3}
%
We choose the spatial domain $\Dx := (0,3L)\times(0,L)$ with $L=0.5$, the time horizon $T = 0.8$ and the Reynolds number $Re = L\nu^{-1} \in \{1600, 3200\}$.
For $K \in \N$ such that $K$ {is an even integer}, we consider the {random} perturbation function
\begin{equation*}
f(\omega; \bx) = \sum_{k=0}^{K/2} Y_{2k}(\omega) \sin(2\pi k(x_2 + Y_{2k+1}(\omega))), 
\ \forall \bx\in\Dx,
\end{equation*}
where $Y_{j} \sim \mathcal{U}[-1,1]$ are independent, identically distributed random variables.
The initial velocity $\bu_{0}$ is given by
\begin{equation*}
\label{eq:uqIncompNSTest3InitSol}
 \begin{aligned}
  {u}_{0;1}(\omega; \bx) &= (1 + \gamma_1 f(\omega; \bx))\ 4 u_{\max} \frac{x_2(L - x_2)}{L^2},\\
  {u}_{0;2}(\omega; \bx) &= \gamma_2 f(\omega; \bx) \frac{x_2(L - x_2)}{L^2},
 \end{aligned}
\end{equation*}
with $u_{\max} = 1.5$, $\gamma_1 = \gamma_2 = 0.025$ and $K = 10$.
We fix the source {term} $\bb{f} = \bb{0}$. 
The bottom and top boundaries are fixed.
We impose inflow conditions, same as the initial velocity, at the left boundary
and outflow conditions on the right boundary. Note that $\bu_{0}(\omega) \in H$.

\begin{figure}[htpb]
\centering
  \begin{subfigure}[ht]{.5\textwidth}
    \centering\includegraphics[width=0.9\textwidth]{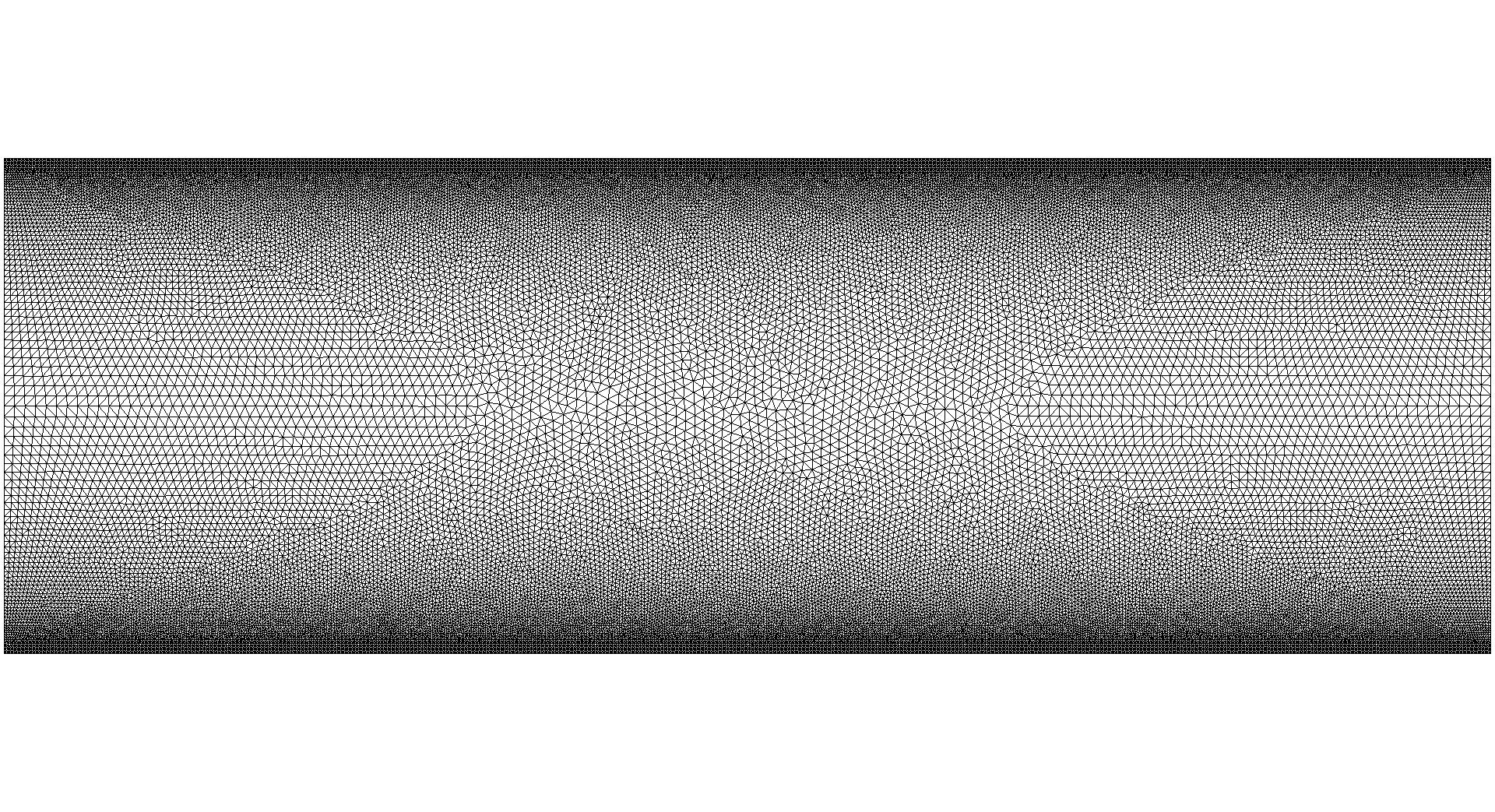}
  \end{subfigure}%
  \begin{subfigure}[ht]{.5\textwidth}
    \centering\includegraphics[width=0.9\textwidth]{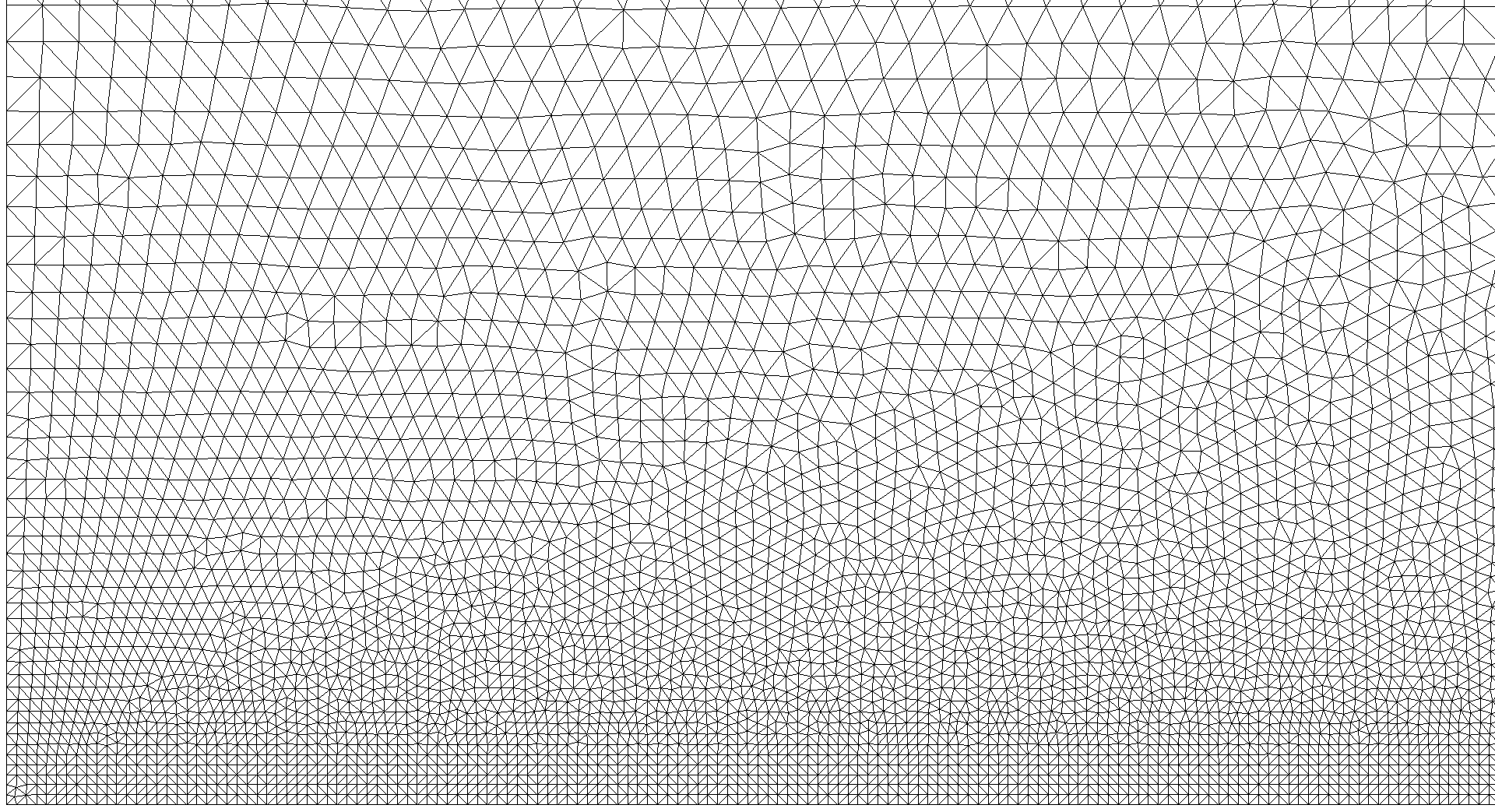}
  \end{subfigure}
  \caption[An example of the unstructured mesh used for simulating channel flow]{\small {(left) Unstructured mesh for the spatial domain $\Dx := (0,3L)\times(0,L)$ with $L=0.5$ that corresponds to level $\ellx = 1$, (right) zoomed-in view of the bottom-left region of the mesh. The minimum and maximum mesh element sizes are approximately $0.0015$ and $0.013$, respectively.}}
 \label{fig:uqIncompNSTest3Mesh}
\end{figure}
\begin{table}[htpb]
\centering
 \begin{tabular}{c | c}
 Mesh level $\ellx$ & Number of time steps\\
 \hline
 0 & 400  \\
 1 & 800 \\
 2 & 1600 \\
 3 & 2500 \\
 \end{tabular}
 \caption{\small {
 Number of time steps used in the numerical solver for different mesh resolutions in the channel flow problem described in {\S}\ref{ss:uqIncompNSTest3}.
 }}
 \label{tab:uqIncompNSTest3NSteps}
\end{table}
We generate unstructured meshes corresponding to resolution levels $\ellx = 0, 1$ using the Gmsh library \cite{Gmsh}, the mesh for level $\ellx = 1$ is shown in Figure \ref{fig:uqIncompNSTest3Mesh}. Additional meshes are generated by conforming uniform refinement of the mesh with resolution level $\ellx=1$.
The number of time steps used in the solver corresponding to different mesh levels are given in Table~\ref{tab:uqIncompNSTest3NSteps}.
A visualization of the initial conditions is shown in the figure below.
%
\begin{figure}[!htb]
\centering
  \begin{subfigure}[ht]{.45\textwidth}
    \centering\includegraphics[width=\textwidth]{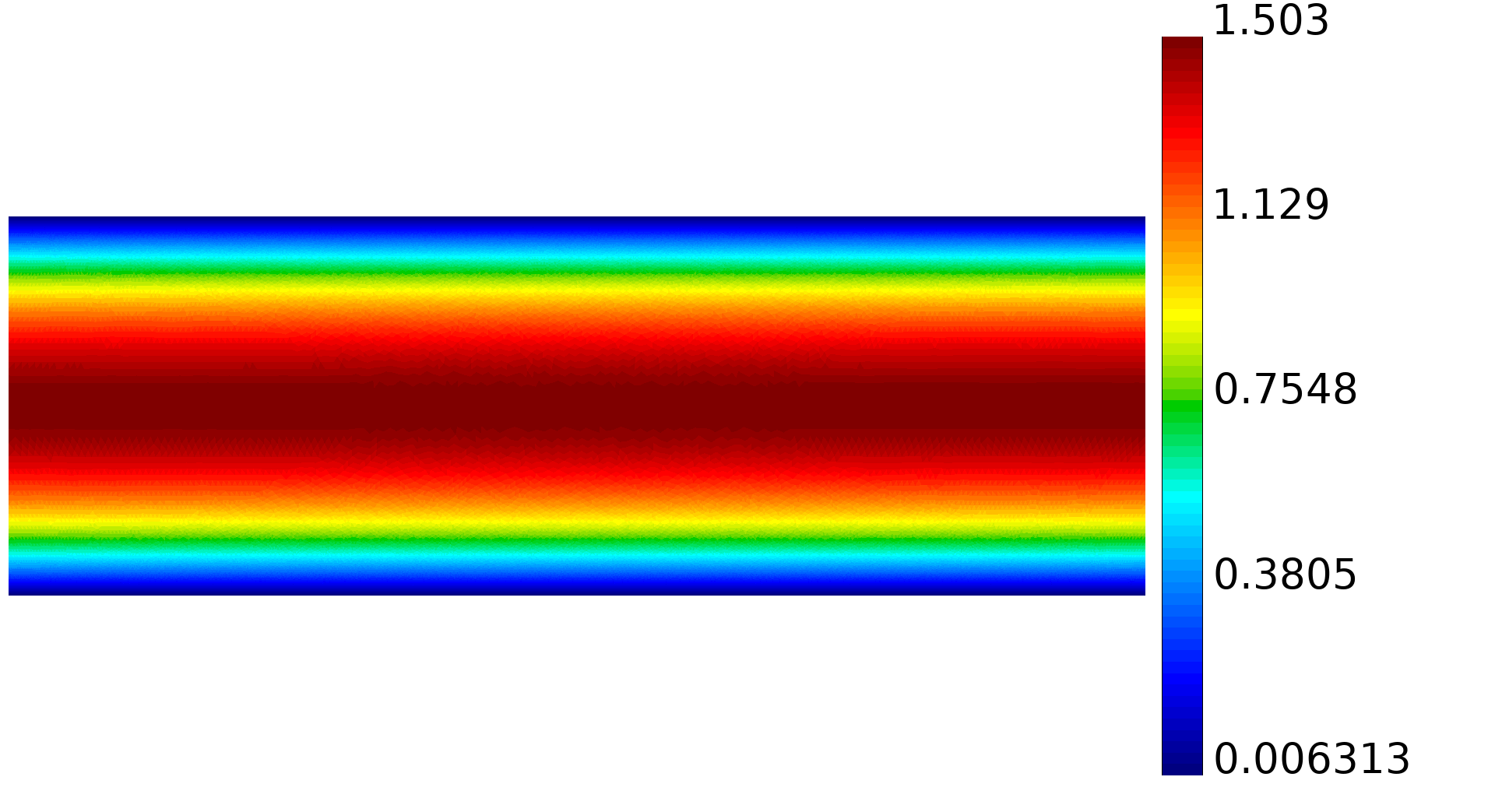}
  \end{subfigure}%
  \begin{subfigure}[ht]{.45\textwidth}
    \centering\includegraphics[width=\textwidth]{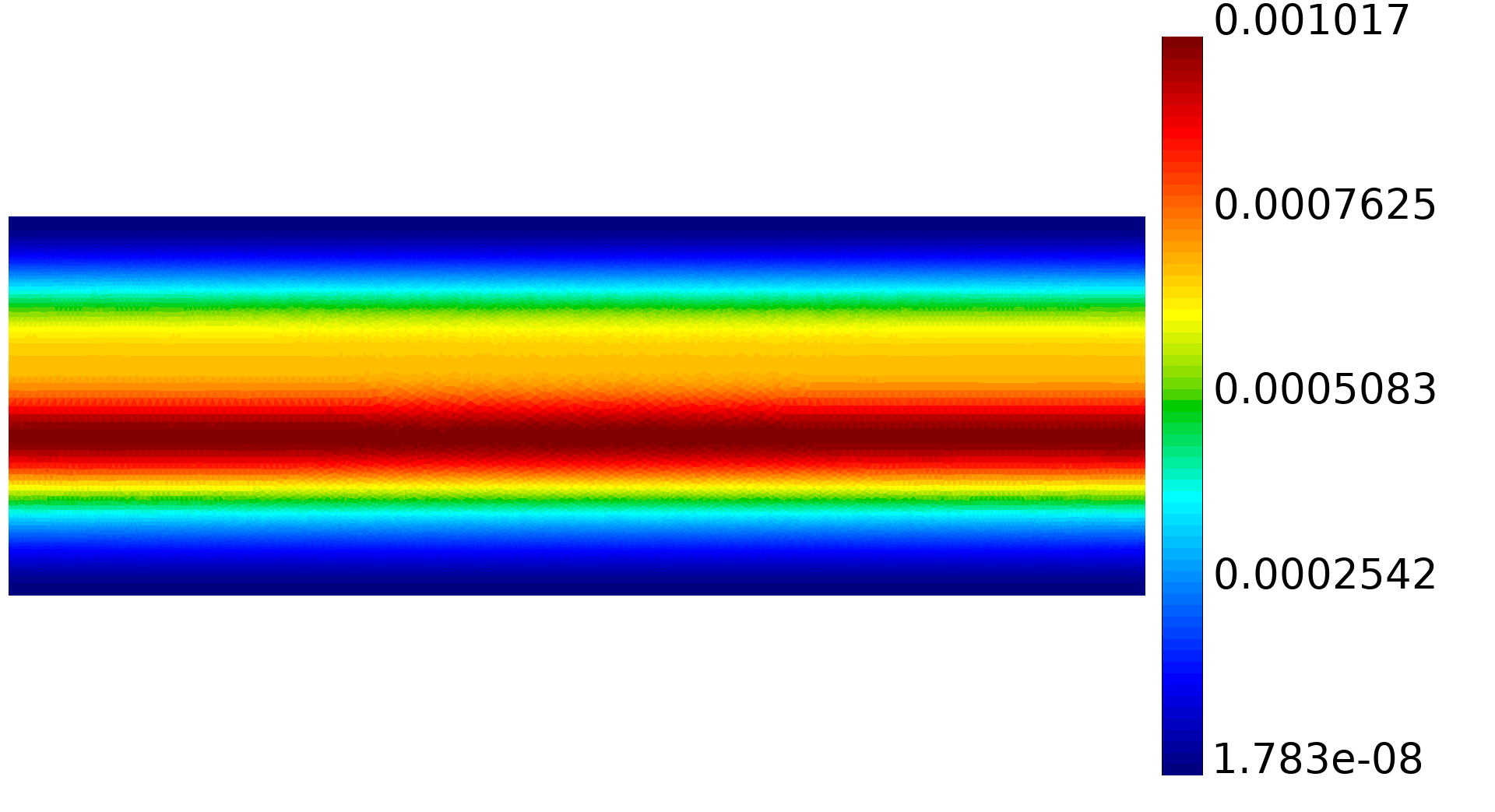}
  \end{subfigure}
  \begin{subfigure}[ht]{.45\textwidth}
    \centering\includegraphics[width=\textwidth]{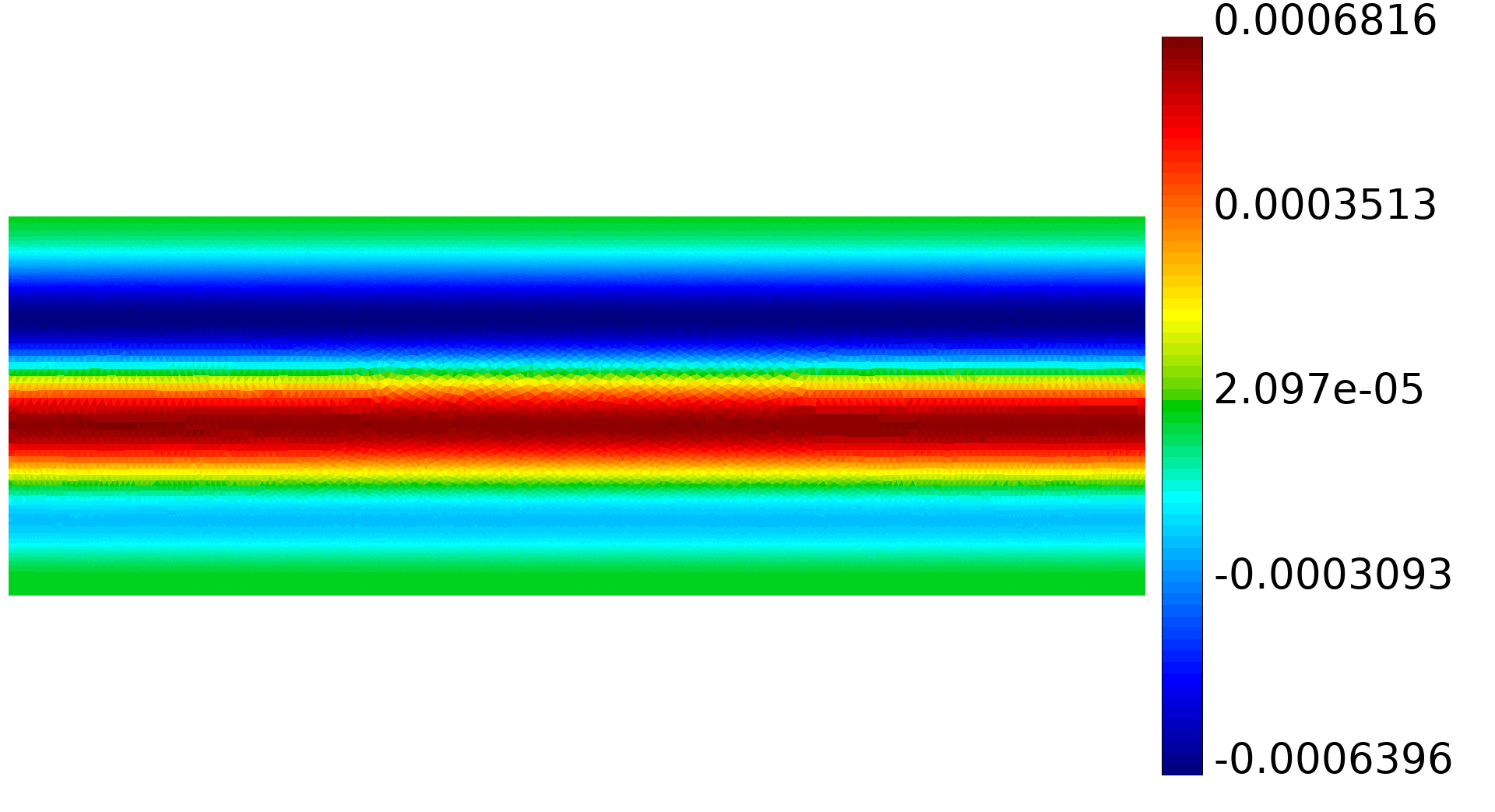}
  \end{subfigure}%
  \begin{subfigure}[ht]{.45\textwidth}
    \centering\includegraphics[width=\textwidth]{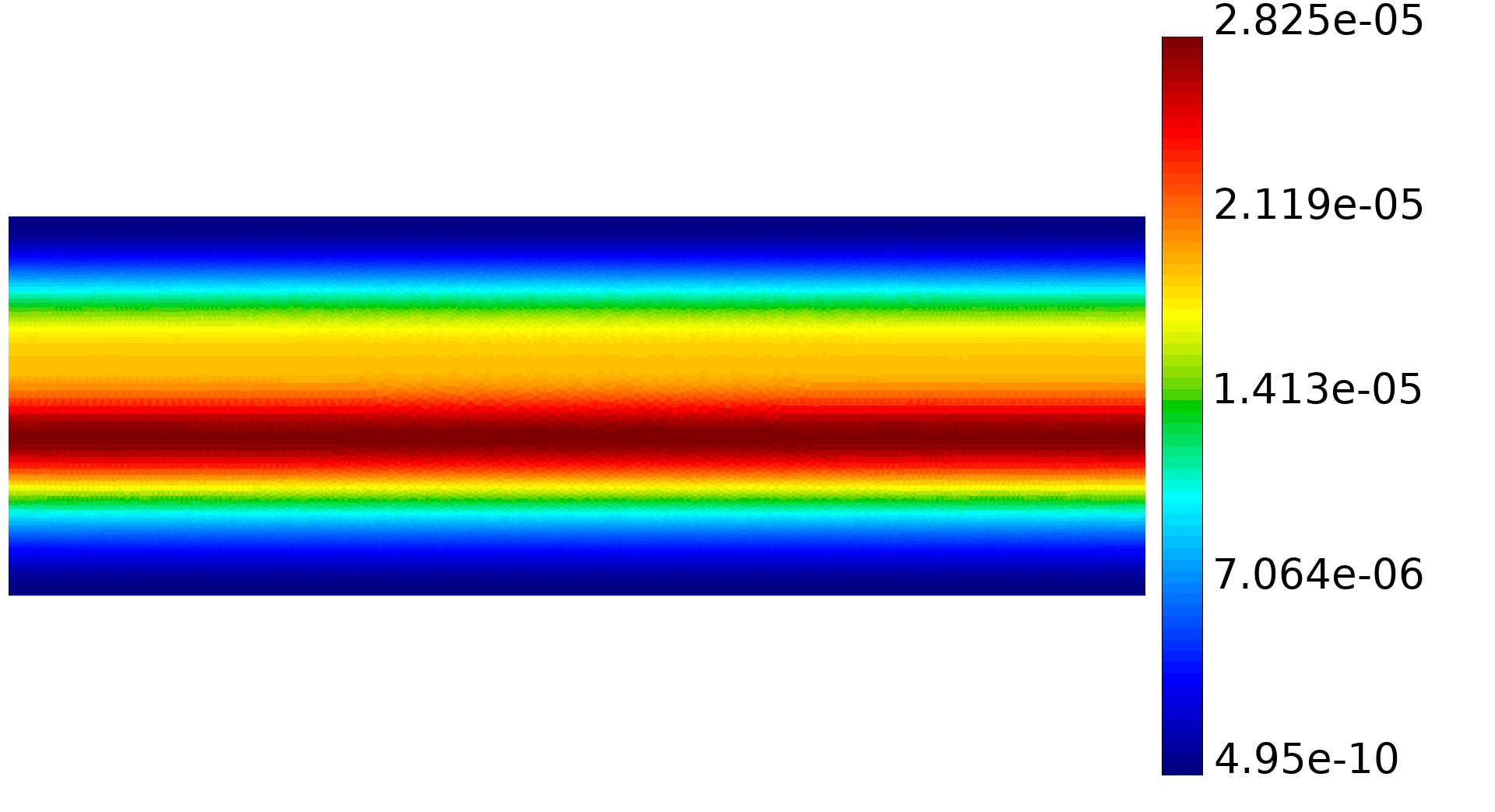}
  \end{subfigure}
  \caption[Channel flow, initial conditions]{\small {Channel flow, 
  as described in {\S}\ref{ss:uqIncompNSTest3};
  (left column) mean and (right column) variance of the initial conditions,
  (top row) horizontal and (bottom row) vertical components of velocity,
  with $120$ samples on the unstructured mesh with resolution level $\ellx = 1$.
  }}
 \label{fig:uqIncompNSTest3Init}
\end{figure}

\begin{figure}[!htb]
\centering
  \begin{subfigure}[ht]{.5\textwidth}
    \centering\includegraphics[width=\textwidth]{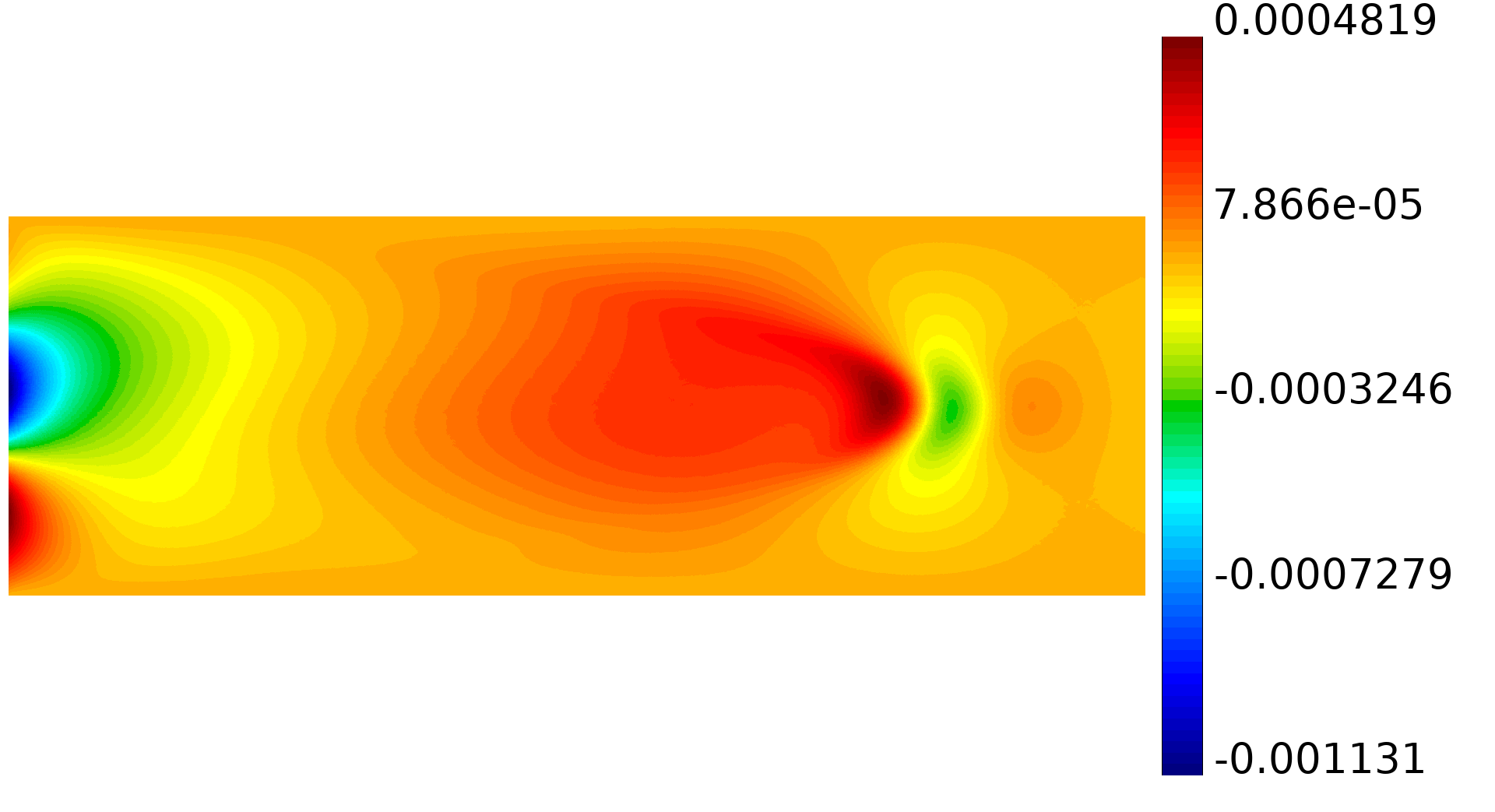}
  \end{subfigure}%
  \begin{subfigure}[ht]{.5\textwidth}
    \centering\includegraphics[width=\textwidth]{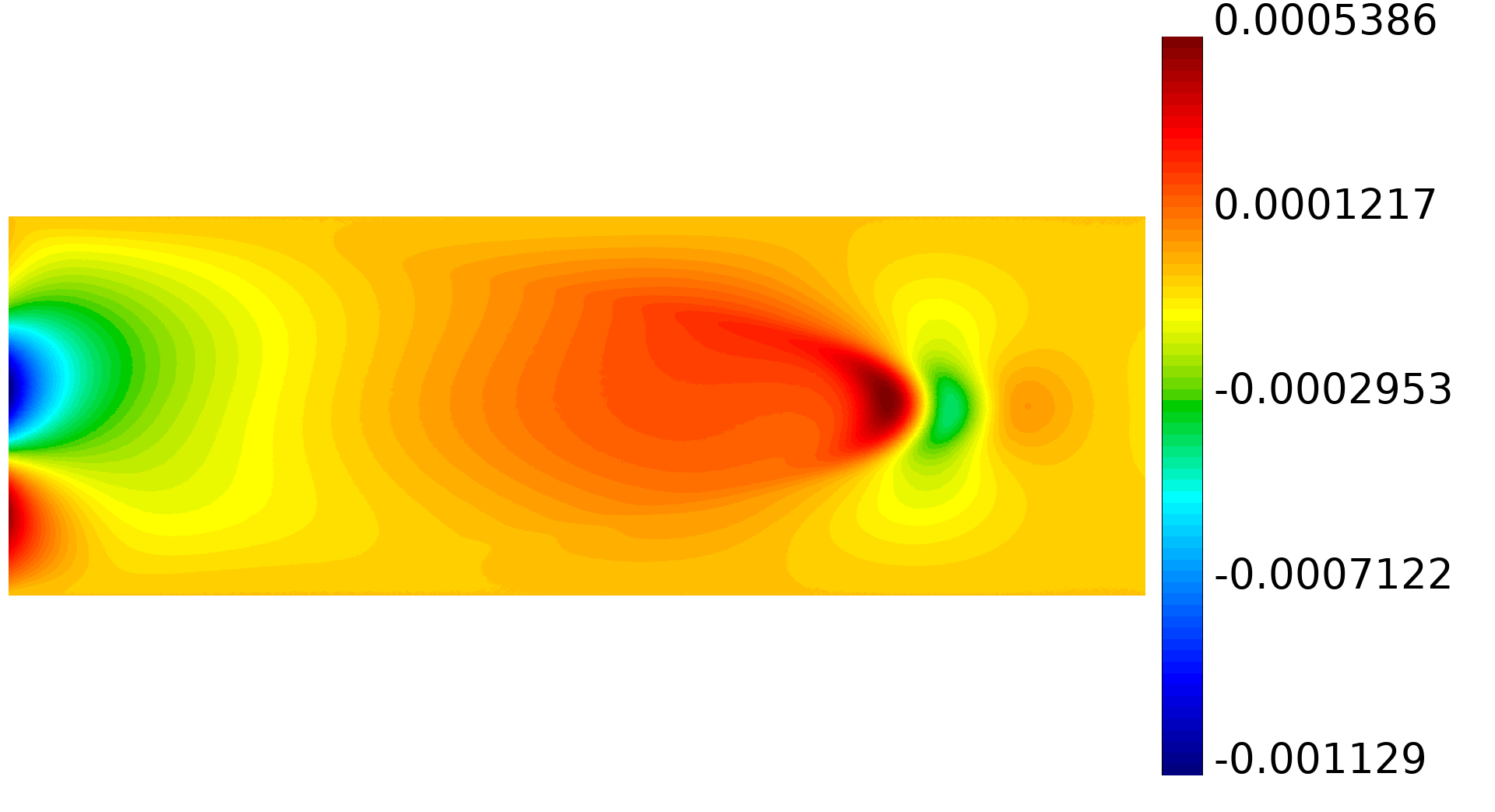}
  \end{subfigure}
  \begin{subfigure}[ht]{.5\textwidth}
    \centering\includegraphics[width=\textwidth]{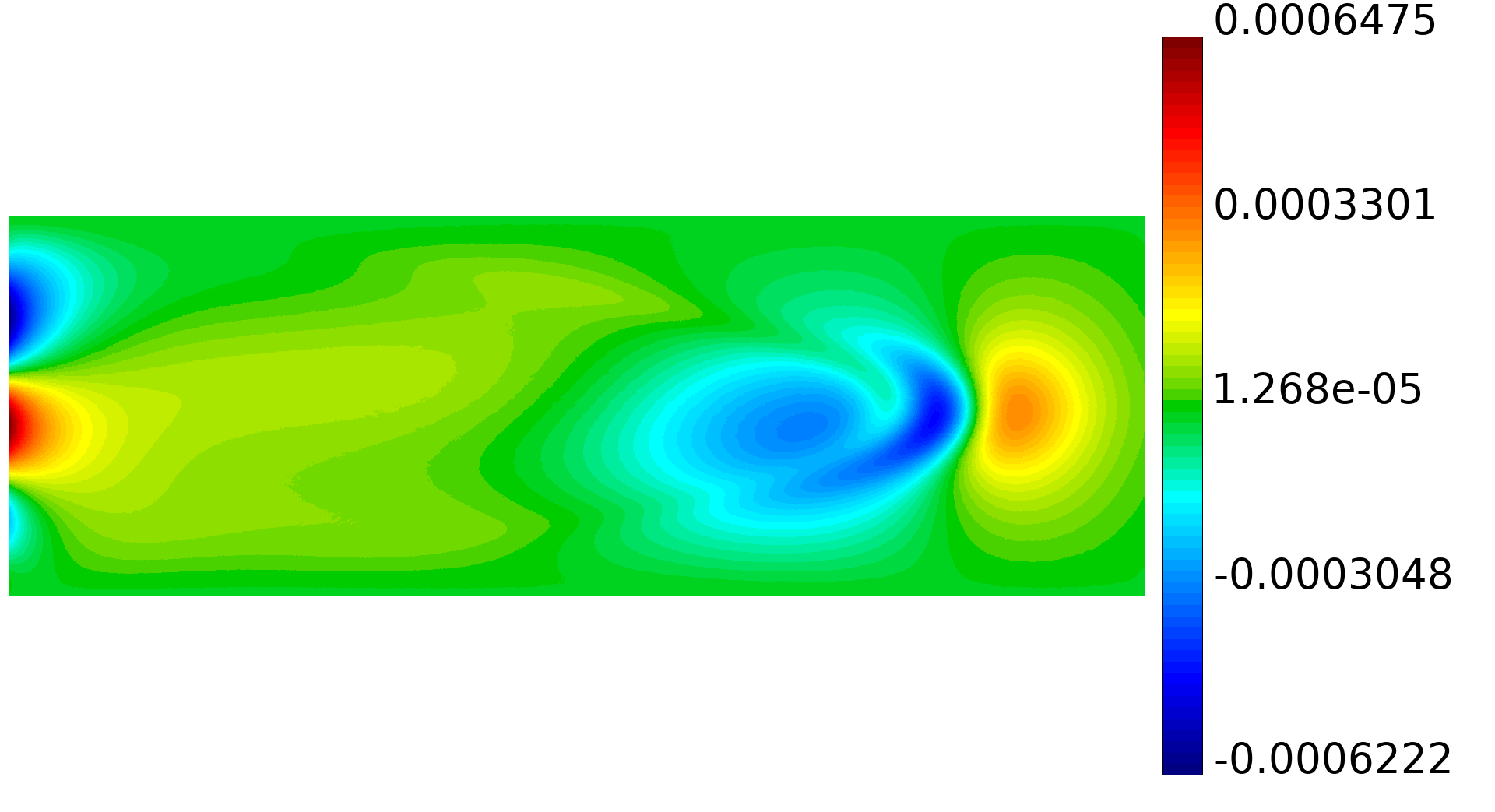}
  \end{subfigure}%
  \begin{subfigure}[ht]{.5\textwidth}
    \centering\includegraphics[width=\textwidth]{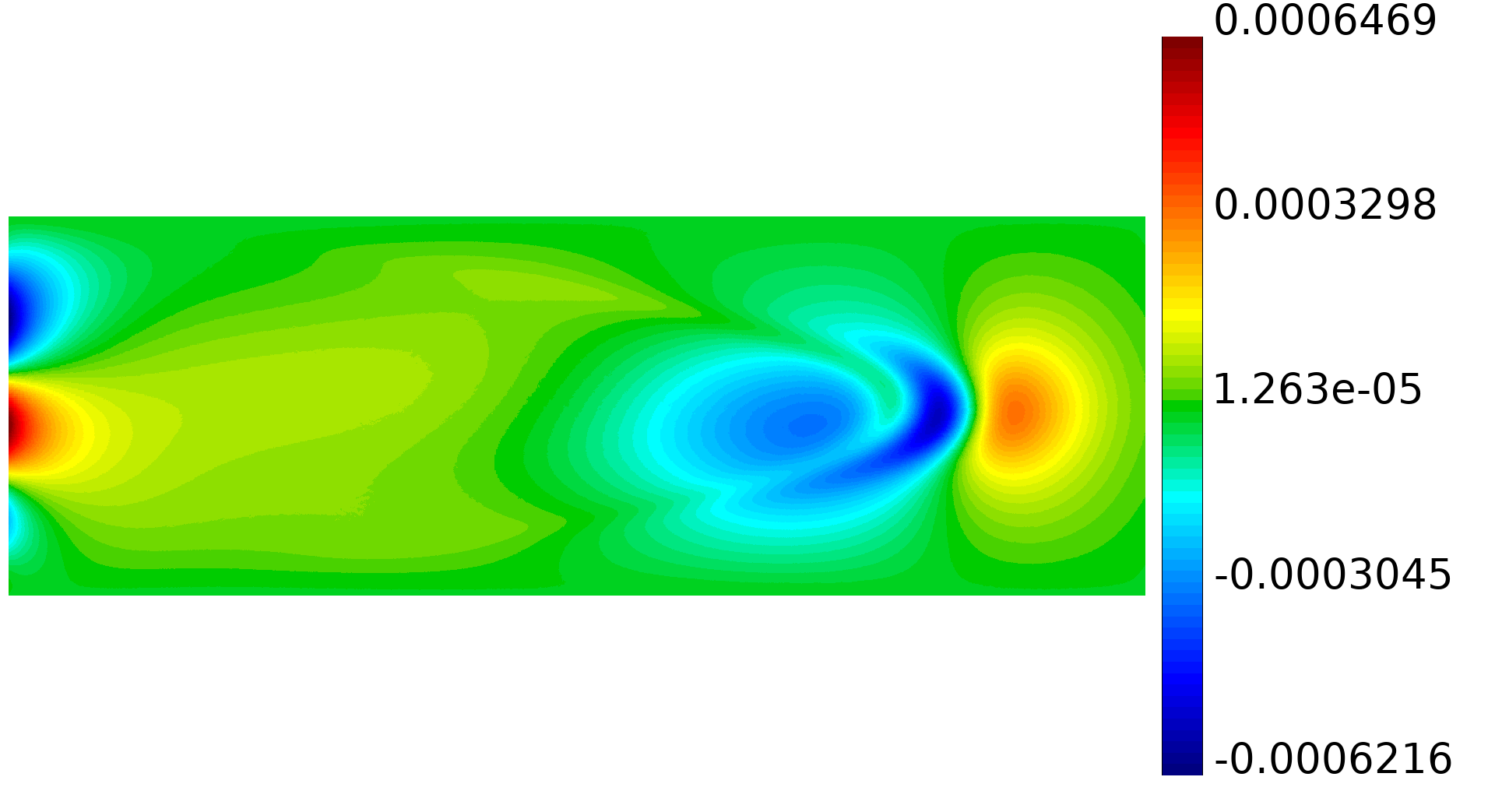}
  \end{subfigure}
  \begin{subfigure}[ht]{.5\textwidth}
    \centering\includegraphics[width=\textwidth]{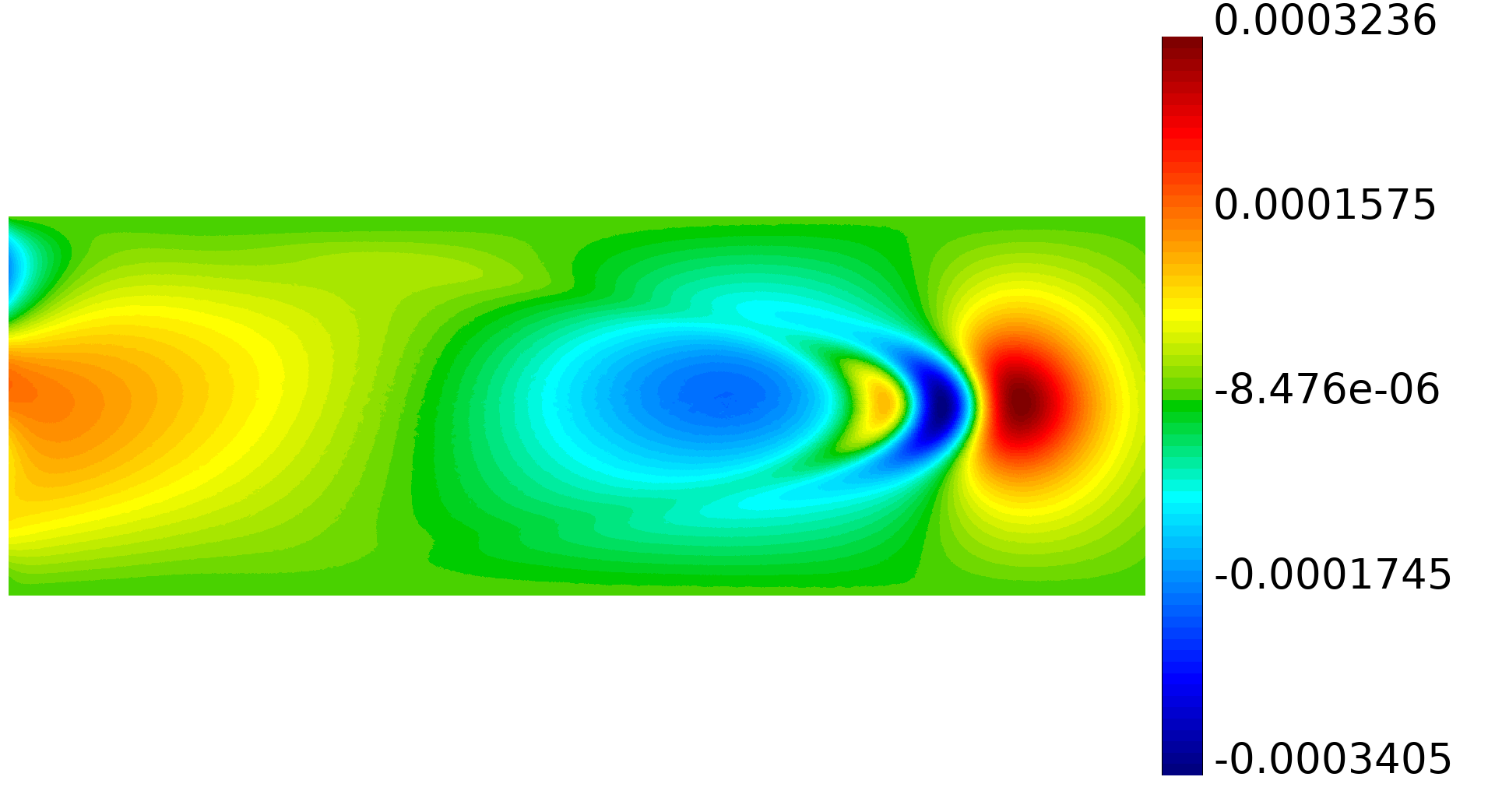}
  \end{subfigure}%
  \begin{subfigure}[ht]{.5\textwidth}
    \centering\includegraphics[width=\textwidth]{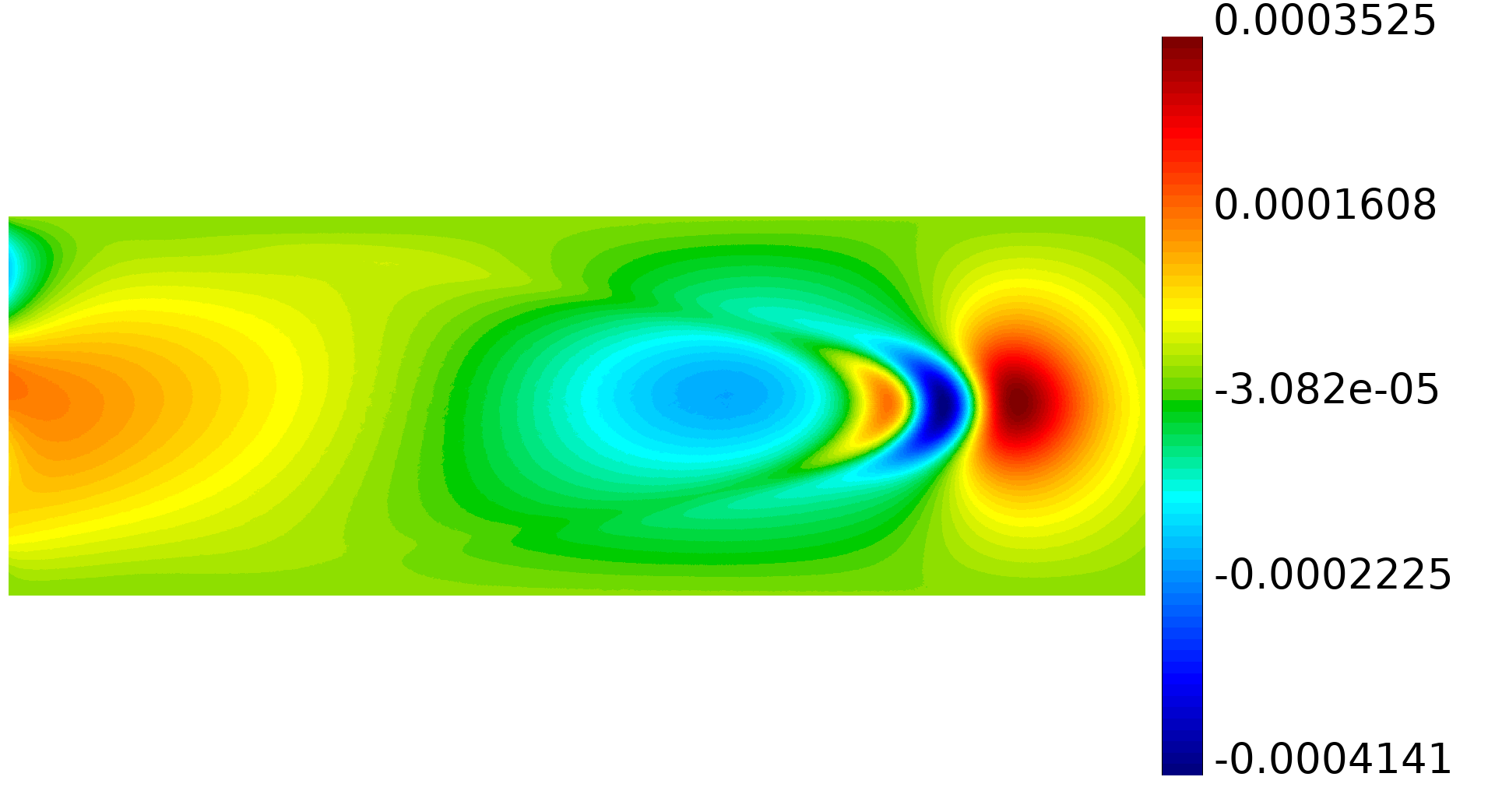}
  \end{subfigure}
  \begin{subfigure}[ht]{.5\textwidth}
    \centering\includegraphics[width=\textwidth]{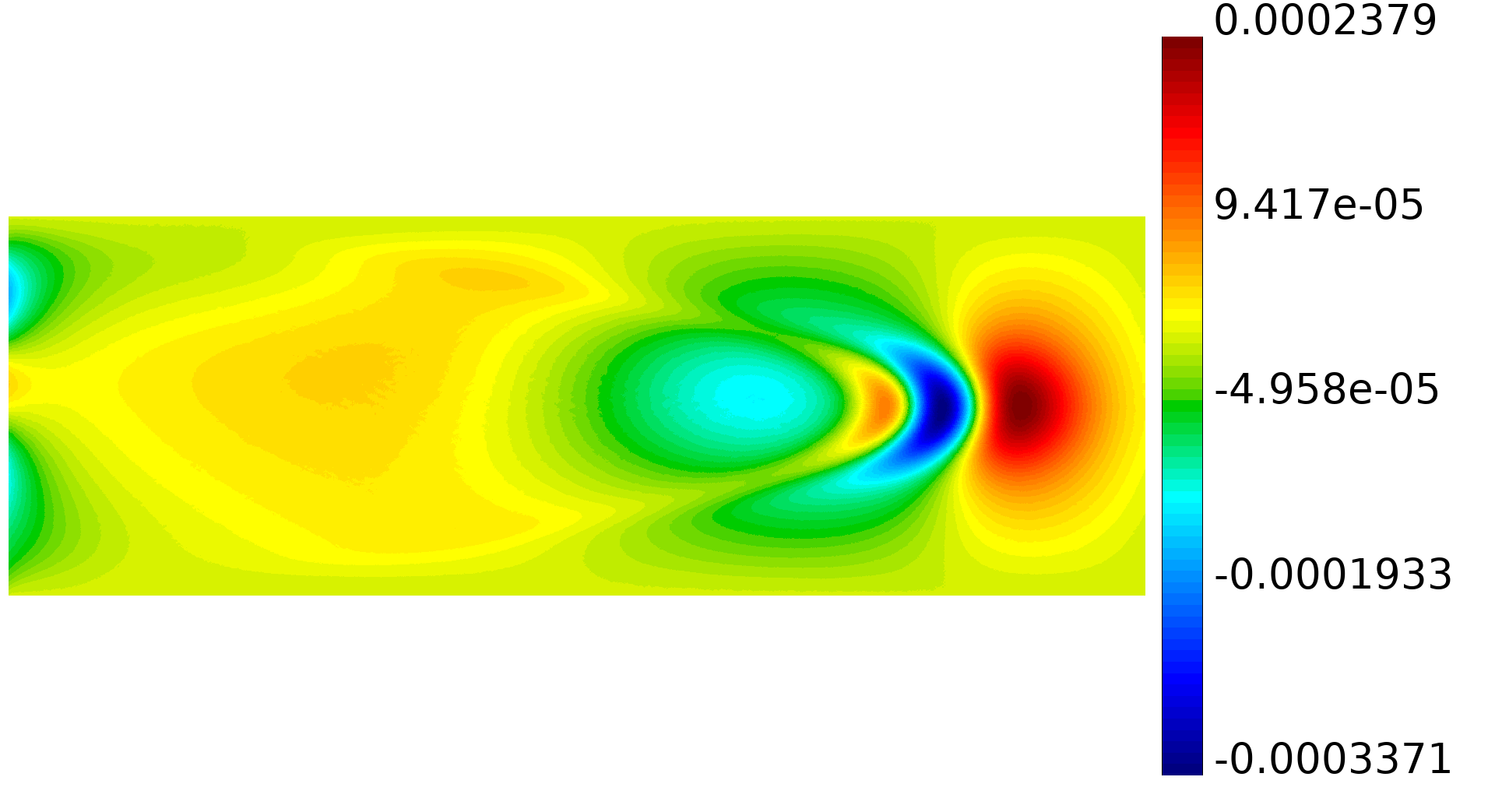}
  \end{subfigure}%
  \begin{subfigure}[ht]{.5\textwidth}
    \centering\includegraphics[width=\textwidth]{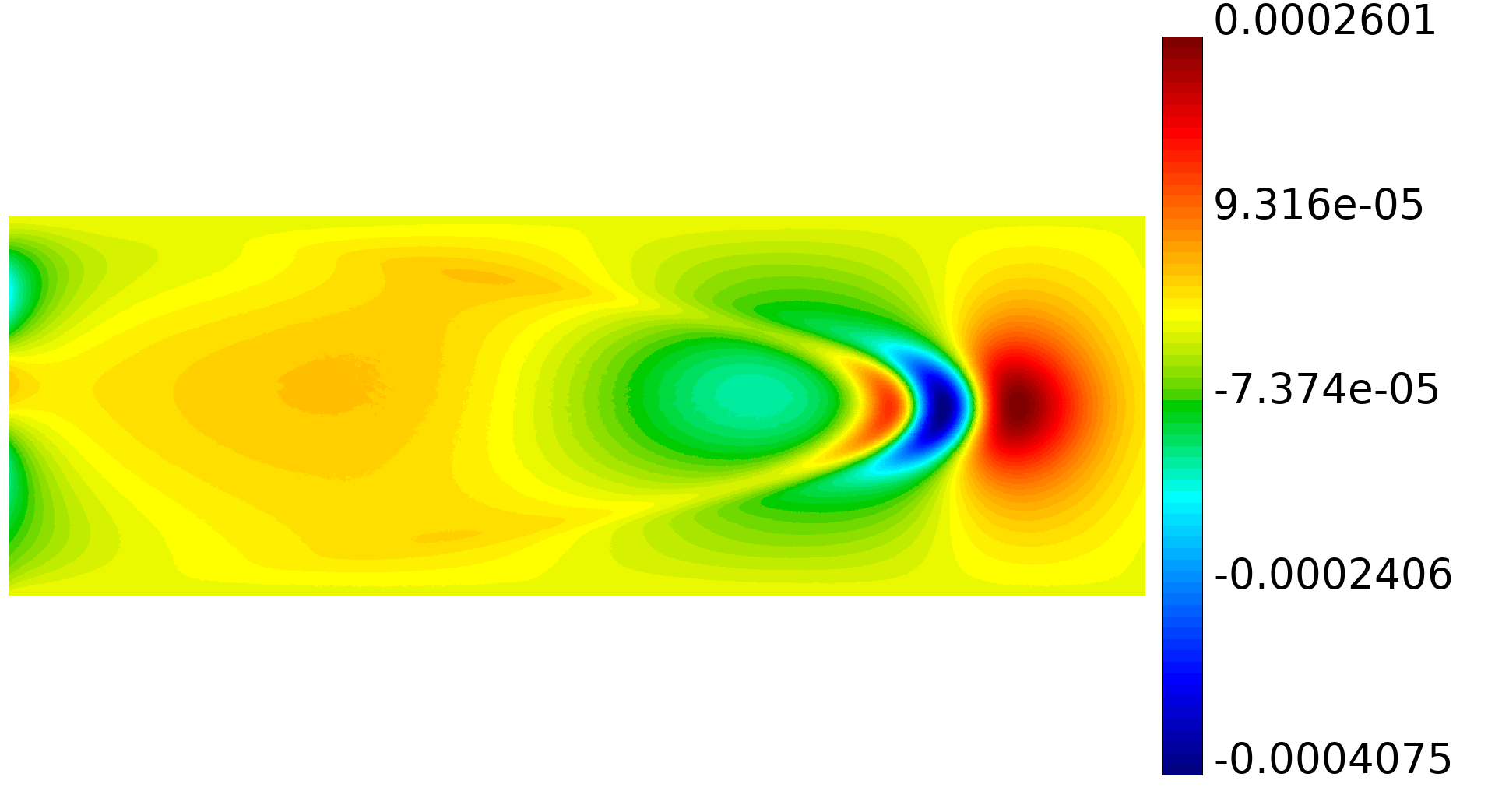}
  \end{subfigure}
  \caption[Channel flow mean of vertical velocity]{\small {Channel flow,
  as described in {\S}\ref{ss:uqIncompNSTest3};
  mean of vertical velocity at $T=0.8$, for (left column) $Re=1600$ and (right column) $Re=3200$, with mesh level $\ellx = 3$ and, 
  from top to bottom, $60, 120, 240$ and $480$ samples.
  }}
 \label{fig:uqIncompNSTest3VisVyMean}
\end{figure}

\begin{figure}[!htb]
\centering
  \begin{subfigure}[ht]{.5\textwidth}
    \centering\includegraphics[width=\textwidth]{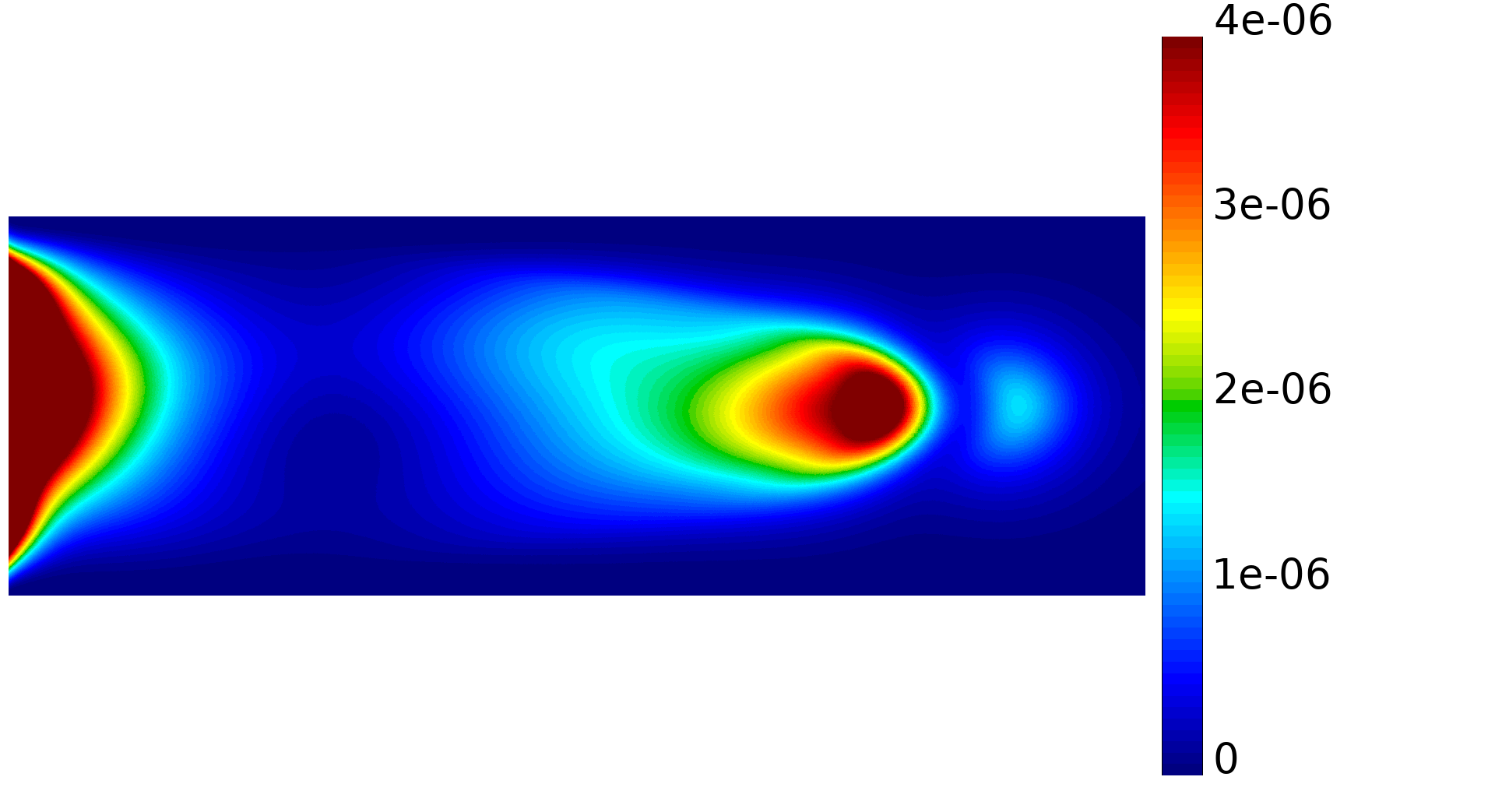}
  \end{subfigure}%
  \begin{subfigure}[ht]{.5\textwidth}
    \centering\includegraphics[width=\textwidth]{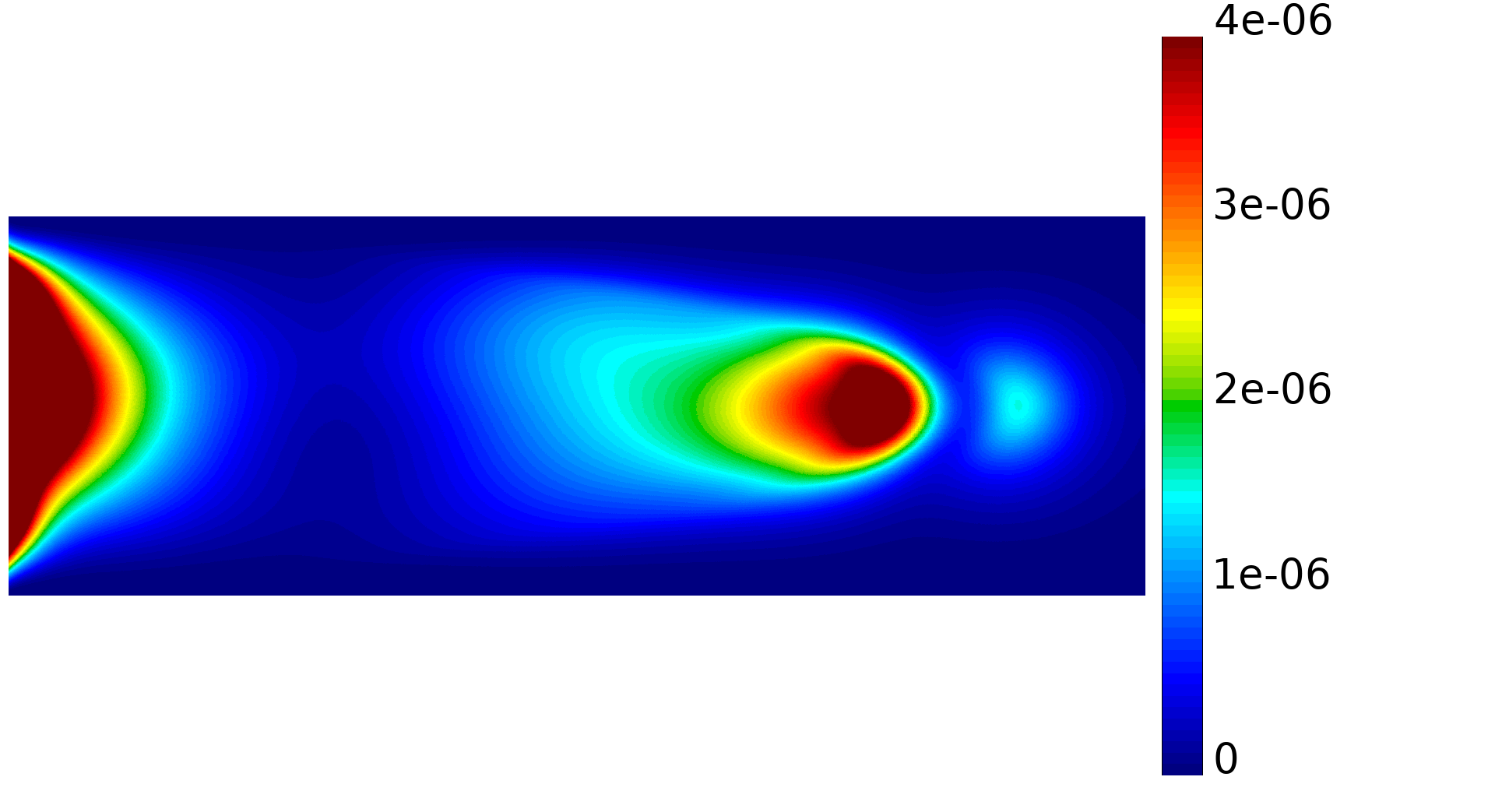}
  \end{subfigure}
  \begin{subfigure}[ht]{.5\textwidth}
    \centering\includegraphics[width=\textwidth]{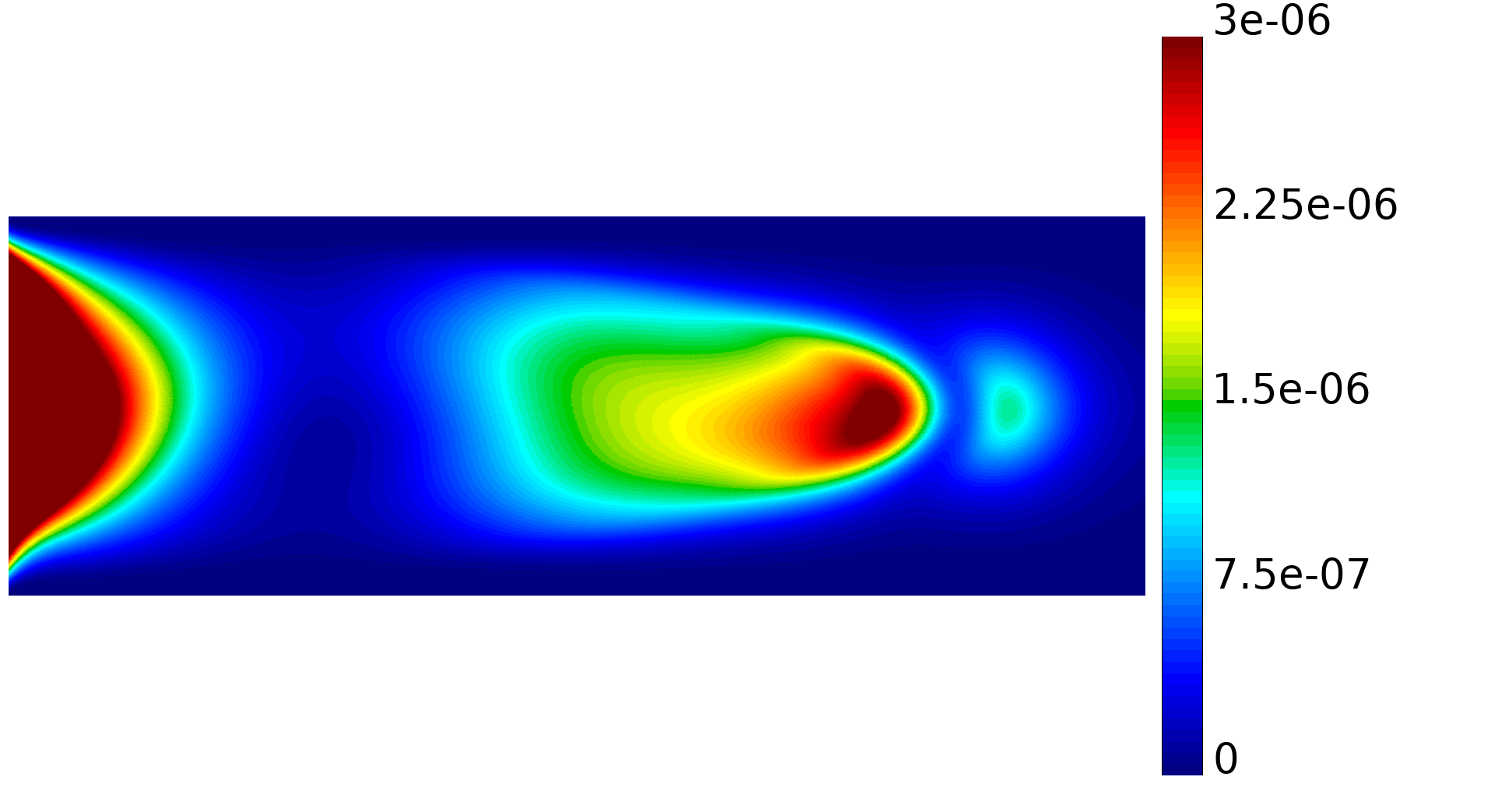}
  \end{subfigure}%
  \begin{subfigure}[ht]{.5\textwidth}
    \centering\includegraphics[width=\textwidth]{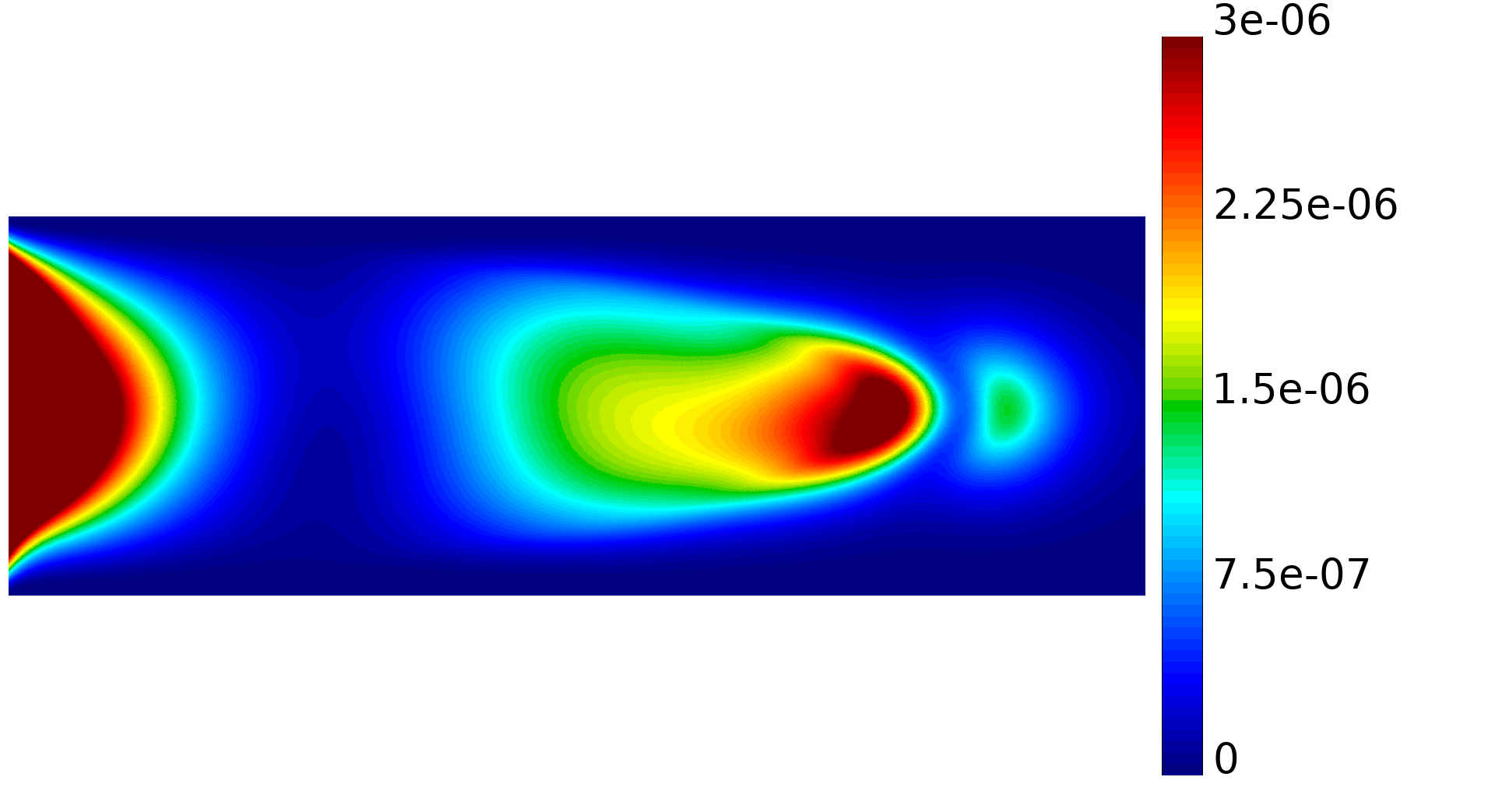}
  \end{subfigure}
  \begin{subfigure}[ht]{.5\textwidth}
    \centering\includegraphics[width=\textwidth]{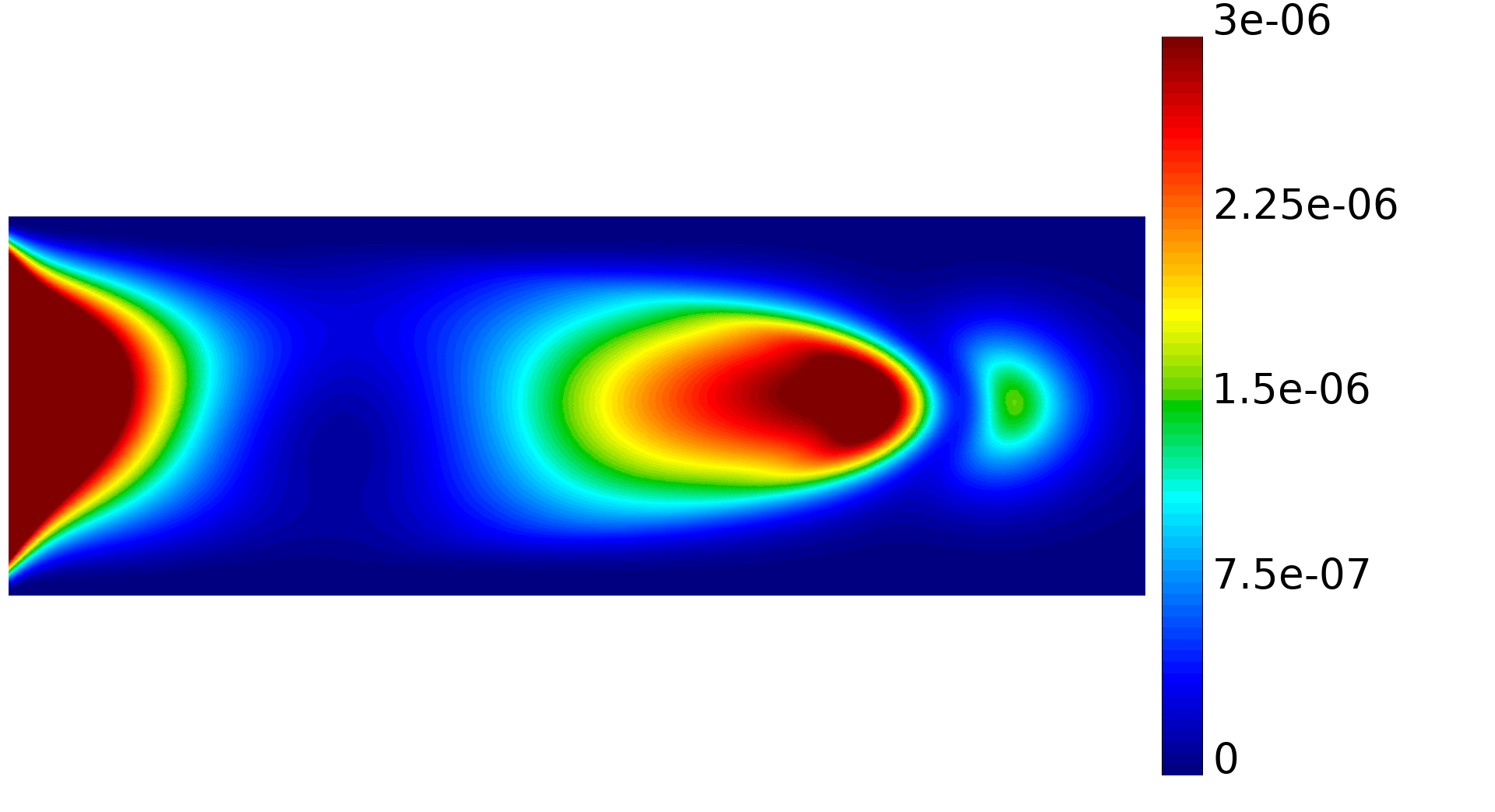}
  \end{subfigure}%
  \begin{subfigure}[ht]{.5\textwidth}
    \centering\includegraphics[width=\textwidth]{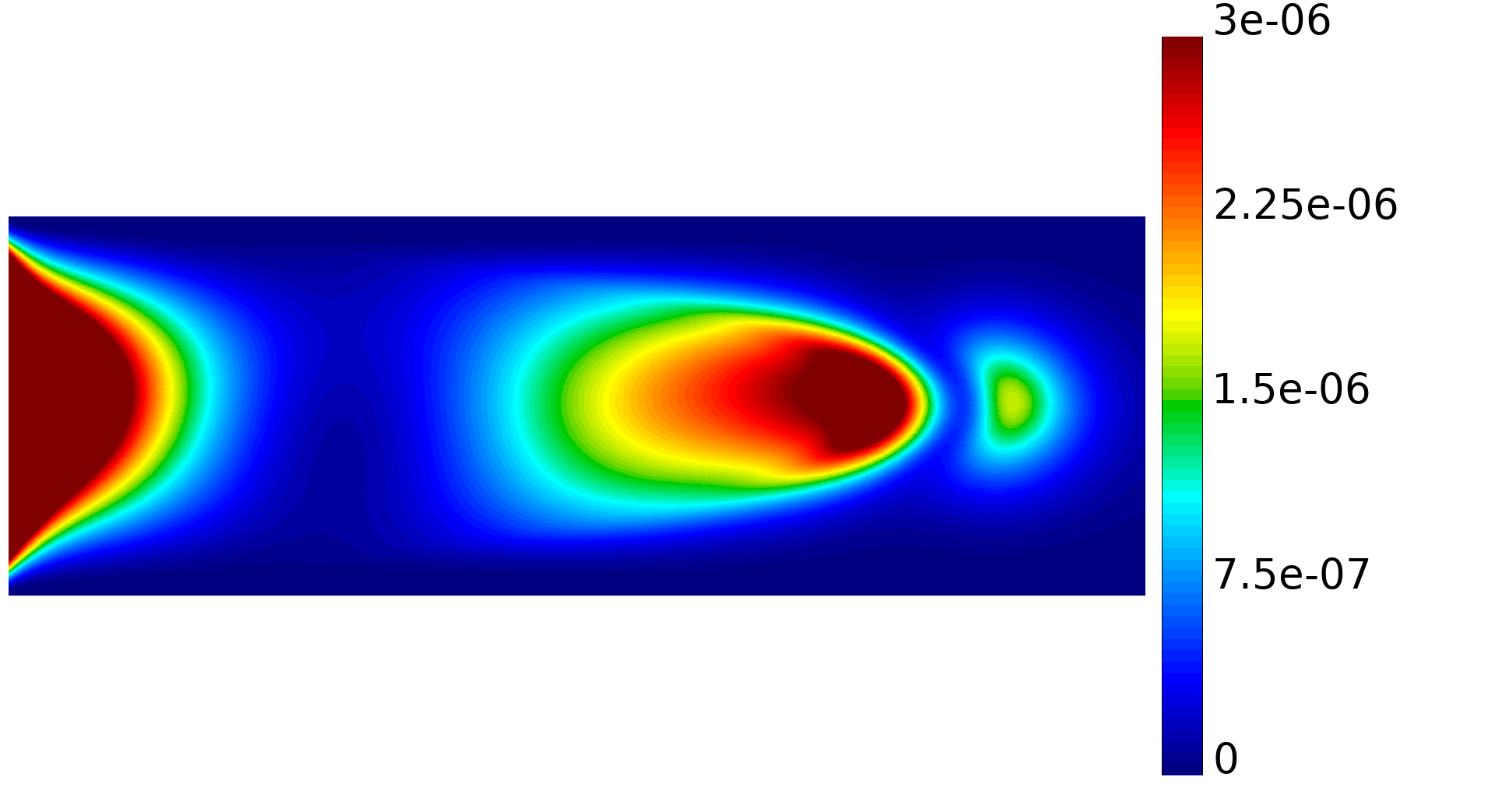}
  \end{subfigure}
  \begin{subfigure}[ht]{.5\textwidth}
    \centering\includegraphics[width=\textwidth]{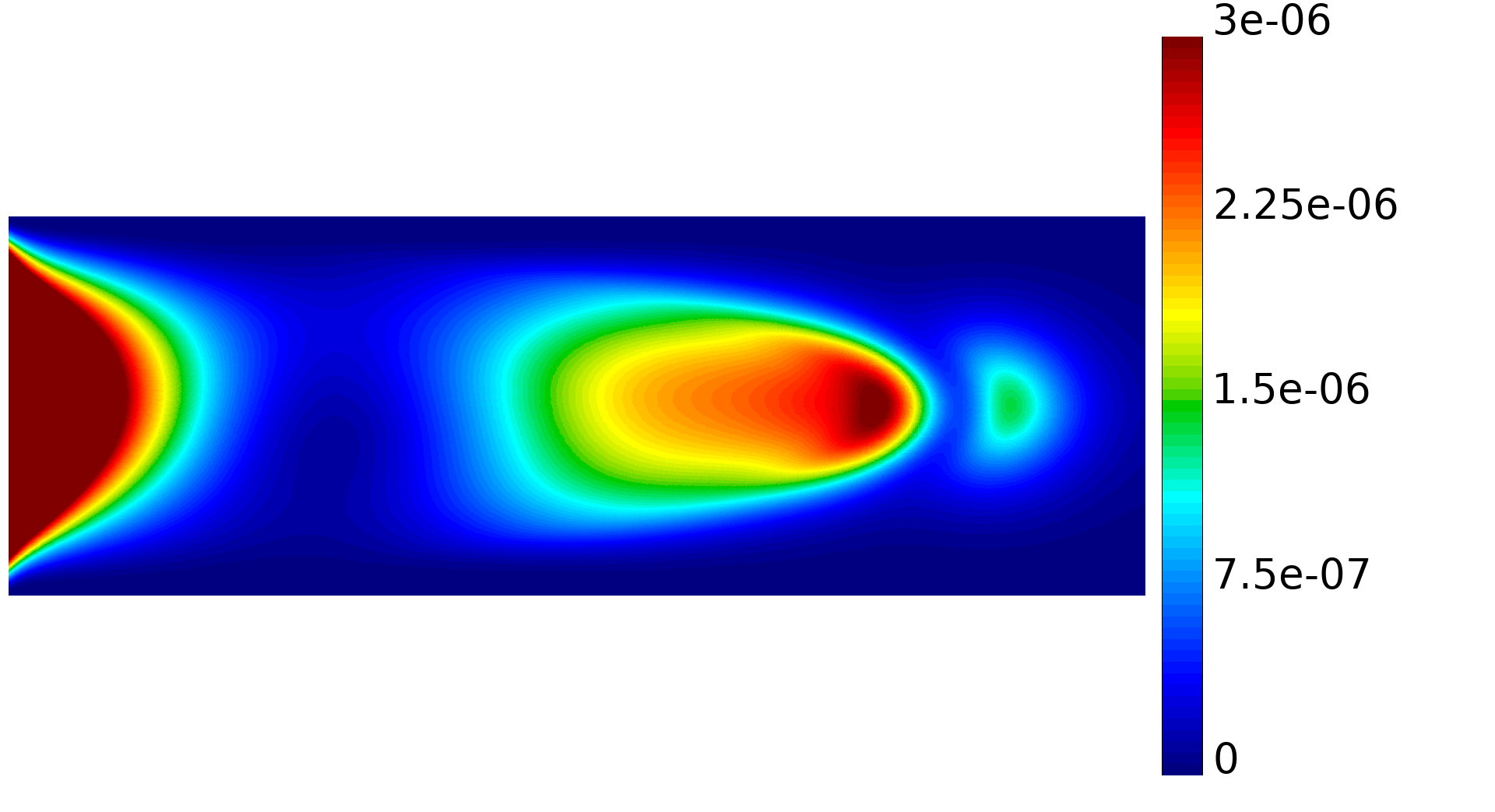}
  \end{subfigure}%
  \begin{subfigure}[ht]{.5\textwidth}
    \centering\includegraphics[width=\textwidth]{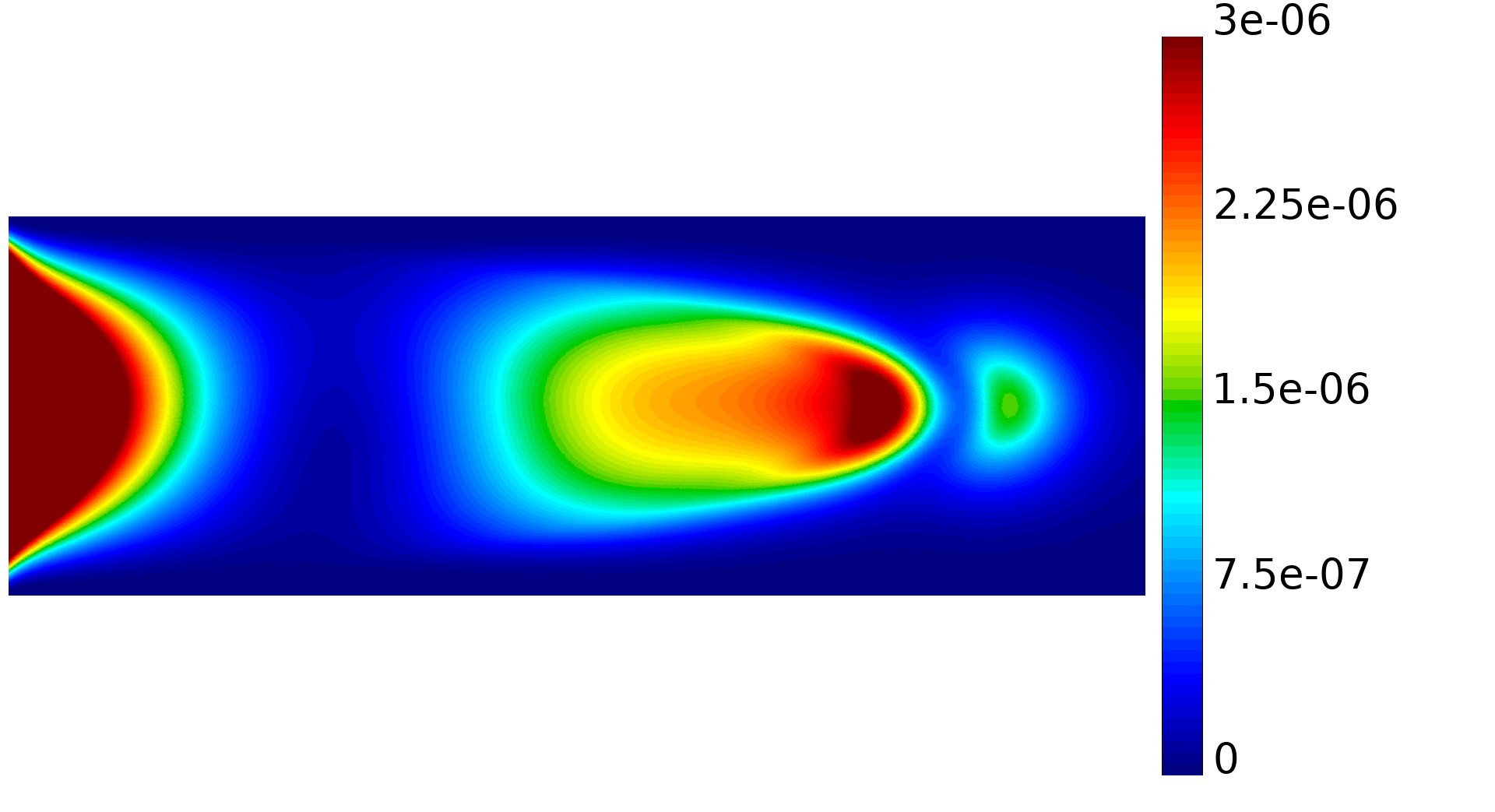}
  \end{subfigure}
  \caption[Channel flow re-scaled variance of vertical velocity]{\small {Channel flow,
  as described in {\S}\ref{ss:uqIncompNSTest3};
  re-scaled variance of vertical velocity at $T=0.8$, for (left column) $Re=1600$ and (right column) $Re=3200$, with mesh level $\ellx = 3$ and, 
  from top to bottom, $60, 120, 240$ and $480$ samples.
  }}
 \label{fig:uqIncompNSTest3VisVyVarRescaled}
\end{figure}

\begin{figure}[!htb]
\centering
  \begin{subfigure}[t]{.5\textwidth}
    \centering\includegraphics[width = \textwidth]{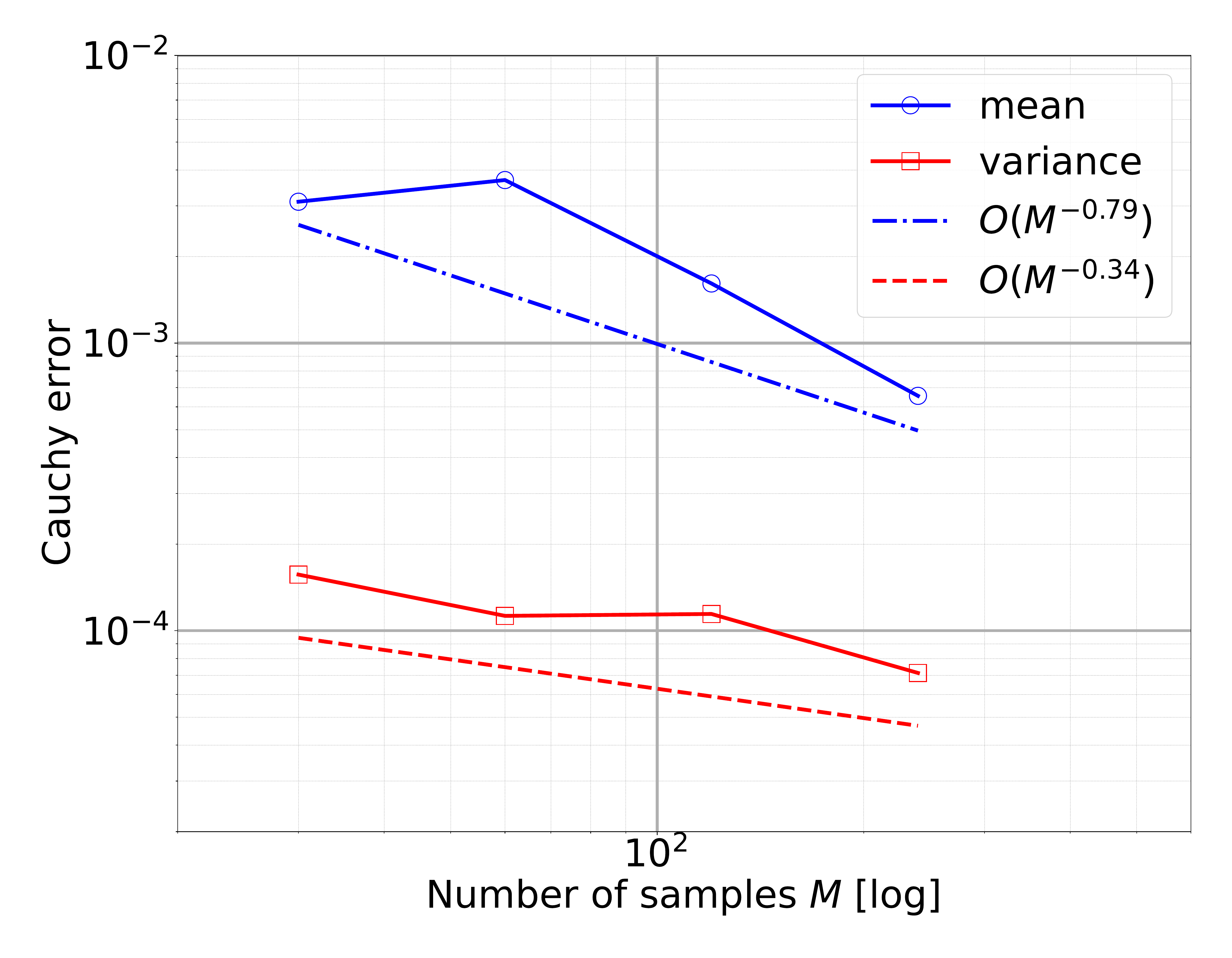}
  \end{subfigure}
 \caption[Channel flow, $Re=1600$, convergence of mean and variance]{\small {Channel flow with $Re=1600$, 
  as described in {\S}\ref{ss:uqIncompNSTest3};
  Cauchy error of mean and variance measured in $\ltwo{\Dx}$ norm
  between samples $30, 60, 120, 240$ and $480$,
  for mesh level $\ellx = 3$, at time $T=0.8$.
 }}
 \label{fig:uqIncompNSTest3Re1600Convg}
\end{figure}

\begin{figure}[!htb]
\centering
  \begin{subfigure}[t]{.5\textwidth}
    \centering\includegraphics[width = \textwidth]{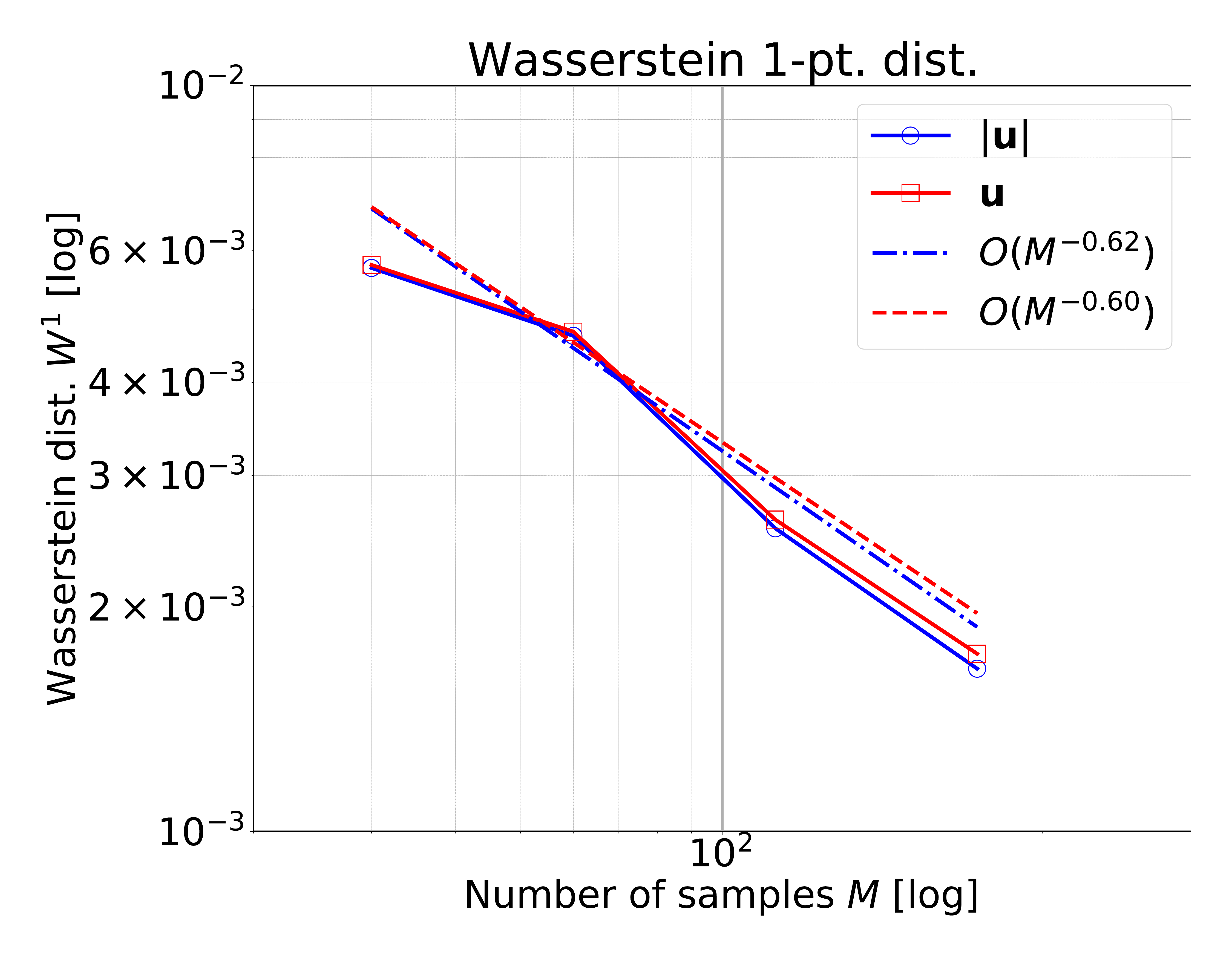}
  \end{subfigure}%
  \begin{subfigure}[t]{.5\textwidth}
    \centering\includegraphics[width = \textwidth]{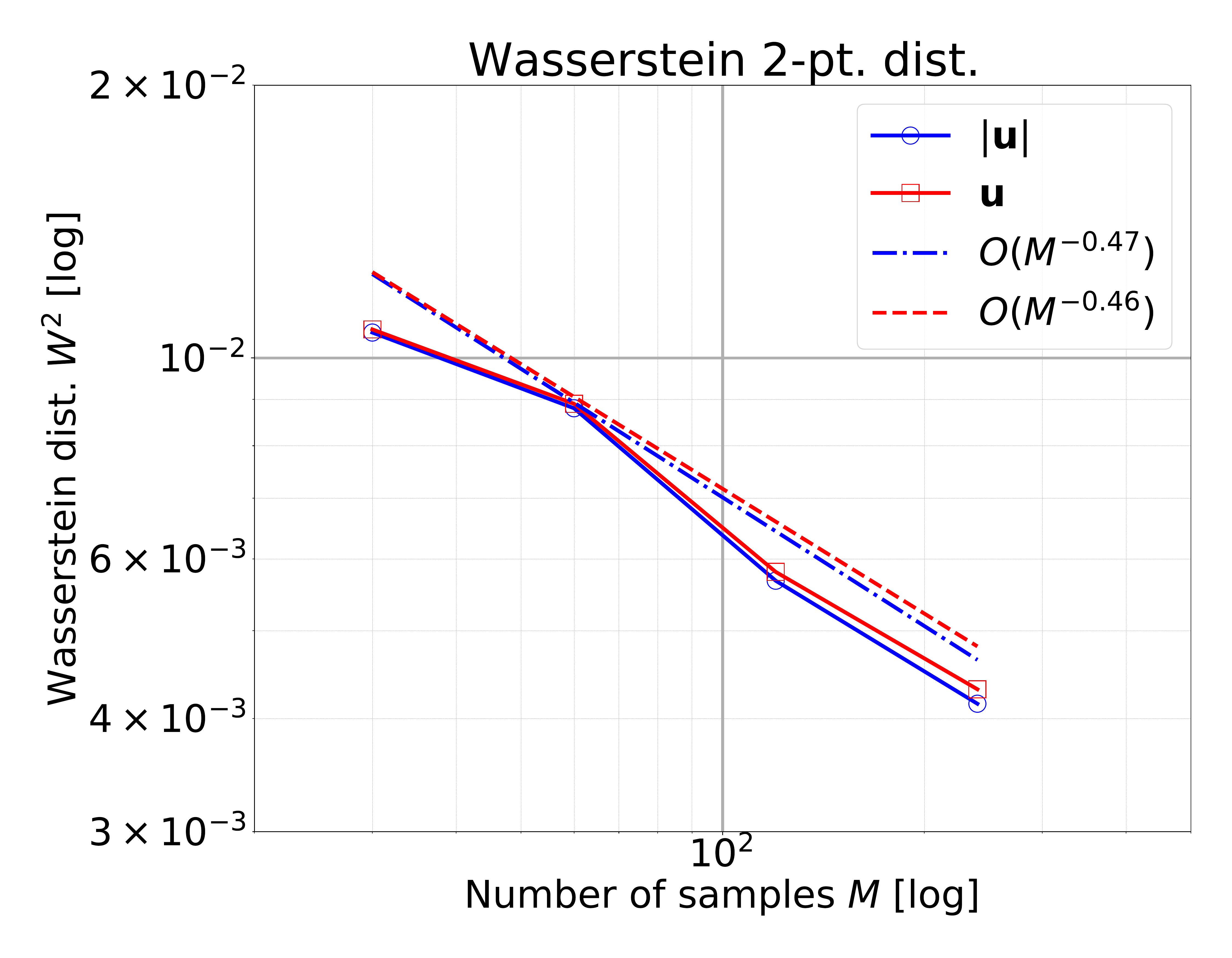}
  \end{subfigure}
 \caption[Channel flow, $Re=1600$, convergence of Wasserstein distances]{\small {Channel flow with $Re=1600$, 
  as described in {\S}\ref{ss:uqIncompNSTest3};
  Wasserstein distances of velocity field $\bu$ and its magnitude $\abs{\bu}$
  measured between samples $30, 60, 120, 240$ and $480$,
  for mesh level $\ellx = 3$, at time $T=0.8$.
  }}
 \label{fig:uqIncompNSTest3Re1600Wd}
\end{figure}

\begin{figure}[!htb]
\centering
  \begin{subfigure}[t]{.5\textwidth}
    \centering\includegraphics[width = \textwidth]{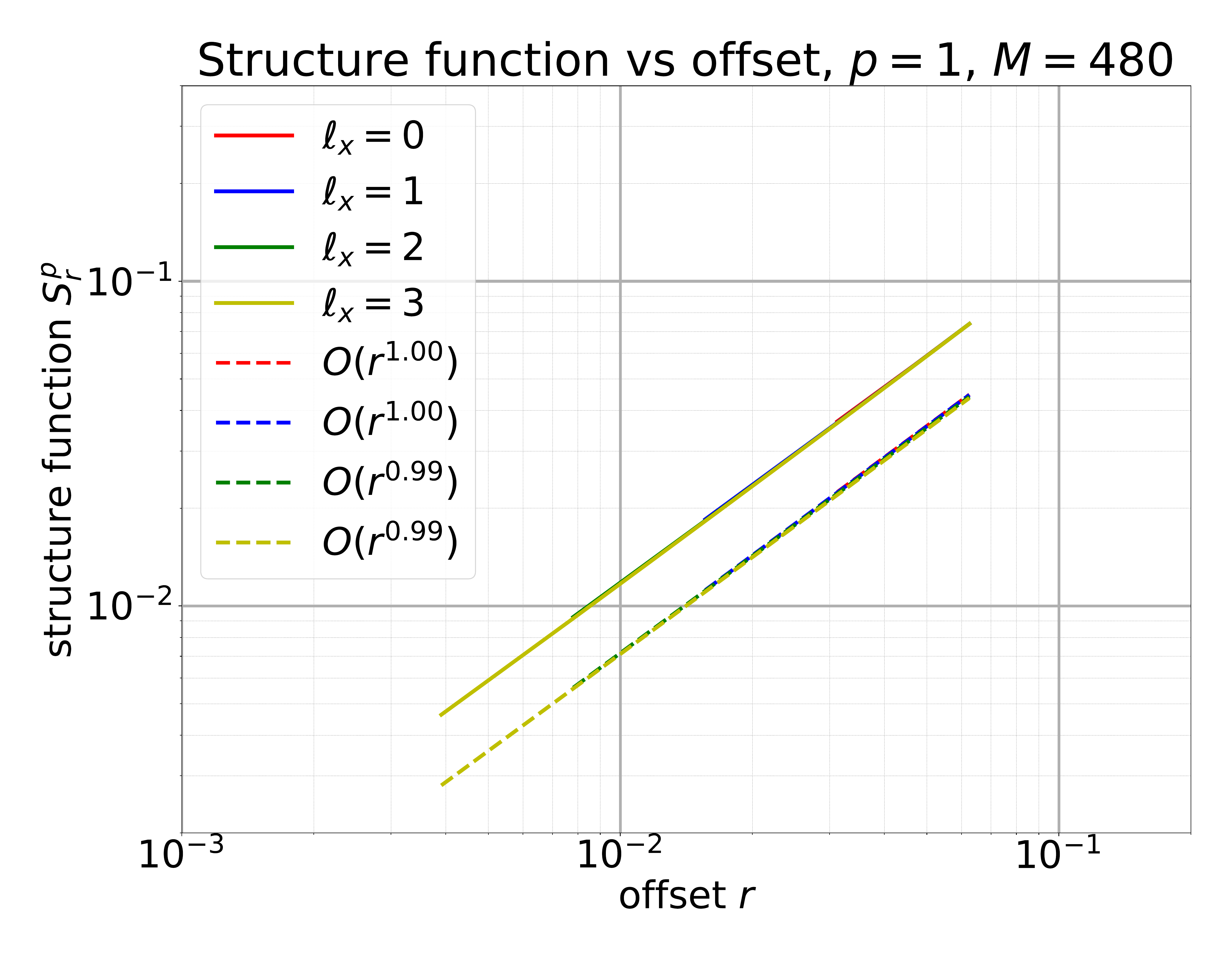}
  \end{subfigure}
  \begin{subfigure}[t]{.5\textwidth}
    \centering\includegraphics[width = \textwidth]{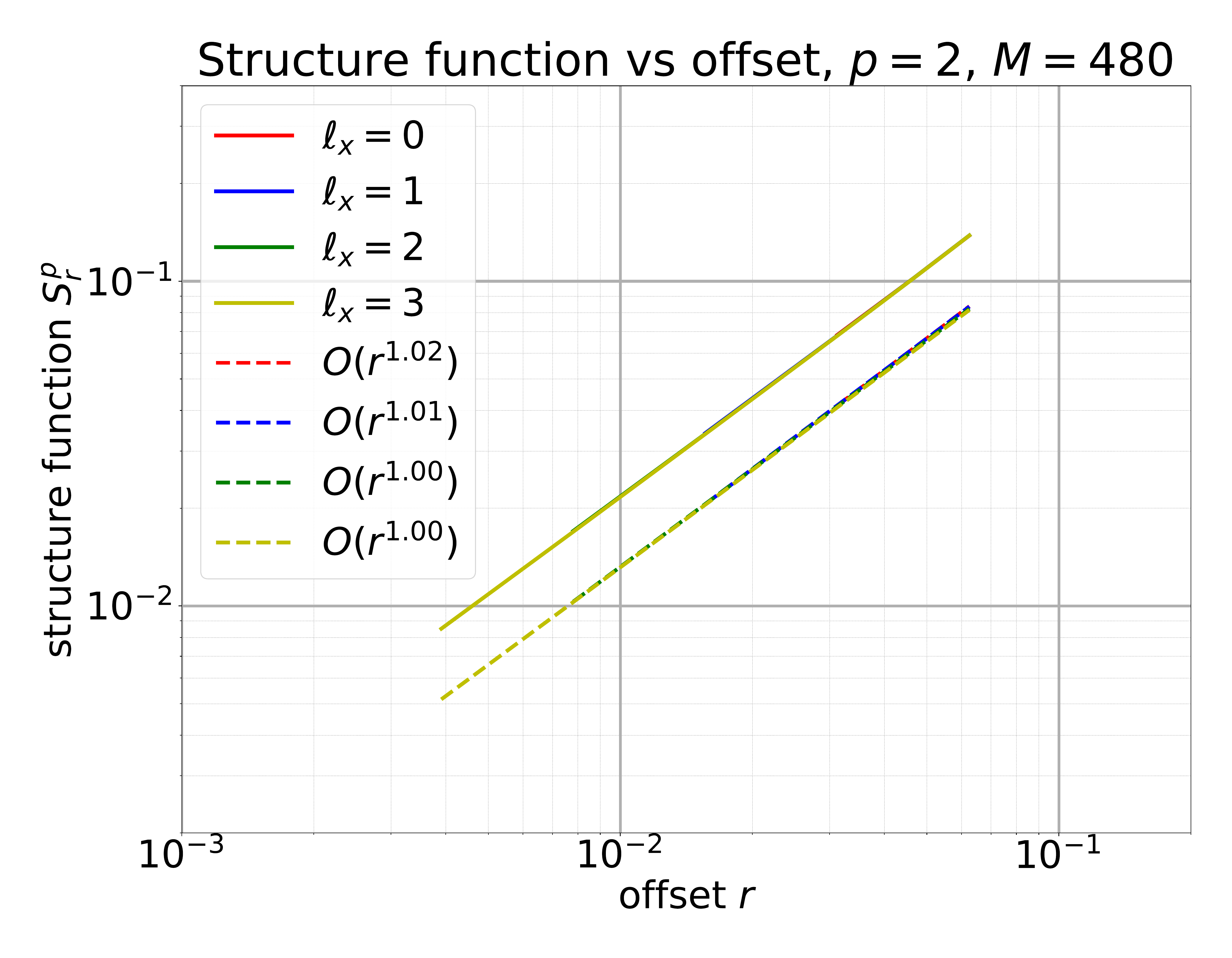}
  \end{subfigure}
  \begin{subfigure}[t]{.5\textwidth}
    \centering\includegraphics[width = \textwidth]{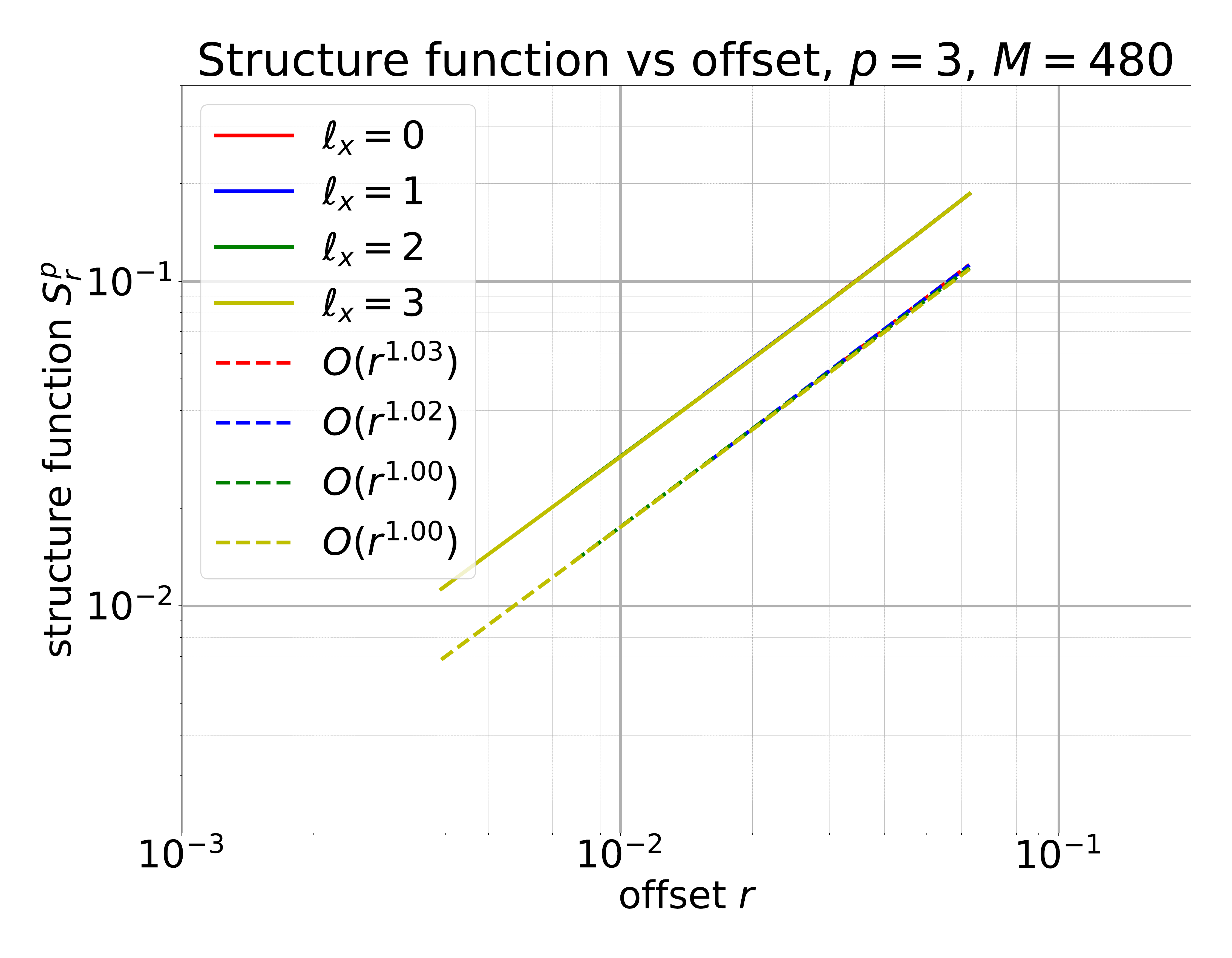}
  \end{subfigure}
 \caption[Channel flow, $Re=1600$, convergence of structure function]{\small {Channel flow with $Re=1600$, 
  as described in {\S}\ref{ss:uqIncompNSTest3};
  Approximation of the structure function of the velocity field with
  $480$ samples, at time $T=0.8$.
 }}
 \label{fig:uqIncompNSTest3Re1600Scube}
\end{figure}

\begin{figure}[!htb]
\centering
  \begin{subfigure}[t]{.5\textwidth}
    \centering\includegraphics[width = \textwidth]{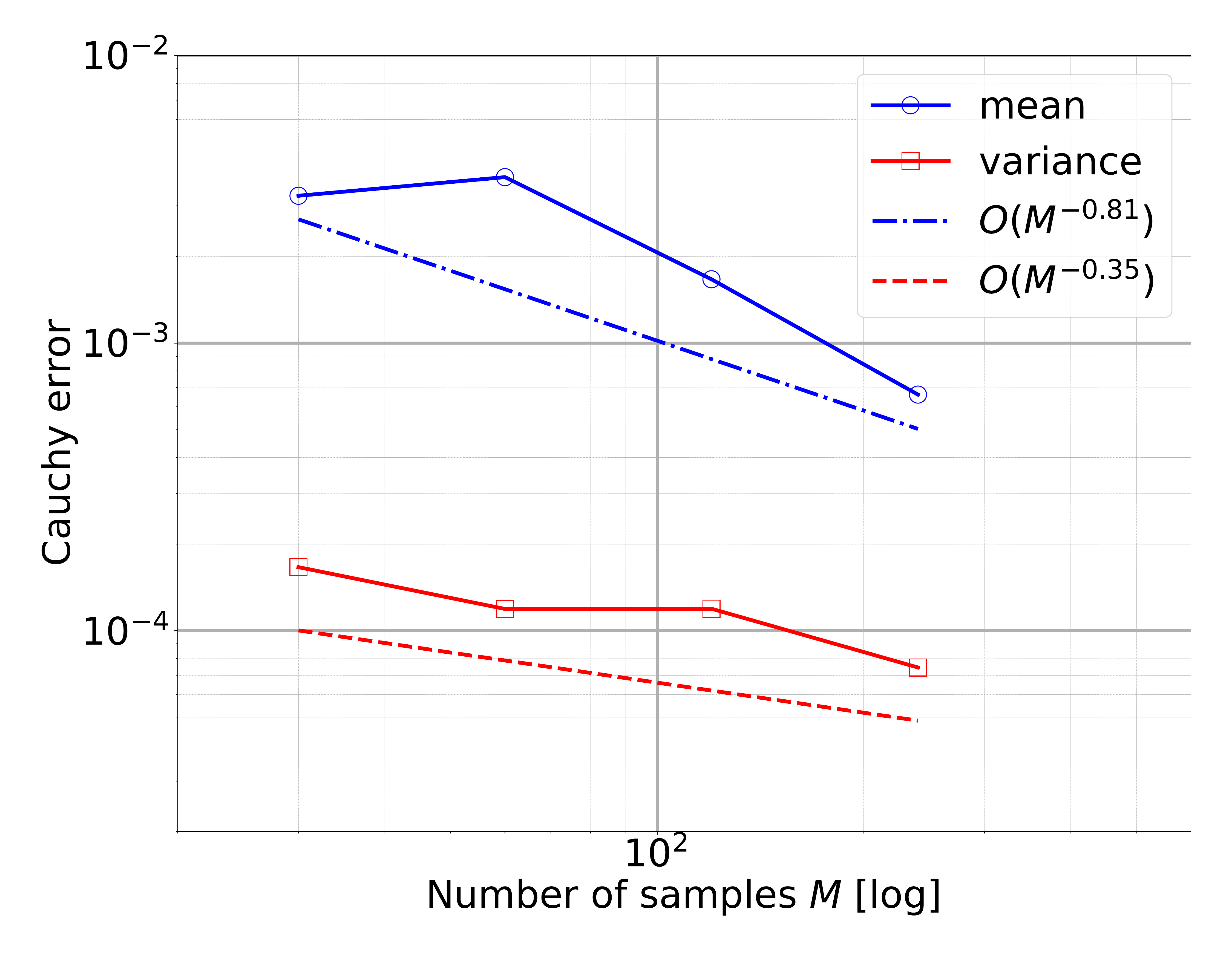}
  \end{subfigure}
 \caption[Channel flow, $Re=3200$, convergence of mean and variance]{\small {Channel flow with $Re=3200$, 
  as described in {\S}\ref{ss:uqIncompNSTest3};
  Cauchy error of mean and variance measured in $\ltwo{\Dx}$ norm
  between samples $30, 60, 120, 240$ and $480$,
  for mesh level $\ellx = 3$, at time $T=0.8$.
 }}
 \label{fig:uqIncompNSTest3Re3200Convg}
\end{figure}

\begin{figure}[!htb]
\centering
  \begin{subfigure}[t]{.5\textwidth}
    \centering\includegraphics[width = \textwidth]{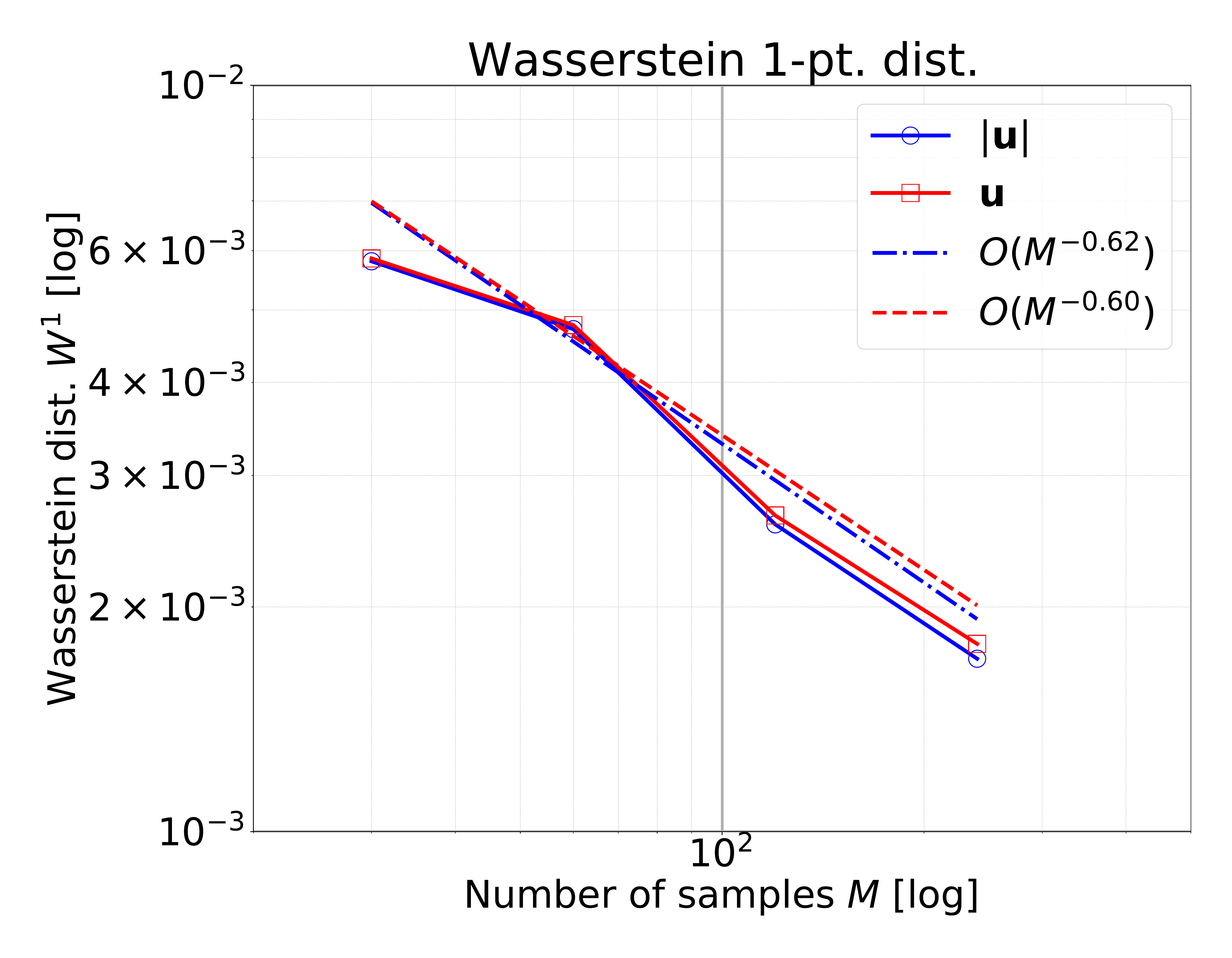}
  \end{subfigure}%
  \begin{subfigure}[t]{.5\textwidth}
    \centering\includegraphics[width = \textwidth]{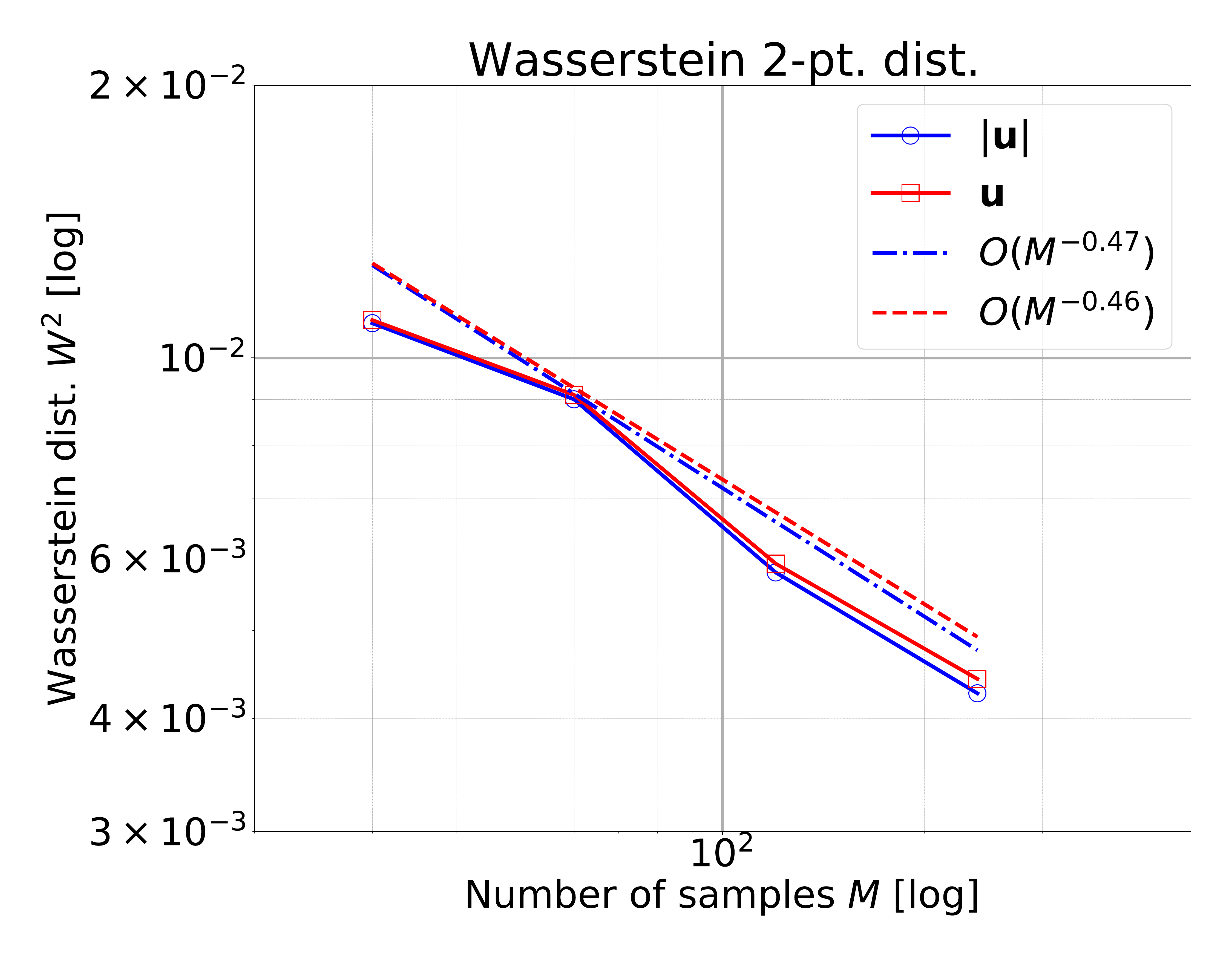}
  \end{subfigure}
 \caption[Channel flow, $Re=3200$, convergence of Wasserstein distances]{\small {Channel flow with $Re=3200$, 
  as described in {\S}\ref{ss:uqIncompNSTest3};
  Wasserstein distances of velocity field $\bu$ and its magnitude $\abs{\bu}$ 
  measured between samples $30, 60, 120, 240$ and $480$,
  for mesh level $\ellx = 3$, at time $T=0.8$.
  }}
 \label{fig:uqIncompNSTest3Re3200Wd}
\end{figure}

\begin{figure}[!htb]
\centering
  \begin{subfigure}[t]{.5\textwidth}
    \centering\includegraphics[width = \textwidth]{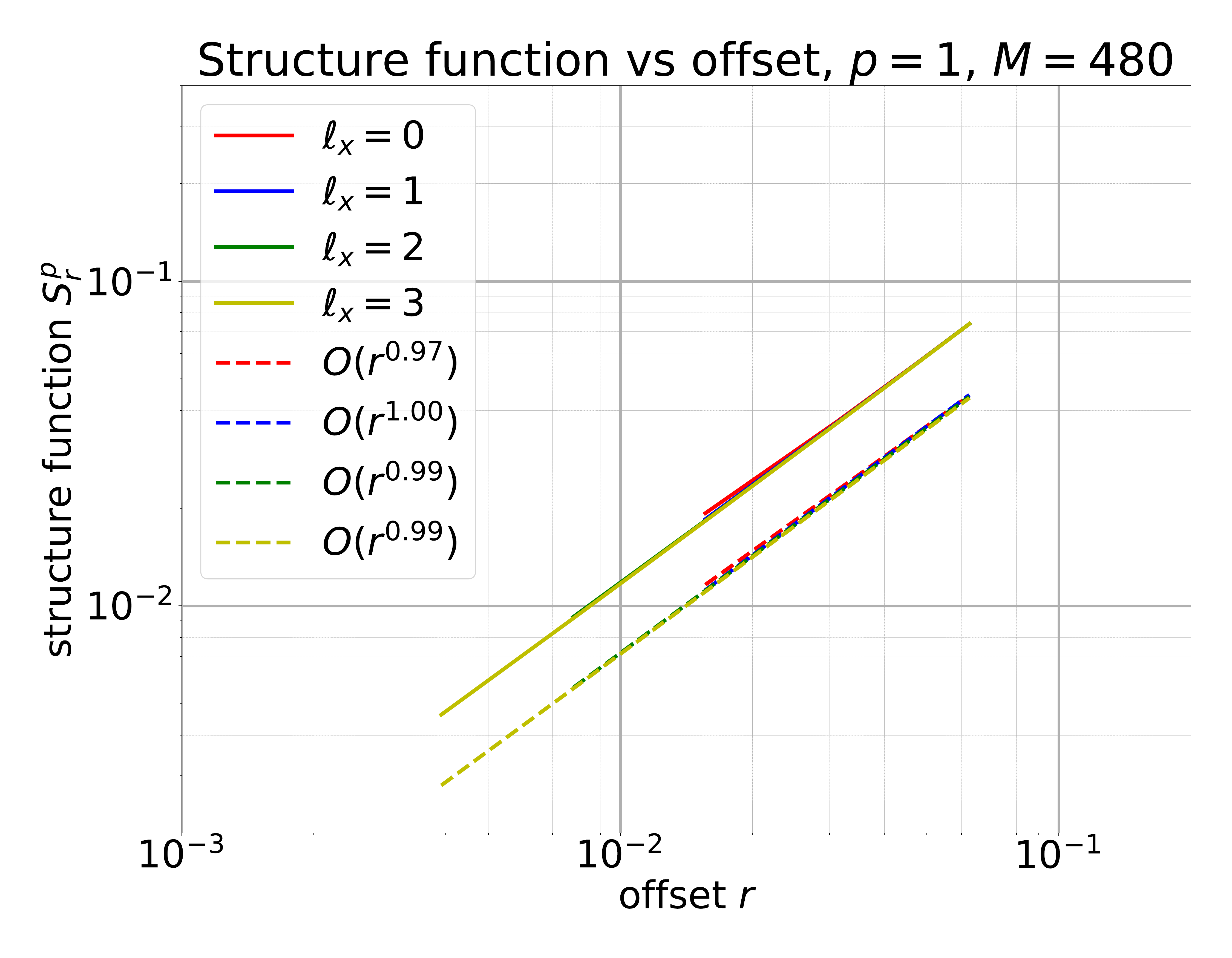}
  \end{subfigure}
  \begin{subfigure}[t]{.5\textwidth}
    \centering\includegraphics[width = \textwidth]{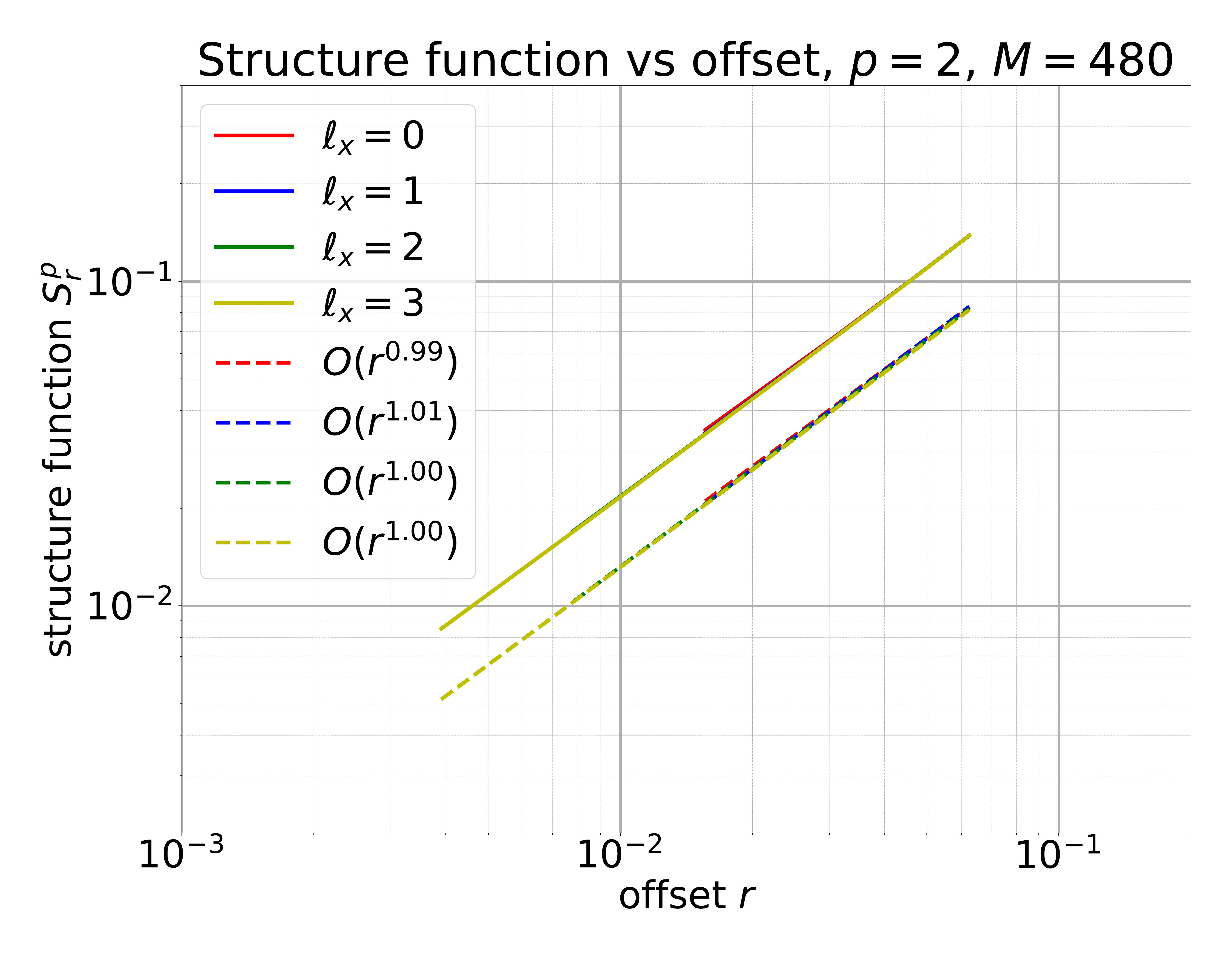}
  \end{subfigure}
  \begin{subfigure}[t]{.5\textwidth}
    \centering\includegraphics[width = \textwidth]{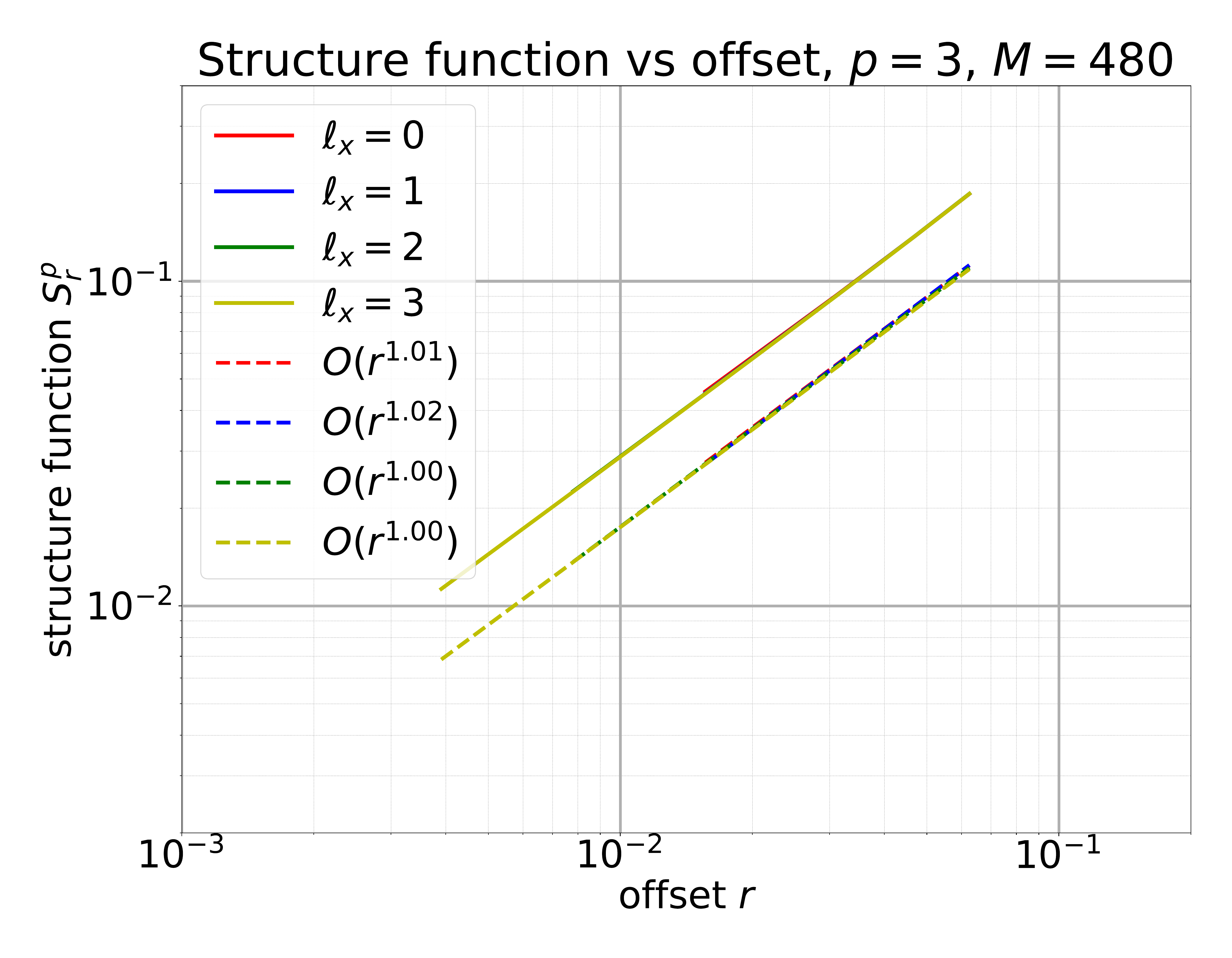}
  \end{subfigure}
 \caption[Channel flow, $Re=3200$, convergence of structure function]{\small {Channel flow with $Re=3200$, 
  as described in {\S}\ref{ss:uqIncompNSTest3};
  Approximation of the structure function of the velocity field with
  $480$ samples, at time $T=0.8$.
 }}
 \label{fig:uqIncompNSTest3Re3200Scube}
\end{figure}

In this experiment, we study the convergence of the mean, variance and Wasserstein distances with respect to the number of samples for a fixed mesh level, and also when both the number of samples and mesh level are increased at the same time. 
As the results for different mesh levels were observed to be similar, we omit many of them for brevity and show the results only for the mesh level $\ellx = 3$ in figures \ref{fig:uqIncompNSTest3VisVyMean} -- \ref{fig:uqIncompNSTest3Re3200Scube}.
The visualizations of the mean of the horizontal velocity are very similar irrespective of the number of samples and the mesh levels, they show a parabolic profile with slight variations. Thus, they have not been shown separately.

The observations are as follows:
\begin{itemize}[nolistsep]
 \item Figure~\ref{fig:uqIncompNSTest3VisVyMean} suggest the formation of vortices in the flow because of interaction with the domain boundary. We can observe some differences in the flow features between Reynolds number $3200$ and $1600$, the flow is more chaotic for $Re=3200$ that leads to higher variance as given in Table~\ref{tab:uqIncompNSTest3CompareRe}.
 \item From figures~\ref{fig:uqIncompNSTest3VisVyMean} and  \ref{fig:uqIncompNSTest3VisVyVarRescaled}, the mean and variance appear to converge as the number of samples is increased, irrespective of the Reynolds number. 
 This observation is further supported by the Cauchy convergence of these statistics in figures~\ref{fig:uqIncompNSTest3Re1600Convg} and \ref{fig:uqIncompNSTest3Re3200Convg}.
 \item We measure Wasserstein distances $W^{1}$ and $W^{2}$ (see \eqref{eq:WassersteinDistW12}) between two ensembles $\bb{U}^{h,M}_{T}$ and $\bb{U}^{h,2M}_{T}$. From figures~\ref{fig:uqIncompNSTest3Re1600Wd} and \ref{fig:uqIncompNSTest3Re3200Wd}, we observe that these distances decrease as we increase the number of samples, which indicates convergence.
 \item When the mesh level and the number of samples are increased simultaneously, the statistics and Wasserstein distances still converge, see tables~\ref{tab:uqIncompNSTest3Re1600} and \ref{tab:uqIncompNSTest3Re3200}.
 \item From figures~\ref{fig:uqIncompNSTest3Re1600Scube} and \ref{fig:uqIncompNSTest3Re3200Scube}, we observe that structure functions converge with a rate $1$ (approx.), which is expected as the solution remains smooth (cf. Remark~\ref{rem:structureFnScaling}).
\end{itemize}

\begin{table}[!htb]
\begin{center}
\begin{tabular}{c| c | c}
\multicolumn{1}{c|}{}
& \multicolumn{2}{c}{$Re$}\\
\hline
  & $1600$ & $3200$ \\
\hline
 Mean     & 9.4839E-01 & 9.4839E-01 \\
 Variance & 5.0132E-04 & 5.2386E-04
\end{tabular}
\caption[Channel flow, mean and variance $Re=1600$ vs $Re=3200$]{\small {Channel flow, 
  as described in {\S}\ref{ss:uqIncompNSTest3};
  Mean and variance in $\ltwo{\Dx}$ norm for mesh level $\ellx = 3$ and $M=480$ samples, at time $T=0.8$.
 }}
\label{tab:uqIncompNSTest3CompareRe}
\end{center}
\end{table}

\begin{table}[!htb]
\begin{center}
\begin{tabular}{c| c | c| c c | c c}
\multicolumn{3}{c|}{}
& \multicolumn{2}{c|}{$\abs{\bu}$}
& \multicolumn{2}{c}{$\bu$}\\
\hline
$\ellx$, $M$ & Mean & Variance & $W^{1}$ & $W^{2}$ & $W^{1}$ & $W^{2}$\\
\hline
 $0$, $60$  & - & - & - & - & - & -\\
 $1$, $120$ & 6.2959E-03 & 1.0938E-04 & 4.7060E-03 & 8.9483E-03 & 4.7648E-03 & 9.0448E-03\\
 $2$, $240$ & 1.6118E-03 & 1.1144E-04 & 2.5590E-03 & 5.6916E-03 & 2.6279E-03 & 5.8191E-03\\
 $3$, $480$ & 6.4555E-04 & 7.1313E-05 & 1.6540E-03 & 4.1550E-03 & 1.7321E-03 & 4.3106E-03
\end{tabular}
\caption[Channel flow, $Re=1600$, convergence wrt mesh resolution]{\small {Channel flow with $Re=1600$, 
  as described in {\S}\ref{ss:uqIncompNSTest3};
  Cauchy error of mean and variance in $\ltwo{\Dx}$ norm and Wasserstein distances measured between ensembles, at time $T=0.8$.
 }}
\label{tab:uqIncompNSTest3Re1600}
\end{center}
\end{table}

\begin{table}[!htb]
\begin{center}
\begin{tabular}{c| c | c| c c | c c}
\multicolumn{3}{c|}{}
& \multicolumn{2}{c|}{$\abs{\bu}$}
& \multicolumn{2}{c}{$\bu$}\\
\hline
$\ellx$, $M$ & Mean & Variance & $W^{1}$ & $W^{2}$ & $W^{1}$ & $W^{2}$\\
\hline
 $0$, $60$  & - & - & - & - & - & -\\
 $1$, $120$ & 6.3397E-03 & 1.1594E-04 & 4.8123E-03 & 9.1556E-03 & 4.8712E-03 & 9.2548E-03\\
 $2$, $240$ & 1.6178E-03 & 1.1930E-04 & 2.5880E-03 & 5.8111E-03 & 2.6603E-03 & 5.9408E-03\\
 $3$, $480$ & 6.5222E-04 & 7.4728E-05 & 1.7059E-03 & 4.2674E-03 & 1.7859E-03 & 4.4257E-03
\end{tabular}
\caption[Channel flow, $Re=3200$, convergence wrt mesh resolution]{\small {Channel flow with $Re=3200$, 
  as described in {\S}\ref{ss:uqIncompNSTest3};
  Cauchy error of mean and variance in $\ltwo{\Dx}$ norm and Wasserstein distances measured between ensembles, at time $T=0.8$.
 }}
\label{tab:uqIncompNSTest3Re3200}
\end{center}
\end{table}

%
\section{Conclusions} \label{s:uqIncompNSConclusions}
%
We conclude this chapter by summarizing the main results in the following:
\begin{itemize}[nolistsep]
 \item We computed the statistical solutions of two-dimensional incompressible Navier-Stokes equations for problems with non-periodic boundary conditions, in particular, the lid-driven cavity problem with a moving top boundary and the channel flow problem with an outflow boundary. Through our numerical experiments, we have gathered ample evidence to conclude that the approximate statistical solutions converge, where we compute these approximations using our MC-FEM solver.
 \item We have developed a novel algorithm to numerically approximate structure functions on unstructured meshes for two-dimensional rectangular domains. We can extend this algorithm to more general two-dimensional domains with some modifications, see Remark~\ref{rem:computeStructureFnsGeneralDomain}. Moreover, the algorithm can be extended for three-dimensional rectangular cuboid domains, albeit an efficient parallel version would be necessary in this case.
\end{itemize}

Numerical computation of an individual solution of the NSE for high Reynolds numbers flows with non-periodic boundary using FEM is expensive because it requires meshes with a high resolution near the boundary. An alternative to the H(div) scheme used in the present work is the more computationally efficient H(div)-HDG scheme with static condensation \cite{LLSc2018, LSc2016}. Then, combined with MC or QMC sampling, one can simulate more challenging problems, such as flow past a cylinder for high Reynolds numbers where the vortex-shedding phenomenon is observed \cite[Chapter~9]{PMi2015}. The flow past a cylinder is an interesting problem to investigate in the vanishing viscosity limit.

\clearpage
\section*{Acknowledgements} 
The author acknowledges support from the project ModCompShock, funded by the European Union's Horizon 2020 research and innovation programme under the Marie Sklodowska-Curie grant agreement number 642768, for the research presented here. The author is grateful to Prof. Dr. Siddhartha Mishra and Prof. Dr. Christoph Schwab, both associated with Seminar for Applied Mathematics (ETH Zurich), for the insightful discussions and their guidance.

\bibliographystyle{plain}

\input{references}
\end{document}

%% file: macros.tex
\numberwithin{equation}{section}		
\numberwithin{figure}{section}			
\numberwithin{table}{section}			

\usepackage[greek.polutoniko,english]{babel}

\theoremstyle{plain}
\newtheorem{theorem}{\sffamily Theorem}[section]

\newtheorem{lemma}[theorem]{\sffamily Lemma}

\newtheorem{remark}[theorem]{\sffamily Remark}
\newtheorem{definition}[theorem]{\sffamily Definition}

\newcommand{\der}[2]{\partial_{#2} #1}
\newcommand{\dder}[2]{\partial_{#2#2} #1}
\newcommand{\abs}[1]{\left\lvert #1\right\rvert}
\newcommand{\norm}[1]{\left\lVert#1\right\rVert}

\renewcommand{\div}{\operatorname{div}}

\newcommand{\mrm}[1]{\mathrm{#1}}

\newcommand{\R}{\mathbb{R}}

\newcommand{\N}{\mathbb{N}}

\newcommand{\E}{\mathbb{E}}
\newcommand{\F}{\mathbb{F}}
\newcommand{\Nz}{\mathbb{N}_{0}}
\newcommand{\bbP}{\mathbb{P}}

\newcommand{\bb}[1]{\mathbf{#1}}

\newcommand{\bx}{\mathbf{x}}
\newcommand{\by}{\mathbf{y}}

\newcommand{\bu}{\mathbf{u}}
\newcommand{\bv}{\mathbf{v}}
\newcommand{\bw}{\mathbf{w}}

\newcommand{\btau}{\boldsymbol{\tau}}

\newcommand{\bn}{\mathbf{n}}

\newcommand{\calI}{\mathcal{I}}
\newcommand{\calS}{\mathcal{S}}

\newcommand{\calM}{\mathcal{M}}

\newcommand{\calD}{\mathcal{D}}
\newcommand{\calA}{\mathcal{A}}
\newcommand{\calF}{\mathcal{F}}
\newcommand{\calC}{\mathcal{C}}

\newcommand{\calT}{\mathcal{T}}

\newcommand{\calB}{\mathcal{B}}
\newcommand{\calP}{\mathcal{P}}

\newcommand{\frakD}{\mathfrak{D}}

\newcommand{\Dx}{\calD}
\newcommand{\Dt}{\calI}
\newcommand{\Dxb}{\partial\Dx}
\newcommand{\DxbD}{{\Dxb}_{\mathrm{D}}}

\newcommand{\tn}{t_{n}}
\newcommand{\Dxbo}{{\Dxb}_{\mathrm{out}}}

\newcommand{\Dxm}[1]{\calT_{#1}}

\newcommand{\Dxmh}{\Dxm{h}}
\newcommand{\hx}{h}
\newcommand{\Kx}{K}
\newcommand{\Kxb}{\partial\Kx}
\newcommand{\Kb}{\partial K}
\newcommand{\hKx}{h_{\Kx}}

\newcommand{\hFx}{h_{\Fx}}

\newcommand{\Fx}{F}
\newcommand{\Dxmf}[1]{\calF_{#1}}
\newcommand{\Dxmfh}{\Dxmf{\hx}}
\newcommand{\Dxmfhb}{\Dxmf{\hx}^{\mrm{bdr}}}
\newcommand{\DxmfhbD}{\Dxmf{\hx}^{\mrm{bdr},\mrm{D}}}
\newcommand{\Dxmfhbo}{\Dxmf{\hx}^{\mrm{bdr},\mrm{out}}}
\newcommand{\Dxmfhi}{\Dxmf{\hx}^{\mrm{int}}}

\newcommand{\avg}[1]{\{\!\{#1\}\!\}}
\newcommand{\jump}[1]{[\![#1]\!]}

\newcommand{\jjump}[1]{[\![\![#1]\!]\!]}

\newcommand{\hdiv}[1]{H(\div;#1)}

\newcommand{\hdivDx}{\hdiv{\Dx}}

\newcommand{\rt}[2]{RT_{#2}(#1)}
\newcommand{\globRt}[1]{{RT_{#1}}(\Dxmh)}

\newcommand{\linf}[1]{L^{\infty}(#1)}

\newcommand{\ltwo}[1]{L^2(#1)}
\newcommand{\ltwoz}[1]{L^2_0(#1)}
\newcommand{\hone}[1]{H^1(#1)}

\newcommand{\lfSource}[1]{\ell_h(#1)}
\newcommand{\bfMasse}{\calM_h}
\newcommand{\bfMass}[2]{\calM_h(#1,#2)}
\newcommand{\bfDive}{\calB_h}
\newcommand{\bfDiv}[2]{\calB_h(#1,#2)}

\newcommand{\bfConvUe}{\calC_h^{up}}
\newcommand{\bfConvU}[3]{\calC_h^{up}(#1;#2,#3)}

\newcommand{\bfDiffSipe}{\calA_h^{sip}}
\newcommand{\bfDiffSip}[2]{\calA_h^{sip}(#1,#2)}

\newcommand{\nDx}{\bn_{\Dx}}

\newcommand{\nxKx}{\bn_{\Kx}}

\newcommand{\half}{\frac{1}{2}}

\newcommand{\ellx}{\ell_{x}}

\newcommand{\Vk}{V_{k}}
\newcommand{\Vzk}{V_{0,k}}
\newcommand{\Wk}{Q_{k}}
\newcommand{\Wzk}{Q_{0,k}}

\newcommand{\penParam}{\sigma}

\newcommand{\buh}{\bu_{h}}
\newcommand{\bvh}{\bv_{h}}

\renewcommand{\vec}[1]{\underline{#1}}
\newcommand{\mat}[1]{\underline{\bb{#1}}}

%% file: statSolsNSE.bbl
\begin{thebibliography}{10}

\bibitem{mfem}
{MFEM}: Modular finite element methods library.
\newblock \url{mfem.org}.

\bibitem{BCGl1989}
John~B. Bell, Phillip Colella, and Harland~M. Glaz.
\newblock A second-order projection method for the incompressible
  {N}avier-{S}tokes equations.
\newblock {\em J. Comput. Phys.}, 85(2):257--283, 1989.

\bibitem{BGLi2005}
Michele Benzi, Gene~H. Golub, and J\"{o}rg Liesen.
\newblock Numerical solution of saddle point problems.
\newblock {\em Acta Numer.}, 14:1--137, 2005.

\bibitem{BCBBr2018}
Crist\'{o}bal Bertoglio, Alfonso Caiazzo, Yuri Bazilevs, Malte Braack, Mahdi
  Esmaily, Volker Gravemeier, Alison~L. Marsden, Olivier Pironneau, Irene~E.
  Vignon-Clementel, and Wolfgang~A. Wall.
\newblock Benchmark problems for numerical treatment of backflow at open
  boundaries.
\newblock {\em Int. J. Numer. Methods Biomed. Eng.}, 34(2):e2918, 34, 2018.

\bibitem{BLMSc2013}
Hester Bijl, Didier Lucor, Siddhartha Mishra, and Christoph Schwab, editors.
\newblock {\em Uncertainty quantification in computational fluid dynamics},
  volume~92 of {\em Lecture Notes in Computational Science and Engineering}.
\newblock Springer, Heidelberg, 2013.

\bibitem{BBFo2013}
Daniele Boffi, Franco Brezzi, and Michel Fortin.
\newblock {\em Mixed finite element methods and applications}, volume~44 of
  {\em Springer Series in Computational Mathematics}.
\newblock Springer, Heidelberg, 2013.

\bibitem{GriSPy}
Martin Chalela, Emanuel Sillero, Luis Pereyra, Mario~Alejandro García, Juan~B.
  Cabral, Marcelo Lares, and Manuel Merchán.
\newblock Grispy: A python package for fixed-radius nearest neighbors search,
  2020.

\bibitem{Ch1967}
Alexandre~Joel Chorin.
\newblock The numerical solution of the {N}avier-{S}tokes equations for an
  incompressible fluid.
\newblock {\em Bull. Amer. Math. Soc.}, 73:928--931, 1967.

\bibitem{Ch1969}
Alexandre~Joel Chorin.
\newblock On the convergence of discrete approximations to the
  {N}avier-{S}tokes equations.
\newblock {\em Math. Comp.}, 23:341--353, 1969.

\bibitem{CKSc2007}
Bernardo Cockburn, Guido Kanschat, and Dominik Sch\"{o}tzau.
\newblock A note on discontinuous {G}alerkin divergence-free solutions of the
  {N}avier-{S}tokes equations.
\newblock {\em J. Sci. Comput.}, 31(1-2):61--73, 2007.

\bibitem{DEr2012}
Daniele~Antonio Di~Pietro and Alexandre Ern.
\newblock {\em Mathematical aspects of discontinuous {G}alerkin methods},
  volume~69 of {\em Math\'{e}matiques \& Applications (Berlin) [Mathematics \&
  Applications]}.
\newblock Springer, Heidelberg, 2012.

\bibitem{Do2015}
S.~Dong.
\newblock A convective-like energy-stable open boundary condition for
  simulations of incompressible flows.
\newblock {\em J. Comput. Phys.}, 302:300--328, 2015.

\bibitem{DKCh2014}
S.~Dong, G.~E. Karniadakis, and C.~Chryssostomidis.
\newblock A robust and accurate outflow boundary condition for incompressible
  flow simulations on severely-truncated unbounded domains.
\newblock {\em J. Comput. Phys.}, 261:83--105, 2014.

\bibitem{FKLLSc2019}
Niklas Fehn, Martin Kronbichler, Christoph Lehrenfeld, Gert Lube, and
  Philipp~W. Schroeder.
\newblock High-order dg solvers for underresolved turbulent incompressible
  flows: A comparison of l2 and h(div) methods.
\newblock {\em International Journal for Numerical Methods in Fluids},
  91(11):533--556, 2019.

\bibitem{FLMi2017}
U.~S. Fjordholm, S.~Lanthaler, and S.~Mishra.
\newblock Statistical solutions of hyperbolic conservation laws: foundations.
\newblock {\em Arch. Ration. Mech. Anal.}, 226(2):809--849, 2017.

\bibitem{pot}
R{'e}mi Flamary and Nicolas Courty.
\newblock Pot python optimal transport library, 2017.

\bibitem{Fo1972}
C.~Foia\c{s}.
\newblock Statistical study of {N}avier-{S}tokes equations. {I}, {II}.
\newblock {\em Rend. Sem. Mat. Univ. Padova}, 48:219--348 (1973); ibid. 49
  (1973), 9--123, 1972.

\bibitem{FMRTe2001}
C.~Foia\c{s}, O.~Manley, R.~Rosa, and R.~Temam.
\newblock {\em Navier-{S}tokes equations and turbulence}, volume~83 of {\em
  Encyclopedia of Mathematics and its Applications}.
\newblock Cambridge University Press, Cambridge, 2001.

\bibitem{FPr1976}
C.~Foia\c{s} and G.~Prodi.
\newblock Sur les solutions statistiques des \'{e}quations de
  {N}avier-{S}tokes.
\newblock {\em Ann. Mat. Pura Appl. (4)}, 111:307--330, 1976.

\bibitem{Fr2006}
W~Franklin.
\newblock Nearest point query on 184,088,599 points in e3 with a uniform grid.
  p.
\newblock 02 2006.

\bibitem{Fr1995}
Uriel Frisch.
\newblock {\em Turbulence}.
\newblock Cambridge University Press, Cambridge, 1995.
\newblock The legacy of A. N. Kolmogorov.

\bibitem{Ga2014}
Gabriel~N. Gatica.
\newblock {\em A simple introduction to the mixed finite element method}.
\newblock SpringerBriefs in Mathematics. Springer, Cham, 2014.
\newblock Theory and applications.

\bibitem{Gmsh}
Christophe Geuzaine and Jean-Fran\c{c}ois Remacle.
\newblock Gmsh: {A} 3-{D} finite element mesh generator with built-in pre- and
  post-processing facilities.
\newblock {\em Internat. J. Numer. Methods Engrg.}, 79(11):1309--1331, 2009.

\bibitem{GRa1986}
Vivette Girault and Pierre-Arnaud Raviart.
\newblock {\em Finite element methods for {N}avier-{S}tokes equations},
  volume~5 of {\em Springer Series in Computational Mathematics}.
\newblock Springer-Verlag, Berlin, 1986.
\newblock Theory and algorithms.

\bibitem{GLSk2014}
W.~{Gropp}, E.~{Lusk}, and A.~{Skjellum}.
\newblock {\em Using MPI in Simple Programs}, pages 23--68.
\newblock 2014.

\bibitem{GSSe2017}
Johnny Guzm\'{a}n, Chi-Wang Shu, and Fil\'{a}nder~A. Sequeira.
\newblock {$\rm H(div)$} conforming and {DG} methods for incompressible
  {E}uler's equations.
\newblock {\em IMA J. Numer. Anal.}, 37(4):1733--1771, 2017.

\bibitem{Ho1951}
Eberhard Hopf.
\newblock \"{U}ber die {A}nfangswertaufgabe f\"{u}r die hydrodynamischen
  {G}rundgleichungen.
\newblock {\em Math. Nachr.}, 4:213--231, 1951.

\bibitem{JLMNRe2017}
Volker John, Alexander Linke, Christian Merdon, Michael Neilan, and Leo~G.
  Rebholz.
\newblock On the divergence constraint in mixed finite element methods for
  incompressible flows.
\newblock {\em SIAM Rev.}, 59(3):492--544, 2017.

\bibitem{Ko1941a}
A.~Kolmogoroff.
\newblock The local structure of turbulence in incompressible viscous fluid for
  very large {R}eynold's numbers.
\newblock {\em C. R. (Doklady) Acad. Sci. URSS (N.S.)}, 30:301--305, 1941.

\bibitem{Ko1941b}
A.~N. Kolmogoroff.
\newblock On degeneration of isotropic turbulence in an incompressible viscous
  liquid.
\newblock {\em C. R. (Doklady) Acad. Sci. URSS (N. S.)}, 31:538--540, 1941.

\bibitem{KLa1966}
A.~Krzhivitski and O.~A. Ladyzhenskaya.
\newblock A grid method for the {N}avier-{S}tokes equations.
\newblock {\em Soviet Physics Dokl.}, 11:212--213, 1966.

\bibitem{LMPa2021b}
S.~Lanthaler, S.~Mishra, and C.~Par{\'e}s-Pulido.
\newblock On the conservation of energy in two-dimensional incompressible
  flows.
\newblock {\em Nonlinearity}, 34(2):1084--1135, 2021.

\bibitem{LMPa2021a}
S.~Lanthaler, S.~Mishra, and C.~Par{\'e}s-Pulido.
\newblock Statistical solutions of the incompressible {E}uler equations.
\newblock {\em Math. Models Methods Appl. Sci.}, 31(2):223--292, 2021.

\bibitem{LLSc2018}
Philip~L. Lederer, Christoph Lehrenfeld, and Joachim Sch{\"{o}}berl.
\newblock Hybrid discontinuous {G}alerkin methods with relaxed {$H({\rm
  div})$}-conformity for incompressible flows. {P}art {I}.
\newblock {\em SIAM J. Numer. Anal.}, 56(4):2070--2094, 2018.

\bibitem{LSc2016}
Christoph Lehrenfeld and Joachim Sch\"{o}berl.
\newblock High order exactly divergence-free hybrid discontinuous {G}alerkin
  methods for unsteady incompressible flows.
\newblock {\em Comput. Methods Appl. Mech. Engrg.}, 307:339--361, 2016.

\bibitem{Le2018}
Filippo Leonardi.
\newblock {\em Numerical methods for ensemble based solutions to incompressible
  flow equations}.
\newblock PhD thesis, ETH Z{\"u}rich, 2018.
\newblock Dis\_no 25171, Prof. Dr. Siddhartha Mishra.

\bibitem{LMSc2016}
Filippo Leonardi, Siddhartha Mishra, and Christoph Schwab.
\newblock Numerical approximation of statistical solutions of planar,
  incompressible flows.
\newblock {\em Math. Models Methods Appl. Sci.}, 26(13):2471--2523, 2016.

\bibitem{Le1933}
Jean Leray.
\newblock {\em \'{E}tude de diverses \'{e}quations int\'{e}grales non
  lin\'{e}aires et de quelques probl\`emes que pose l'hydrodynamique}.
\newblock NUMDAM, [place of publication not identified], 1933.

\bibitem{Le1934}
Jean Leray.
\newblock Sur le mouvement d'un liquide visqueux emplissant l'espace.
\newblock {\em Acta Math.}, 63(1):193--248, 1934.

\bibitem{LTa1997}
Doron Levy and Eitan Tadmor.
\newblock Non-oscillatory central schemes for the incompressible {$2$}-{D}
  {E}uler equations.
\newblock {\em Math. Res. Lett.}, 4(2-3):321--340, 1997.

\bibitem{Ly2020}
Kjetil Lye.
\newblock {\em Computation of statistical solutions of hyperbolic systems of
  conservation laws}.
\newblock PhD thesis, ETH Z{\"u}rich, 2020.
\newblock Dis\_no 26728, Prof. Dr. Siddhartha Mishra.

\bibitem{MPRo1990}
Y.~Maday, Anthony~T. Patera, and Einar~M. R\o~nquist.
\newblock An operator-integration-factor splitting method for time-dependent
  problems: application to incompressible fluid flow.
\newblock {\em J. Sci. Comput.}, 5(4):263--292, 1990.

\bibitem{PMi2015}
John W.~Mitchell Philip J.~Pritchard.
\newblock {\em Fox and McDonald's introduction to fluid mechanics}.
\newblock John Wiley \& Sons, New York, 8 edition, 2015.

\bibitem{SSc1986}
Youcef Saad and Martin~H. Schultz.
\newblock G{MRES}: a generalized minimal residual algorithm for solving
  nonsymmetric linear systems.
\newblock {\em SIAM J. Sci. Statist. Comput.}, 7(3):856--869, 1986.

\bibitem{Sc2005}
Ren\'{e}~L. Schilling.
\newblock {\em Measures, integrals and martingales}.
\newblock Cambridge University Press, New York, 2005.

\bibitem{SJLLLSc2019}
Philipp~W. Schroeder, Volker John, Philip~L. Lederer, Christoph Lehrenfeld,
  Gert Lube, and Joachim Sch\"{o}berl.
\newblock On reference solutions and the sensitivity of the 2d kelvin-helmholtz
  instability problem.
\newblock {\em Comput. Math. Appl.}, 77(4):1010--1028, 2019.

\bibitem{SLe2018}
Philipp~W. Schroeder, Christoph Lehrenfeld, Alexander Linke, and Gert Lube.
\newblock Towards computable flows and robust estimates for inf-sup stable
  {FEM} applied to the time-dependent incompressible {N}avier-{S}tokes
  equations.
\newblock {\em SeMA J.}, 75(4):629--653, 2018.

\bibitem{SLu2018}
Philipp~W. Schroeder and Gert Lube.
\newblock Divergence-free {$H({\rm div})$}-{FEM} for time-dependent
  incompressible flows with applications to high {R}eynolds number vortex
  dynamics.
\newblock {\em J. Sci. Comput.}, 75(2):830--858, 2018.

\bibitem{Su2015}
T.~J. Sullivan.
\newblock {\em Introduction to uncertainty quantification}, volume~63 of {\em
  Texts in Applied Mathematics}.
\newblock Springer, Cham, 2015.

\bibitem{VWe1996}
Aad~W. van~der Vaart and Jon~A. Wellner.
\newblock {\em Weak convergence and empirical processes}.
\newblock Springer Series in Statistics. Springer-Verlag, New York, 1996.
\newblock With applications to statistics.

\end{thebibliography}
